\documentclass[opre]{informs3}

\OneAndAHalfSpacedXI 

\usepackage{endnotes}
\let\footnote=\endnote

\usepackage{algorithm}
\usepackage{algorithmic}
\usepackage{booktabs}
\usepackage{subfigure}
\usepackage{pdflscape}
\usepackage[section]{placeins}
\usepackage{enumerate}
\usepackage{bbm}
\def\qed{$\square$}

\usepackage{natbib}
 \bibpunct[, ]{(}{)}{,}{a}{}{,}%
 \def\bibfont{\small}%
\TheoremsNumberedThrough     
\ECRepeatTheorems

\EquationsNumberedThrough    

\MANUSCRIPTNO{OPRE-2015-10-573.R2}

\begin{document}


\RUNAUTHOR{Chan, Shen and Siddiq}

\RUNTITLE{Robust Defibrillator Deployment}

\TITLE{Robust Defibrillator Deployment Under  
Cardiac Arrest Location Uncertainty via Row-and-Column Generation}

\ARTICLEAUTHORS{%
\AUTHOR{Timothy C. Y. Chan}
\AFF{Department of Mechanical and Industrial Engineering, University of Toronto, Toronto, Canada, \\ \EMAIL{tcychan@mie.utoronto.ca}.} 

\AUTHOR{Zuo-Jun Max Shen}
\AFF{Department of Industrial Engineering and Operations Research, University of California, Berkeley, USA, \EMAIL{maxshen@berkeley.edu}.}

\AUTHOR{Auyon Siddiq}
\AFF{Department of Industrial Engineering and Operations Research, University of California, Berkeley, USA, \EMAIL{auyon.siddiq@berkeley.edu}.}
} 

\ABSTRACT{%
Sudden cardiac arrest is a significant public health concern. Successful treatment of cardiac arrest is extremely time sensitive, and use of an automated external defibrillator (AED) where possible significantly increases the probability of survival. Placement of AEDs in public locations can improve survival by enabling bystanders to treat victims of cardiac arrest prior to the arrival of emergency medical responders, thus shortening the time between collapse and treatment. However, since the exact locations of future cardiac arrests cannot be known {\it a priori}, AEDs must be placed strategically in public locations to ensure their accessibility in the event of an out-of-hospital cardiac arrest emergency. In this paper, we propose a data-driven optimization model for deploying AEDs in public spaces while accounting for uncertainty in future cardiac arrest locations. Our approach involves discretizing a continuous service area into a large set of scenarios, where the probability of cardiac arrest at each location is itself uncertain. We model uncertainty in the spatial risk of cardiac arrest using a polyhedral uncertainty set that we calibrate using historical cardiac arrest data. We propose a solution technique based on row-and-column generation that exploits the structure of the uncertainty set, allowing the algorithm to scale gracefully with the total number of scenarios. Using real cardiac arrest data from the City of Toronto, we conduct an extensive numerical study on AED deployment public locations. We find that hedging against cardiac arrest location uncertainty can produce AED deployments that outperform a intuitive sample average approximation by 9 to 15\%, and cuts the performance gap with respect to an ex-post model by half. Our findings suggest that accounting for cardiac arrest location uncertainty can lead to improved accessibility of AEDs during cardiac arrest emergencies and the potential for improved survival outcomes.

}%


\KEYWORDS{robust optimization; row-and-column generation; healthcare operations; facility location.}

\maketitle

%


\vspace{-10mm}
\section{Introduction}\label{sec:introduction}
Sudden cardiac arrest is a leading cause of death in North America, and is responsible for over 400,000 deaths each year \citep{HSF,mozaffarian2016executive}. Chances of survival have been estimated to decrease by as much as 10 percent with each minute of delay in treatment \citep{larsen1993}, making treatment of cardiac arrest extremely time sensitive. Currently, the probability of survival from an out-of-hospital cardiac arrest is low, with only 5-10\% of out-of-hospital cardiac arrest victims surviving to hospital discharge \citep{weisfeldt2010}.

Treatment of cardiac arrest involves cardiopulmonary resuscitation (CPR) and electrical shocks -- known as defibrillation -- from an  automated external defibrillator (AED). AEDs are portable devices that can automatically assess heart rhythms and perform defibrillation if necessary, making their use by non-professionals (lay responders) in a cardiac arrest emergency viable \citep{gundry1999comparison}. Treatment with a defibrillator following cardiac arrest is vital: the probability of survival from cardiac arrest has been estimated to be as high as 40-70\% if defibrillation is administered within three minutes of the victim's collapse \citep{valenzuela2000,page2000use,caffrey2002public}. As a consequence, there is growing interest in the development of public access defibrillation (PAD) programs, which are organized efforts to place defibrillators in public areas, such as malls, coffee shops, restaurants, and subway stations, with the hopes of improving easy access to defibrillators in the event of a cardiac arrest emergency. With appropriate positioning, public AEDs can enable lay responders to treat cardiac arrest victims while awaiting the arrival of professional emergency medical responders. Indeed, PAD programs have been shown to decrease the time to treatment and improve survival outcomes \citep{public2004public}, and it is estimated that their widespread implementation may save between 2,000 and 4,000 lives annually in the United States \citep{hazinski2005lay}. Thus, the strategic positioning of public AEDs has a crucial role to play in strengthening the overall response to cardiac arrest.

Despite evidence that supports improved chances of survival with PAD programs, AEDs are used infrequently. For example, only 8\% of public location cardiac arrests in Toronto, Canada involve the use of an AED  by a bystander \citep{sun2016overcoming}. Although the reasons for this low usage rate are multi-faceted, it is clear that AEDs that are poorly positioned may be located too far from the locations of cardiac arrest emergencies for them to be accessed and used by lay responders. Recently, several mobile phone applications have been developed that aim to improve AED usage during public location cardiac arrests. These applications notify volunteer lay responders of the location of nearby AEDs so that they can be retrieved and used quickly (e.g. \cite{pulsepoint}). In order to be effective, these emerging technologies -- and public access defibrillation in general -- rely critically on AEDs being placed strategically in public locations. The question of where to place public access AEDs remains relatively open in the medical literature \citep{portner2004out, folke2009}.


In this paper, we present a data-driven optimization model for deploying AEDs in public locations. We consider a setting where cardiac arrests occur continuously over a bounded service area. Our model addresses a key challenge in the AED deployment problem, which is the uncertainty in the locations of future cardiac arrests during the deployment of AEDs. Since the chances of survival from cardiac arrest diminish rapidly with each passing minute, it is critical that witnesses of cardiac arrest are able to quickly retrieve and use nearby AEDs when needed. Due to this extreme time sensitivity, protecting against cardiac arrest location uncertainty can help improve survival outcomes by shortening the distance (and thus the response time) to the nearest AED during a cardiac arrest emergency.

Although cardiac arrests are not uncommon within the general population, they are extremely rare events for any {\it particular} building or geographic location. Thus, any estimate of the risk of cardiac arrest at a specific location is likely to be subject to considerable uncertainty. To account for this uncertainty while still leveraging historical data, we take a distributionally robust optimization approach in which service is optimized for the worst-case spatial distribution of cardiac arrests. Our approach involves discretizing the continuous service area into a large set of scenarios, where the probability vector that describes the probability of the next cardiac arrest arriving at each scenario is only known to reside within a polyhedral uncertainty set. We construct the uncertainty set by representing the service area as the union of a set of {\it uncertainty regions}, where only the aggregate probability of the next cardiac arrest arriving in each uncertainty region (i.e., group of scenarios) is known. An important consequence of modeling uncertainty in this manner is that it induces sparsity in the worst-case distribution (with respect to the number of scenarios), which we show can be exploited through an efficient row-and-column algorithm to obtain optimal solutions to the robust problem.

For a given AED deployment and a set of cardiac arrest locations, one can construct a {\it distance distribution} which describes the histogram of distances between the cardiac arrests and the nearest AED in the given deployment. Inspired by the use of quantile-based targets in emergency medical services \citep{pell2001,pons2002}, we use a conditional value-at-risk (CVaR) objective function to directly optimize the tail of the worst-case distance distribution, which allows us to mitigate the risk of unacceptably long distances between cardiac arrest locations and their nearest AEDs. An additional benefit of a CVaR-based objective function is that it permits more flexibility with respect to selecting a performance metric to optimize, e.g., we can minimize the mean or tail of the worst-case distribution by using an appropriate parameterization of the model. 

We summarize our contributions as follows.

\begin{enumerate}
\item
We present a robust optimization model for AED deployment that accounts for uncertainty in the spatial distribution of cardiac arrests. Our model permits cardiac arrests to arrive anywhere within a continuous service area, which we approximate through a fine discretization. We use historical cardiac arrest data to construct a distributional uncertainty set that consists of a set of generally-shaped uncertainty regions, where only the probability of the next cardiac arrest event occurring within each uncertainty region is known. We also present an auxilliary result that shows that our choice of uncertainty set and parameters unifies the canonical $p$-median and $p$-center problems, as well as their robust analogues.

 \item
We present a bound on the approximation error induced by the discretization of the service area, which depends on the granularity of the discretization and the configuration of the uncertainty regions. The bound can be used to quantify the suboptimality of the solutions produced by our discretized model with respect to the underlying continuous problem.

 \item
We propose a row-and-column generation algorithm for solving the robust AED deployment problem. Our proposed algorithm operates by exploiting the structure of the uncertainty set, which allows it to substantially outperform a standard Benders-based algorithm and a single-stage MIP reformulation. Notably, our algorithm decouples the size of the master problem from the number of scenarios used to discretize the service area, which allows us to efficiently solve problems where the discretization is extremely fine. Combined with the bounds on the discretization error, the row-and-column generation algorithm allows us to efficiently obtain provably near-optimal solutions to the underlying continuous problem.

 \item
 We present results from an extensive numerical study on the placement of public AEDs using real cardiac arrest data from Toronto, Canada. Our results demonstrate that accounting for uncertainty in cardiac arrest locations can decrease the distance between cardiac arrest victims and the nearest AED by 9-15\%, under both typical and worst-case realizations of the uncertainty.
  \end{enumerate}

\section{Related Literature}
Our paper primarily contributes to the medical literature on cardiac arrest and public access defibrillation, and the operations research literature on facility location under uncertainty. We also briefly review the connection between our model and the general literature on optimization under uncertainty. 

{\bf Public access defibrillation.} Most of the literature in the medical community on AED location has focused on identifying building types that have a high incidence of cardiac arrest, such as shopping malls \citep{becker1998public,engdahl2005localization,brooks2013determining}.  However, relying on building type alone to guide AED locations has several drawbacks. First, it does not account for cardiac arrests that occur outdoors, since these events cannot be assigned to a building type. Second, a high risk building type in one city may not be high risk in another. Lastly, identifying high risk buildings does not address the question of how to optimally deploy a limited number of defibrillators. By contrast, our proposed approach is guided by historical cardiac arrest data and is agnostic to building type, which makes it generalizable to any city.

There is a recent but limited body of work on optimization-based approaches to placing AEDs in public locations \citep{siddiq2012,chan2013,chan2014,sun2016overcoming}. Our paper differs from and extends this literature in a few substantive ways. First, previous approaches have solely optimized over exact historical cardiac arrest locations. While these approaches serve as a good first step toward improving public access defibrillation, they do not explicitly account for uncertainty in future cardiac arrest locations. Our paper is the first to explicitly incorporate the uncertainty in cardiac arrest locations into the optimization model. We find that accounting for cardiac arrest location uncertainty leads to improved outcomes based on several performance metrics, compared to a baseline method of optimizing over historical locations only.  Second, existing AED models focus on maximizing cardiac arrest {\it coverage}, which is the number of cardiac arrests that are within some fixed distance (e.g., 100m) of an AED. While this approach can be effective for the {\it median} cardiac arrest patient, it can have the adverse effect of neglecting cardiac arrest victims at the tail of the distance distribution. By contrast, our modeling approach involves a risk-based objective function (inspired by ambulance response time guidelines), which allows us to mitigate the risk of unacceptably large distances between cardiac arrest victims and the nearest AED. Lastly, existing work on AED location has focused on identifying general hotspots for cardiac arrest within a large service area (e.g., \cite{chan2014} identify historical clusters of cardiac arrests within a city). We complement this work by taking a more tactical approach to the AED deployment problem, by focusing on the precise placement of AEDs within a small service area consisting of several city blocks.

{\bf Facility location under uncertainty.}  Facility location problems under uncertainty have received considerable attention, particularly with respect to uncertainty in demand node weights and edge lengths \citep{owen1998,snyder2006,baron2011}, and more recently in the risk of service disruptions at the facilities \citep{lim2010,cui2010,shen2011}. Our model differs from the existing literature in that we consider the arrival of demand points over a continuous, planar service area, rather than being restricted to a small set of locations. As a result, we do not aggregate demand into a small set of known locations, which is a common approach in facility location models related to emergency response \citep{cho2014simultaneous,erkut2008,brotcorne2003} and in general \citep{francis2000aggregation}. Previous work on demand location uncertainty in facility location problems is limited and has generally focused on the placement of a single facility \citep{cooper1978,drezner1989,averbakh2005}. To the best of our knowledge, this paper is the first to extend the literature on demand location uncertainty to the general case of siting multiple facilities, which we show can be done in a computationally tractable manner. 

Our paper is most closely related to work by \cite{baron2011}, who also present a robust optimization model for facility location. In their paper, demand is again restricted to a set of known locations, but the vector of demands over all locations is only known to live within a given uncertainty set. In contrast, we model total demand as being known (i.e., normalized to 1), and consider uncertainty in how this demand is distributed in space.

Our work is also related to facility interdiction problems, which focus on protecting against potential disruptions at facilities. Facility interdiction models, which, like ours, take the form of min-max optimization problems, were first considered by \cite{church2007protecting} and later by \cite{scaparra2008bilevel}, \cite{liberatore2011analysis} and \cite{losada2012stochastic}. With respect to CVaR, the only other work to use it in a facility location context is the \textit{mean-excess} model presented in \citet{chen2006}, which focuses on demand level uncertainty.

{\bf Optimization under uncertainty.} 
Our main formulation takes the form of a robust optimization problem in which a large set of scenario probabilities are uncertain and only known to reside within a given polyhedron. In the special case where there is no uncertainty (i.e. the uncertainty set is a singleton representing the true distribution), our formulation simplifies to an extensive form stochastic program with a large set of scenarios \citep{birge2011introduction}. Other papers which model uncertainty in probabilities using polyhedral uncertainty sets include \cite{chan2013adaptive}, who consider radiation therapy with uncertainty in patient breathing patterns, and \cite{farias2013nonparametric}, who consider uncertainty in consumer choice probabilities. \cite{ben2013robust} present a general framework which uses $\phi$-divergences to model distributional uncertainty sets, and obtain polyhedral sets as a special case. For a review of the broader robust optimization literature, we refer the reader to \cite{bental2009} and \cite{bertsimas2011}.

Since the source of uncertainty in our model stems from an unknown discrete probability distribution, our model can also be interpreted as an instance of distributionally robust optimization, which refers to optimization problems in which only a partial description of a probability distribution is available \citep{scarf1958min,dupacova1980minimax,birge1987computing,shapiro2002minimax,delage2010distributionally,goh2010distributionally,xu2012distributional}. Recently, \cite{wiesemann2013distributionally} proposed a canonical framework for distributionally robust convex optimization, and show that it subsumes well-known robust optimization models. A common approach to modeling distributional uncertainty sets is to impose constraints on the moments of the distribution (e.g., \cite{ghaoui2003worst,delage2010distributionally}). Since we instead model distributional uncertainty using polyhedral sets, from a technical perspective our model is more aligned with the literature on robust linear optimization (e.g., \cite{bertsimas2004}).

Due to our use of a CVaR objective function, our main formulation takes the form of a two-stage (min-max-min) robust optimization model, where the inner minimization problem is used to compute CVaR. \cite{gabrel2011} and \cite{zeng2013} also present two-stage models for facility location problems, where the focus is on uncertainty in demand level rather than location. Other settings in which two-stage robust optimization models have been considered include network flows \citep{atamturk2007,ordonez2007}, power systems \citep{bertsimas2013,zeng2012}, and military applications, including the defense of critical infrastructure \citep{brown2006,alderson2011} and network interdiction \citep{brown2009}. 

Lastly, our paper has some similarities with work by \cite{carlsson2013robust}, in that they also consider a distributionally robust optimization problem in a spatial setting where demand points arrive in a continuous service area. Assuming the mean and covariance of the demand distribution is known, they focus on a vehicle routing problem where the goal is to partition the service area into a set of subregions such that the worst-case load for any vehicle across all subregions is minimized. By comparison, we assume that a fixed set of subregions is given, and the only distributional information available is the probability of an arrival in each subregion. 
 

\section{Model}\label{sec:model}
Let $\mathcal{A} \subset \mathbb{R}^2$ represent the continuous service area over which cardiac arrests arrive. Let $\mathcal{I}$ be an index set for $m$ candidate sites, and let ${\bf y} \in \{0,1\}^m$ be a binary vector, where $y_i=1$ indicates the presence of an AED in the $i^{th}$ candidate site. Let the function $d: \mathcal{A} \times \mathcal{A} \rightarrow \mathbb{R}_+$ measure the distance between two locations in $\mathcal{A}$. With a slight abuse of notation, let ${\bf y}(a) \in \mathcal{A}$ represent the nearest facility to a point  $a \in \mathcal{A}$, so that $d(a,{\bf y}(a))$ is the distance between a point $a \in \mathcal{A}$ and its nearest facility.

Let $\xi$ be a random vector representing the location of the next cardiac arrest event within $\mathcal{A}$, with distribution $\mu$. The distribution $\mu$ can be interpretted as describing the risk of a cardiac arrest event occurring over $\mathcal{A}$.  Without loss of generality, we assume $P_\mu(\xi \in \mathcal{A}) = 1$. For a set of AED locations ${\bf y}$, $d(\xi,{\bf y}(\xi))$ represents the distance between the next demand arrival and its nearest facility. Since $d(\xi,{\bf y}(\xi))$ is itself random due to $\xi$, we shall refer to the distribution of $d(\xi, {\bf y}(\xi))$ as the {\it distance distribution} induced by ${\bf y}$. We can think of the distribution of $d(\xi,{\bf y}(\xi))$ as fully describing the performance of the AED deployment.

As discussed in Section 1, we wish to place AEDs in a manner that mitigates the risk of large distances between cardiac arrest victims and the nearest AED. Restated, our goal is to identify an AED deployment ${\bf y}$ that controls the right-tail of the distance distribution, $d(\xi,{\bf y}(\xi)$. Since placing an AED in every possible candidate site is prohibitively expensive, we assume throughout that up to $P$ AEDs are available for deployment. If the distribution $\mu$ is known, then the problem of minimizing the CVaR of the distance distribution is given by 
\begin{align} \underset{{\bf y} \in \mathbf{Y}}{\text{minimize}} \;\; &\text{CVaR}_{\mu}[d(\xi,{\bf y}(\xi))] \label{model:muknown} 
\end{align}
where 
$$\textbf{Y} = \left\{ {\bf y} \in \{0,1\}^m \; \bigg| \; \underset{i \in \mathcal{I}}{\sum} y_i = P\right\}.$$ 
CVaR is closely related to the value-at-risk (VaR) measure, both of which have origins in the finance literature (see \cite{rockafellar2000,rockafellar2002}). For a given loss distribution, $\beta$-VaR represents the smallest value $\alpha$ such that the loss does not exceed $\alpha$ with probability $1-\beta$. By contrast, $\beta$-CVaR represents the expected loss conditional on the loss exceeding $\beta$-VaR.

In practice, the distribution of cardiac arrest locations is unlikely to be known precisely. Since -- for a given location -- cardiac arrests are low probability events, optimizing AED locations solely with respect to limited historical data (effectively a sample average approximation) may result in a poor deployment of AEDs. Moreover, any estimate of $\mu$ from historical data is also likely to be subject to uncertainty, especially if a limited amount of data is available.

To account for this uncertainty, we instead assume that the distribution $\mu$ is only known to belong to an uncertainty set, $\mathcal{U}$. A key question at this juncture is how to structure the set $\mathcal{U}$. To obtain a formulation that is both tractable and effective at capturing uncertainty in $\mu$, we take the following approach. Suppose the service area $\mathcal{A}$ can be divided into (possibly overlapping) {\it uncertainty regions} $\mathcal{A}_1, \mathcal{A}_2, \ldots , \mathcal{A}_{|\mathcal{J}|}$, where $\mathcal{A} =  \mathcal{A}_1 \cup \mathcal{A}_2 \cup \ldots \cup \mathcal{A}_{|\mathcal{J}|}$. Suppose now that the only available information about the distribution $\mu$ is the probability of the next cardiac arrest occuring in each of the regions $\mathcal{A}_1,\ldots,\mathcal{A}_{|\mathcal{J}|}$, which we denote by $\lambda_1,\ldots,\lambda_{|\mathcal{J}|}$. The worst-case CVaR can now be minimized by solving the problem
\begin{align} \underset{{\bf y} \in \mathbf{Y}}{\text{min}} \; \underset{{\bf \mu \in \mathcal{U} }}{\text{max}} \;\; & \text{CVaR}_{\mu}[d(\xi,{\bf y}(\xi))], \label{minmaxcvar}
 \end{align}
where $$\mathcal{U} = \{  \mu \; | \;  \mathbb{P}_{\mu}(\xi \in \mathcal{A}_1) = \lambda_1, \ldots,\mathbb{P}_{\mu}(\xi \in \mathcal{A}_{|\mathcal{J}|}) = \lambda_{|\mathcal{J}|}\},$$
is the distributional uncertainty set. Note that since the uncertainty regions are permitted to overlap, the probabilities $\lambda_1,\ldots,\lambda_{|\mathcal{J}|}$ may sum to a value greater than 1. We assume throughout that the uncertainty regions are given. In practice, the design of the uncertainty regions $\mathcal{A}_j$ should be guided by several considerations such as tractability, availability of historical data, and model interpretability (we discuss these considerations in greater detail in Section \ref{sec:casestudy}). 

There are several reasons behind our choice of the uncertainty set $\mathcal{U}$. First, we posit that obtaining estimates of the cardiac arrest risk at an aggregate (i.e., uncertainty region) level is less onerous than estimating the distribution $\mu$ directly, especially if data on historical cardiac arrest locations is limited. Estimating the parameters $\lambda_1,\ldots,\lambda_{|\mathcal{J}|}$ in lieu of $\mu$ itself allows us to partially capture the underlying distribution while accounting for uncertainty. This approach is well-aligned with the notion that aggregate forecasts tend to be more reliable than point forecasts when dealing with uncertainty \citep{simchi2004, sheffi2005, bertsimas2006robust, nahmias2009}. In addition, our choice for the structure of $\mathcal{U}$ lends itself to a tractable mathematical programming formulation of \eqref{minmaxcvar}. The key to this tractability is the observation that for any ${\bf y}$, the worst-case distribution in $\mathcal{U}$ is supported on a relatively small number of locations in $\mathcal{A}$. This sparsity can in turn be exploited through a decomposition algorithm, which allows us to solve the robust optimization problem efficiently (see Section 4).  Our choice of $\mathcal{U}$ also provides a natural generalization of classical facility location models, which have a rich and extensive history.  In particular, classical models typically assume demand is concentrated at $n$ weighted demand points (e.g., a weight of $\lambda_j$ for demand point $j = 1, \ldots, n$), which is recovered from our model by simply shrinking each $\mathcal{A}_j$ to a singleton. Thus, a natural interpretation of our distributionally robust model is as an extension of classical models where each demand point lives in a region $\mathcal{A}_j$, wherein the true demand location will be realized after facility siting. We discuss the connection between our model and the classical location literature in more detail in Section \ref{sec:generality} of the electronic companion. In particular, we show that our model unifies the classical $p$-median and $p$-center problems, as well as their robust analogues.

We note here that our choice of uncertainty set bears similarities to the uncertainty set discussed in the paper by \cite{wiesemann2013distributionally}, who present a canonical modeling approach for distributionally robust optimization problems. The authors formulate the uncertainty set using constraints on the probability that a random vector is realized within various convex sets. They show that within their framework, simply checking the whether the distributional uncertainty set is empty can be a strongly NP-hard problem, unless the constituent convex sets satisfy a particular nesting condition. This hardness follows from the generality of the uncertainty set, which can be shown to subsume integer programming, in the absence of the nesting condition.  Our setting is considerably more structured. By discretizing the service area into a large number of scenarios, we obtain a representation of $\mathcal{U}$ that is polyhedral. As a consequence, the worst-case distribution can be identified by solving a linear program. 

 In the remainder of this section, we show how our discretization approach allows us to formulate the optimization problems \eqref{model:muknown} and \eqref{minmaxcvar} as tractable mathematical programs. We then propose a bound on the error introduced by the discretization.

\subsection{Known cardiac arrest distribution}
For ease of exposition, we first consider the case where $\mu$ is known. A standard approach for formulating stochastic programs is to discretize the outcome space into a finite set of {\it scenarios} and to assign a probability to each scenario. Since in our setting, $\xi$ is the random location of the next cardiac arrest event, each scenario (i.e., realization of $\xi$) has a one-to-one correspondence with a particular geographic location in $\mathcal{A}$. Accordingly, let $\Xi \subset \mathcal{A}$ be a discrete set of locations that approximate the continuous service area $\mathcal{A}$, so that $\Xi$ is a discrete approximation of the outcome space of $\xi$, and each element of $\Xi$ corresponds to a distinct location in $\mathcal{A}$. Letting $\mathcal{K}$ be an index set for the scenarios in $\Xi$, we can write $\Xi = \left\{ \xi_1, \xi_2, \ldots, \xi_{|\mathcal{K}|} \right\}$ to represent the set of all possible realizations of $\xi$. Let $u_k$ be the probability that the next cardiac arrest will occur at location $\xi_k$. Since in this discrete setting the possible cardiac arrest locations are restricted to $\Xi$, we have $\sum_{k \in \mathcal{K}} u_k = 1$.

Let the parameter $d_{ik}$ represent the distance between candidate site $i$ and scenario $\xi_k$ based on the distance function $d(\cdot,\cdot)$. Let $z_{ik}$ be an assignment variable that is equal to 1 if a cardiac arrest realized at location $\xi_k$ is assigned to an AED at location $i$. The feasible set of assignments is given by
$$\textbf{Z}(\textbf{y}) = \left\{ z_{ik} \ge 0, \; i \in \mathcal{I}, k \in \mathcal{K} \;\bigg|\; \underset{i \in \mathcal{I}}{\sum} \; z_{ik} = 1, k \in \mathcal{K}; \; z_{ik} \le y_i, \;  i \in \mathcal{I}, k \in \mathcal{K} \right\}.$$ 
Note that based on the discretization, we have $d(\xi_k,{\bf y}(\xi_k)) = \min_{{\bf z} \in {\bf Z}({\bf y})} \sum_{i \in \mathcal{I}} d_{ik}z_{ik}$. Hence, for fixed ${\bf y}$, $\text{CVaR}_{\mu} [d(\xi,{\bf y}(\xi))]$ can be expressed as \citep{rockafellar2000,rockafellar2002}:
 \begin{equation}\text{CVaR}_{\mu} [d(\xi,{\bf y}(\xi))] \approx \underset{{\bf z}, \alpha}{\text{minimize}} \;\; \alpha + \frac{1}{(1-\beta)} \sum_{k \in \mathcal{K}} u_k \max \left\{ \sum_{i \in \mathcal{I}} d_{ik} z_{ik} - \alpha, 0  \right\}. \label{cvarminus}
 \end{equation}
Strictly speaking, the quantity in \eqref{cvarminus} is known as $\beta$-CVaR$^{-}$, which is an approximation of $\beta$-CVaR. This is a commonly used approximation for the case of discrete distributions -- see \cite{rockafellar2000} for further details. The nominal optimization problem \eqref{model:muknown} for a known distribution can now be written as
 \begin{align}\label{model:muknownprogram}
     \underset{\textbf{y},\textbf{z},\alpha}{\text{minimize}}  \;\; & \alpha + \frac{1}{(1-\beta)} \sum_{k \in \mathcal{K}} u_k \max \left\{ \sum_{i \in \mathcal{I}} d_{ik} z_{ik} - \alpha, 0  \right\}  \tag{\textsf{N-AED}}  \\
     \mbox{subject to} \;\; & {\bf z} \in {\bf Z}({\bf y}),  \; {\bf y} \in {\bf Y}, \notag \\
     & \alpha \ge 0, \notag 
     \end{align}
     which is a mixed-integer linear program after linearizing the maximization term in the objective.
     
\subsection{Distributionally robust model for unknown cardiac arrest distribution}
We now extend formulation \ref{model:muknownprogram} to accomodate uncertainty in the distribution $\mu$. Given the uncertainty regions $\mathcal{A}_1,\ldots,\mathcal{A}_{|\mathcal{J}|}$, we construct scenario sets $\Xi_1, \Xi_2, \ldots, \Xi_{|\mathcal{J}|}$, where $$\Xi_j = \{\xi_k, k \in \mathcal{K}\; | \; \xi_k \in \mathcal{A}_j\}$$
and $ \bigcup_{j=1}^n \Xi_j  = \Xi$. Let $\mathcal{K}_j$ index the scenarios in $\Xi_j$, and let $\mathcal{J}$ index the set of uncertainty regions. Since in the discretized setting we assume a cardiac arrest arriving in $\mathcal{A}_j$ can only be realized at one of the locations in $\Xi_j$, we have $P_{\mu}(\xi \in \Xi_j) = P_{\mu}(\xi \in \mathcal{A}_j) = \lambda_j.$   The worst-case CVaR of the distance distribution given a feasible ${\bf y},{\bf z}$ is then given by optimal value of the following max-min problem:
  \begin{align} \max_{{\bf u}} \min_{\alpha} \;\; & \alpha + \frac{1}{(1-\beta)}  \sum_{k \in \mathcal{K}} u_k \max \left\{ \sum_{i \in \mathcal{I}} d_{ik} z_{ik} - \alpha, 0  \right\} \notag  \\
 \text{subject to} \;\; & \sum_{k \in \mathcal{K}} u_k = 1,  \label{model:maxmincvar} \\
 & \sum_{k \in \mathcal{K}_j}u_k = \lambda_j, \quad j \in \mathcal{J}, \notag \\
 & u_k \ge 0, \quad k \in \mathcal{K}, \notag \\
 & \alpha \ge 0. \notag
  \end{align}
For conciseness, define 
$$ \mathbf{U} = \left\{ {\bf u} \in \mathbb{R}_+^{|\mathcal{K}|} ~\bigg|~ \sum_{k \in \mathcal{K}} u_k = 1, \sum_{k \in \mathcal{K}_j} u_k = \lambda_j, \; j \in \mathcal{J} \right\}. $$
Note that the set ${\bf U}$ may be empty if the parameters $\lambda_1,\ldots,\lambda_{|\mathcal{J}|}$ are selected arbitrarily (e.g., consider two regions $\mathcal{A}_1$, $\mathcal{A}_2$ where $\mathcal{A}_1 \cap \mathcal{A}_2 = \emptyset$ but $\lambda_1 = \lambda_2 = 1$). However, if $\lambda$  is estimated appropriately from historical cardiac arrest data, then nonemptiness of ${\bf U}$ is guaranteed (see Proposition \ref{prop:MLE}). We discuss the design of uncertainty regions and the estimation of $\lambda$ further in Section \ref{sec:MLE}. Note also that our discretization scheme allows us to model generally shaped uncertainty regions. This flexibility allows physical obstacles to be incorporated into the uncertainty regions, such as certain private buildings where public cardiac arrests cannot occur. 

 Minimizing the worst-case CVaR over ${\bf y},{\bf z}$, we now arrive at the following two-stage robust optimization problem:
  \begin{align} \label{model:robustcvar2}
  \min_{{\bf y},{\bf z}} \max_{{\bf u}} \min_{\alpha} \;\; & \alpha + \frac{1}{(1-\beta)}  \sum_{k \in \mathcal{K}} u_k \max \left\{ \sum_{i \in \mathcal{I}} d_{ik} z_{ik} - \alpha, 0  \right\} \notag  \\
 \text{subject to} \;\; & {\bf u} \in {\bf U}, {\bf z} \in {\bf Z}({\bf y}), {\bf y} \in {\bf Y}, \tag{\textsf{R-AED}} \\
 & \alpha \ge 0. \notag
  \end{align}
With respect to the size of \textsf{R-AED}, in practice we expect $\mathcal{K}$ to be large (on the order of $10^3$ or larger), so that our discretization procedure reasonably approximates the underlying continuous problem. For example, in our case study (Sections 5 and 6), we take $|\mathcal{K}| = 2,600$. Similarly, the number of candidate sites $|\mathcal{I}|$ in our case study is 120, and the number of uncertainty regions $|\mathcal{J}|$ is 15.
 
It is worth discussing how formulation \textsf{R-AED} compares to other robust location models with demand uncertainty. The model closest to ours is the one proposed by \cite{baron2011}, who consider demand-level uncertainty at each of a fixed set of demand locations. Our model can also be viewed as modeling demand-level uncertainty, although there are some key differences with \cite{baron2011}. In formulation \textsf{R-AED}, the discrete uncertainty set ${\bf U}$ forces the aggregate demand from the locations $\mathcal{K}_j$ to be $\lambda_j$, and normalizes the total demand over all locations in $\mathcal{K}$ to be 1. This is in contrast to the approach of \cite{baron2011}, who model the uncertainty by allowing the vector of demands over all locations to reside within a box or ellipsoidal uncertainty set. This difference in how demand uncertainty is modeled turns out to be consequential. Based on our choice of the (polyhedral) uncertainty set, the worst-case distribution in ${\bf U}$ will generally be supported on a relatively small number of locations, compared to the total number of scenarios. This means that, in the worst-case distribution, the demand $u_k$ in the vast majority of locations in $\mathcal{K}$ will be 0. Moreover, increasing the size of the scenario $\mathcal{K}$ has no effect on the number of locations that support the worst-case distribution in ${\bf U}$. This sparsity and independence from the total number of scenarios in $\mathcal{K}$ is exploited by our row-and-column generation algorithm to obtain solutions to \textsf{R-AED} (Section 4). By contrast, in \cite{baron2011}, the uncertain demand at each location typically has a non-zero lower bound, which leads to all locations having non-zero demand, and thus precludes the application of our decomposition technique.

In the most general case, formulation \ref{model:robustcvar2} can be interpreted as a facility location model with a single demand point (i.e., the location of the ``next" cardiac arrest) and stochastic edge lengths, where $d_{ik}$ gives the length of the edge between the demand point and candidate site $i$ under scenario $k$. Under this interpretation, $\lambda_j$ represents the probability that one of the scenarios in the set $\mathcal{K}_j$ occurs. Our model differs from the existing literature on stochastic edge lengths since we assume that the vector of probabilities ${\bf u}$ is itself uncertain, and only known to reside within a polyhedral uncertainty set.

Formulation \ref{model:robustcvar2} can also be shown to unify the well known $p$-median and $p$-center location problems, as well as their robust analogues under demand location uncertainty. We formalize this connection in Section \ref{sec:generality} of the electronic companion.

\subsection{Discretization error bounds}\label{sec:approximation}
The model \ref{model:robustcvar2} is an approximation in the sense that it restricts the continuous demand distribution to be supported on a discrete set. As a result, the ``true" worst-case CVaR in the continuous setting may be underestimated by the optimal value of formulation \ref{model:robustcvar2}. In this section, we present bounds on this underestimation, which we refer to as the \textit{discretization error}. The discretization error can be interpreted as the optimality gap associated with an optimal solution of formulation \ref{model:robustcvar2} with respect to the underlying continuous problem \eqref{minmaxcvar}. 

To define and prove the discretization error bounds, we first construct a partition of the service area $\mathcal{A}$. The partitioning is required for technical purposes only, and can be applied to any general setting with overlapping uncertainty regions. Note that any configuration of the uncertainty regions $\mathcal{A}_j$ implies a partitioning of the full service area $\mathcal{A}$ into potentially more than $n$ ``subregions". For example, suppose $\mathcal{A} = \mathcal{A}_1 \cup \mathcal{A}_2$, where $\mathcal{A}_1 \cap \mathcal{A}_2 \neq \emptyset$. Then we consider three disjoint subregions: $\mathcal{A}'_1 = \mathcal{A}_1 \setminus \mathcal{A}_2$, $\mathcal{A}'_2 = \mathcal{A}_2 \setminus \mathcal{A}_1$, and $\mathcal{A}'_3 = \mathcal{A}_1 \cap \mathcal{A}_2$.  In general, let $\mathcal{A}'_1,\ldots,\mathcal{A}'_{|\mathcal{R}|}$ represent these subregions, and let $\Xi'_1,\ldots,\Xi'_{|\mathcal{R}|}$ be the discrete counterparts to the $\mathcal{A}'_r$. Let $\mathcal{R}$ be the index set for the subregions and let $\mathcal{K}_r$ index the scenarios in subregion $r$. Since by definition the $\mathcal{A}'_r$ form a partition of $\mathcal{A}$, it follows that $\Xi'_{r} \cap \Xi'_{r'} = \emptyset$ for all $r \neq r'$. Note that in the most general case, we may have $|\mathcal{R}| = 2^{|\mathcal{J}|}$, if every uncertainty region partially overlaps with every other uncertainty region. However, this extreme case is unlikely to arise naturally in the context of AED deployment and other facility location problems. Note also that if there is no overlap among the $\mathcal{A}_j$, then $|\mathcal{R}| = |\mathcal{J}|$. To construct the bound, we first require the following assumptions.
    \begin{assumption} $\Xi'_r = \mathcal{A}'_r \cap \mathcal{L}$ for all $r \in \mathcal{R}$, where $\mathcal{L}$ is a square lattice in the plane with a grid spacing of length $\sigma$.
    \end{assumption}

      \begin{assumption}
      The shape of the subregion $\mathcal{A}'_r$ is such that for any point $a \in \mathcal{A}'_r \setminus \Xi'_r$, at least one of the four points in $\mathcal{L}$ that define the smallest square containing $a$ is inside $\mathcal{A}'_r$.
      \end{assumption}

          \begin{figure}[t]
\centering
\subfigure[Assumption 2 satisfied.]{%
\includegraphics[scale=0.6]{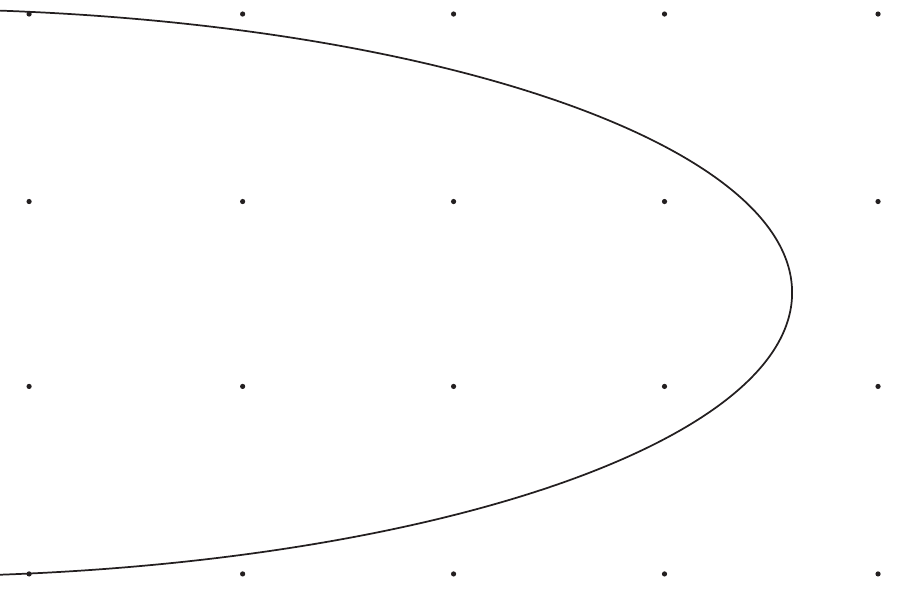}
\label{fig:subfigure1}}
\hspace{10mm}
\subfigure[Assumption 2 violated.]{%
\includegraphics[scale=0.6]{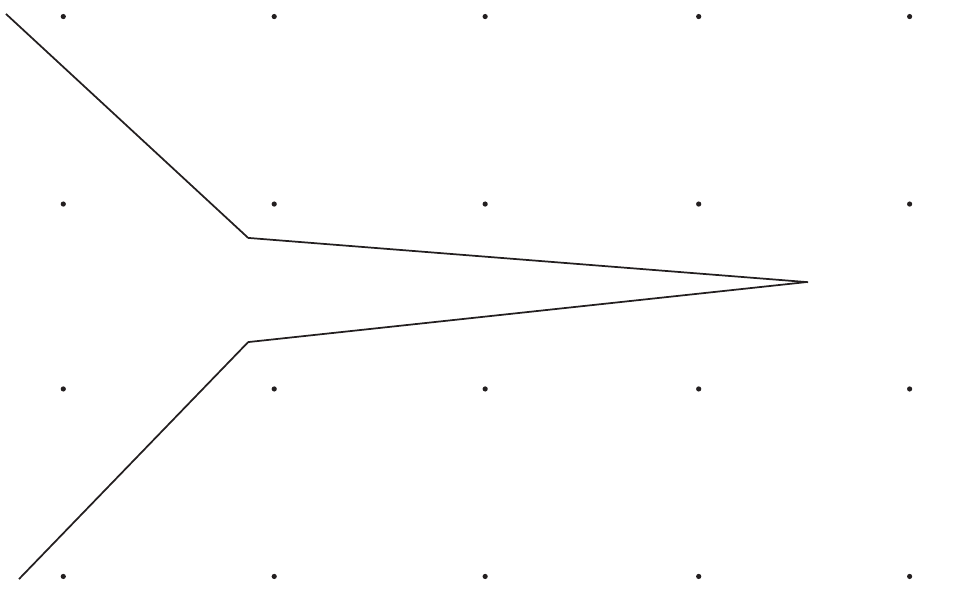}
\label{fig:subfigure2}}
\caption{Examples for Assumption 2.}
\label{fig:assumptionshape}
\end{figure}
Assumption 1 states that the uncertainty regions are discretized using a square lattice, which simplifies the analysis considerably without loss of generality. Assumption 2 is meant to exclude pathologically shaped uncertainty regions relative to the size of the discretization (see Figure \ref{fig:assumptionshape}), and can be satisfied with a sufficiently granular discretization of the uncertainty region. We now define a parameter $\ell(\sigma)$ with the following property:
\begin{align} \label{ell}
\ell(\sigma) \ge \underset{a \in \mathcal{A}'_r}{\operatorname{max}} \; \underset{\xi \in \Xi'_r}{\operatorname{min}} \; d(a,\xi), \quad \text{ for all } r \in \mathcal{R}.
\end{align}
Here, $\ell(\sigma)$ is an upper bound on the maximum distance between any point in $\mathcal{A}'_r$ and its closest scenario location in $\Xi'_r$. To make the expression in \eqref{ell} concrete, we provide some examples of the parameter $\ell(\sigma)$ in Remark 1, depending on the choice of the distance metric. 
\begin{remark}
Let Assumptions 1 and 2 hold. Let $\|\cdot\|_1$, $\|\cdot\|_2$ and $\|\cdot\|_\infty$ be the Euclidean, rectilinear and maximum metric on $\mathbb{R}^2$. 
\begin{enumerate}[1.]
\item
If $d(a,b) = \|a - b\|_1$, then $\ell(\sigma) = 2\sigma$.
\item
If $d(a,b) = \|a - b\|_2$, then $\ell(\sigma) = \sqrt{2}\sigma$.
\item
If $d(a,b) = \|a - b\|_\infty$, then $\ell(\sigma) = \sigma$.
\end{enumerate}
\end{remark}
Note that setting $\ell(\sigma) = 2\sigma$ is also valid for any distance function based on a $p$-norm with $p \ge 1$, since $\| \cdot \|_p$ is non-increasing in $p$. Let $Z_D$ be the optimal value of formulation \ref{model:robustcvar2}, $Z_C$ be the optimal value of the associated continuous problem \eqref{minmaxcvar}, and $$\Delta := \frac{Z_C - Z_D}{Z_C}$$ be the discretization error. We now present the main bound.  

      \begin{theorem}\label{thm:discretizationerror}
         For any $\varepsilon > 0$, if $\ell(\sigma) \le \varepsilon (1 - \beta) Z_D$, then $\Delta \le \varepsilon$.
        \end{theorem}
Theorem \ref{thm:discretizationerror} characterizes how finely the uncertainty regions must be discretized to achieve a certain discretization error. Since $\ell(\sigma)$ vanishes with $\sigma$, the discretization error can be made arbitrarily small by selecting an appropriate value for $\sigma$.  Thus, we can obtain near-optimal solutions to the underlying continuous uncertainty problem by solving formulation \ref{model:robustcvar2} with a sufficiently fine discretization (and thus a sufficiently large number of scenarios). We can also re-arrange the expression in Theorem \ref{thm:discretizationerror} to calculate an upper bound on $Z_C$ as a function of $\sigma$ and $Z_D$. The upper bound $\bar{Z}_C$ is constructed as follows:
    \begin{equation}\label{Zbound}
    Z_{C}  \le \bar{Z}_{C} = Z_D(1+\varepsilon),\;\; \textnormal{where } \varepsilon = \frac{\ell(\sigma)}{Z_D(1 - \beta)}.
    \end{equation}
In the next section, we outline a decomposition technique which enables us to compute a tight bound $\bar{Z}_C$ by solving instances of \ref{model:robustcvar2} with an extremely large number of scenarios. 
\section{Row-and-column generation algorithm}\label{sec:solution}
It is well known that two-stage robust optimization problems with the structure of \eqref{minmaxdual} can be solved using Benders-based decomposition algorithms (see, for example, \cite{brown2006} \cite{brown2009} \cite{alderson2011}). These approaches decompose the formulation into a master problem, which is a relaxation of the original problem, and a subproblem, which generates constraints to be added to the master problem. The master and subproblems are solved iteratively until convergence to a provably optimal solution. However, when applied to \eqref{minmaxdual}, the number of constraints and decision variables of the (integer) master problem using these approaches depends on the number of scenarios $|\mathcal{K}|$. As a consequence, these decomposition approaches can fail to be tractable if $|\mathcal{K}|$ is even moderately sized. We overcome this issue by proposing an algorithm that decouples the size of the master problem from the number of scenarios. A key consequence of this decoupling is that the tractability of the resulting formulation no longer depends on the granularity of the discretization, which, when combined with the bounds in Section \ref{sec:approximation}, allows us to approximate the continuous problem using an extremely large number of scenarios. A similar row-and-column generation algorithm for solving min-max optimization problems was independently proposed by \cite{zeng2013}. A key difference is that unlike the model in \cite{zeng2013}, which has an integer subproblem, our subproblem is a linear program.

We compare our proposed method with two other solution approaches that are natural candidates for solving  formulation \ref{model:robustcvar2} -- a Benders-based row generation algorithm and a mixed integer programming (MIP) reformulation of \ref{model:robustcvar2}.

    \subsection{Algorithm overview} \label{section:rowandcolumn}
We now provide an overview of the row-and-column generation algorithm. First, observe that the inner minimization problem of \eqref{model:maxmincvar} is a linear program. Letting ${\bf p}$ denote the dual vector for this linear program, the dual polyhedron -- which depends on the distribution ${\bf u}$ -- can be written as
\begin{align}\label{Pset}
{\bf P}({\bf u}) = \left\{p_k \ge 0, k \in \mathcal{K} ~\Bigg|~ \sum_{k \in \mathcal{K}} p_k \le 1, \; p_k \le u_k/(1-\beta), k \in \mathcal{K}\right\}. 
\end{align}
The optimization problem \eqref{model:maxmincvar} can now be reformulated as
\begin{align}\label{minmaxdual}
  \underset{\textbf{y},{\bf z}}{\operatorname{min}} \; \underset{\textbf{u},{\bf p}}{\operatorname{max}}\;\; & \sum_{i \in \mathcal{I}}   \sum_{k \in \mathcal{K}} d_{ik} z_{ik} p_k  \notag \\\vspace{0.2in}
 \mbox{subject to} \;\; & {\bf p} \in {\bf P}({\bf u}), {\bf u} \in {\bf U}, \\
    & {\bf z} \in \mathbf{Z}({\bf y}), {\bf y} \in {\bf Y}. \notag
   \end{align}
By introducing the auxilliary variable $t$, we obtain the following epigraph representation of \eqref{minmaxdual}:
\begin{subequations}\label{infinitedim}
\begin{align}
  \underset{\textbf{y},{\bf z},t}{\operatorname{minimize}} \;\; & t  \\
 \mbox{subject to} \;\; & t \ge \sum_{i \in \mathcal{I}}   \sum_{k \in \mathcal{K}} d_{ik} z_{ik} p_k, \quad {\bf p} \in \mathbf{P}({\bf u}), {\bf u} \in {\bf U}, \label{tinfinite} \\
    & {\bf z} \in \mathbf{Z}({\bf y}), {\bf y} \in \mathbf{Y}.
           \end{align}
    \end{subequations}
Note that formulation \eqref{infinitedim} is a semi-infinite optimization problem, since the constraint \eqref{tinfinite} is enforced over all possible ${\bf p}$. To obtain a feasible solution to \eqref{infinitedim}, we may instead solve a relaxation where the constraint \eqref{tinfinite} is enforced over a finite subset of $\mathbf{P}$. Let this subset be ${\bf P}_{|\mathcal{S}|} = \{{\bf p}^1,\ldots,{\bf p}^{|\mathcal{S}|}\}$, and let $\mathcal{S}$ be an index set for the vectors in ${\bf P}_{|\mathcal{S}|}$. To obtain a relaxation of \eqref{tinfinite} with a finite number of constraints, we replace the constraint set represented by \eqref{tinfinite} with the following set of $|\mathcal{S}|$ constraints:
\begin{align}\label{tfinite}
 t \ge \sum_{i \in \mathcal{I}}\sum_{k \in \mathcal{K}} d_{ik}z_{ik}p_k^s, \quad s \in \mathcal{S}.
\end{align}
Although replacing the constraint \eqref{tinfinite} with \eqref{tfinite} in formulation \eqref{infinitedim} produces a mixed-integer linear program, the resulting formulation may still be intractably large, due to the dependence of the number of constraints and variables on the size of the scenario set, $\mathcal{K}$. Consider now the set
\begin{align}
\mathcal{K}_+ = \left\{ k \in \mathcal{K} ~\Bigg|~ \sum_{s \in \mathcal{S}} p_k^s > 0 \right\},
\end{align}
and note that constraint \eqref{tfinite} can be written as
\begin{align}
t \ge \sum_{i \in \mathcal{I}}\sum_{k \in \mathcal{K}_+} d_{ik}z_{ik}p_k^s + \sum_{i \in \mathcal{I}}\sum_{k \in \mathcal{K} \setminus \mathcal{K}_+} d_{ik}z_{ik}p_k^s, \quad s \in \mathcal{S}.
\end{align}
Since by definition we have $\sum_{k \in \mathcal{K} \setminus \mathcal{K}_+} p_k^s = 0$ for all $s \in \mathcal{S}$, we may drop the second term on the right hand side of \eqref{tfinite}, which further simplifies the constraint to 
\begin{align}
t \ge \sum_{i \in \mathcal{I}}\sum_{k \in \mathcal{K}_+} d_{ik}z_{ik}p_k^s, \quad s \in \mathcal{S}. 
\end{align}
Observe that if $k \notin \mathcal{K}_+$, then the value of the corresponding assignment variable $z_{ik}$ has no impact on the objective function, since its coefficient is 0 for all $s \in \mathcal{S}$. We can therefore remove all $z_{ik}$ variables and the constraints $z_{ik} \le y_i$ and $z_{ik} \ge 0$, for $k \in \mathcal{K}\setminus \mathcal{K}_+$. This leads to the following relaxation of \eqref{infinitedim}: 
\begin{align}
  \underset{\textbf{y},{\bf z},t}{\operatorname{minimize}} \;\; & t  \notag \\
 \mbox{subject to} \;\; & t \ge \sum_{i \in \mathcal{I}}   \sum_{k \in \mathcal{K}_+} d_{ik} z_{ik} p_k^s, \quad s \in \mathcal{S}, \notag \\
    & \sum_{i \in \mathcal{I}}y_i = P, \notag \\
    & \sum_{i \in \mathcal{I}}z_{ik} = 1, \quad k \in \mathcal{K}_+, \label{masterprob} \tag{\textsf{R-AED-MP}} \\
    & z_{ik} \le y_i, \quad i \in \mathcal{I}, \; k \in \mathcal{K}_+, \notag \\
    & z_{ik} \ge 0, \quad i \in \mathcal{I}, \; k \in \mathcal{K}_+, \notag  \\
    & y_i \in \{0,1\}, \quad i \in \mathcal{I}. \notag
   \end{align}
The relaxed problem \ref{masterprob} yields a lower bound on the optimal value of the original problem \eqref{minmaxdual}. To tighten the relaxation, we can add a new vector {\bf p} to ${\bf P}_{|\mathcal{S}|}$ and re-solve \ref{masterprob}, with the aim of obtaining an improved lower bound. Let $\bar{{\bf y}}$ be a deployment obtained at a solution of \ref{masterprob}. We then construct a feasible assignment vector $\bar{{\bf z}}$ that assigns each scenario to its nearest AED as follows. Let $i^*(k) = \argmin_{ \{i \in \mathcal{I} | \bar{y}_i = 1 \}} \{d_{ik}\}$, which represents the closest AED to location $k$, given the deployment $\bar{{\bf y}}$. Then for each $k \in \mathcal{K}$, set 
\begin{align} \label{constructz}
\bar{z}_{i^*(k),k} = 1 \text{ and } \bar{z}_{ik} = 0 \text{ for all } i \neq i^*(k).
\end{align} 
It is straightforward to show that the $\bar{\bf z}$ given by \eqref{constructz} is also optimal for \textsf{R-AED-MP}. To identify a new ${\bf p}$, we then solve the following subproblem, where $\bar{{\bf y}},\bar{t}$ is obtained as the incumbent solution to \ref{masterprob} and $\bar{{\bf z}}$ is given by \eqref{constructz}:
\begin{align}
 \underset{\textbf{u},{\bf p}}{\operatorname{maximize}}\;\; & \sum_{i \in \mathcal{I}}   \sum_{k \in \mathcal{K}} d_{ik} \bar{z}_{ik} p_k  \label{subproblem} \tag{\textsf{R-AED-SP}} \\
 \mbox{subject to} \;\;  &{\bf p} \in {\bf P}({\bf u}), {\bf u} \in {\bf U}.   \notag
  \end{align}
To understand why the construction of $\bar{\bf z}$ using \eqref{constructz} is necessary, note that since \textsf{R-AED-MP} minimizes CVaR, the optimal assignment variables ${\bf z}$ obtained at a solution to \textsf{R-AED-MP} may not assign every location $k$ to its nearest AED, since only the locations in the tail of the distance distribution affect the optimal objective function value. As a consequence, solving \textsf{R-AED-SP} using the $\bar{\bf z}$ obtained directly from \textsf{R-AED-MP} may produce a solution to \textsf{R-AED-SP} that does not correspond to a worst-case distribution (with respect to the incumbent deployment ${\bf y}$). While in theory, the algorithm converges in finite time without reconstructing $\bar{\bf z}$ according to \eqref{constructz}, in practice we found that including this additional assignment step results in fewer iterations and a significant speed improvement.

Observe that formulation \textsf{R-AED-SP} is a linear program. Now let ${\bf u}^*,{\bf p}^*$ be an optimal solution to \ref{subproblem}, and let $Z_{\text{SP}}$ be the optimal value. If $\bar{t} \ge Z_{\text{SP}}$, then it follows that $\bar{t} \ge \sum_{i \in \mathcal{I}}   \sum_{k \in \mathcal{K}} d_{ik} \bar{z}_{ik} p_k$ for all ${\bf p} \in {\bf P}({\bf u}), {\bf u} \in {\bf U}$, and thus the incumbent solution is certifiably optimal to the original problem. If $\bar{t} < Z_{\text{SP}}$, then we introduce the variables $z_{ik}$, $k \in \mathcal{K}^*$ to the master problem, as well as the constraints $z_{ik} \ge 0$, $z_{ik} \le y_i$, $k \in \mathcal{K^*}$ and $t \ge \sum_{i \in \mathcal{I}}   \sum_{k \in \mathcal{K}} d_{ik} \bar{z}_{ik} p^*_k$, where $\mathcal{K}^* =  \{k \in \mathcal{K} ~|~ p^*_k > 0\}$. We then solve \ref{masterprob} again to obtain a new solution $\bar{{\bf y}},\bar{{\bf z}},t$, which begins a new iteration of the algorithm. Since there are a finite (but possibly large) number of rows and columns that can be generated, convergence after a finite number of iterations is guaranteed. The algorithm is summarized in Algorithm 1. 
       \begin{algorithm}
       \caption{Row-and-column generation}
       \begin{algorithmic}
       \STATE \textbf{Initialize} Set $\mathcal{S} = \{0\}$,  $s = 1$, and $\epsilon > 0$. Let $\mathcal{K}_+ = \emptyset$. 
       \STATE 1. Solve  \ref{masterprob} to obtain solution $(\bar{\textbf{y}},\bar{t})$ and objective value $Z_{\text{MP}}$.
              \STATE 2. Construct $\bar{{\bf z}}$ according to \eqref{constructz}.
       \STATE 2. Solve  \ref{subproblem} with fixed $\bar{\textbf{z}}$ to obtain $\bar{\textbf{p}}^{s}$ and objective value $Z_{\text{SP}}$.
       \STATE 3. \textbf{If} $\frac{Z_{\text{SP}} - Z_{\text{MP}}}{Z_{\text{SP}}} \le \epsilon$, terminate and return optimal deployment $\bar{\textbf{y}}$ and worst-case CVaR $\bar{t}$.
       \STATE \hspace{4mm} \textbf{else} 
       \hspace{6mm} \begin{enumerate}[i)]
       \item 
       Construct index set 
       $$\mathcal{K}^* = \{k \in \mathcal{K} \setminus \mathcal{K}_+ ~|~ \bar{p}^s_k > 0\}$$
       and set $\mathcal{K}_+ \leftarrow \mathcal{K}_+ \cup \mathcal{K}^*$.
       \item
        Add variables $z_{ik}$, $ i \in \mathcal{I}$, $k \in \mathcal{K}^*$ to master problem \ref{masterprob}.
        \item Add constraints 
        \begin{align*}
        &t \ge \sum_{i \in \mathcal{I}}\sum_{k \in \mathcal{K}_+}d_{ik}z_{ik}\bar{p}_k^s, \\
         &z_{ij}^{k} \le y_i, \quad  i \in \mathcal{I},\;  k \in \mathcal{K}^*, \\
         &\sum_{i \in \mathcal{I}} z_{ij}^{k} = 1, \quad  k \in \mathcal{K}^*, 
         \end{align*}
         to master problem \ref{masterprob}.
         \item 
          Set $\mathcal{S} \leftarrow \mathcal{S} \cup \{ s \} $. Increment $s$ and return to step 1.
        \end{enumerate}
       \end{algorithmic}
       \end{algorithm}
The intuition behind the algorithm is as follows. At a given iteration, the master problem only includes the subset of scenarios that are assigned a non-zero probability mass in the worst case distribution for at least one of the previous iterations. Based on the incumbent solution, the subproblem then identifies the worst-case distribution using the entire scenario set. The scenarios which support the worst-case distribution (obtained by solving the subproblem) are then added to the master problem on an as needed basis, along with the relevant assignment variables and related constraints. Due to this decomposition, at any iteration of the algorithm the size of the master problem \ref{masterprob} depends only on the number of candidate sites, $|\mathcal{I}|$, and the size of the restricted scenario set, $\mathcal{K}_+$, which may be substantially smaller than $\mathcal{K}$. In fact, due to the structure of the uncertainty set ${\bf U}$, the size of the restricted scenario set $\mathcal{K}_+$ can be shown to depend only  on the configuration of the uncertainty regions and the number of iterations completed by the algorithm, and is otherwise independent of the total number of scenarios in $\mathcal{K}$. This result is formalized in the following proposition.
\begin{proposition}\label{mastersize}
Assume that for any $k$ and $k'$ such that $k \neq k'$, we have $d_{ik} \neq d_{ik'}$ for all $i \in \mathcal{I}$. Then the number of constraints and variables in the master problem \textsf{R-AED-MP} in the $s^{th}$ iteration of Algorithm 1 is $O(|\mathcal{I}|(|\mathcal{J}|+1)s)$.
\end{proposition}
The assumption that $k \neq k'$ implies $d_{ik} \neq d_{ik'}$ for all $i \in \mathcal{I}$ is mild, and ensures that for any candidate site $i$, there are no exact ties for the closest location in $\mathcal{K}$ (indeed, this assumption can be relaxed if we instead assume that the solution to \textsf{R-AED-SP} is always a basic feasible solution, which is the case if the simplex method is used to solve \textsf{R-AED-SP}). It is worth noting that the size of the master problem \ref{masterprob} may still depend indirectly on $\mathcal{K}$ if the total number of scenarios is small (i.e., if $|\mathcal{R}|s > |\mathcal{K}|$). Additionally, Proposition 1 only guarantees that the size of the master problem in any given iteration is independent of $\mathcal{K}$, but does not guarantee that the {\it total} number of iterations required will be independent of $\mathcal{K}$. However, our numerical results indicate that increasing the number of scenarios generally has a minimal impact on the performance of the row-and-column generation algorithm, because the worst-case distribution is generally sparse with respect to $\mathcal{K}$. Note that unlike the master problem, the subproblem \ref{subproblem} still depends on the size of $\mathcal{K}$. However, the subproblem is a linear program with $|\mathcal{K}|$ decision variables, and we found that it can be solved efficiently for practically sized problems, even if the set $\mathcal{K}$ is quite large.

 For illustrative purposes, in Section \ref{totalscenarios} of the electronic companion we present a plot which depicts the total number of scenarios and the associated optimality gap for an example instance of the row-and-column generation algorithm.

       \subsection{Benchmark solution approaches}
       We compare the performance of the row-and-column generation algorithm to two benchmark approaches: a mixed-integer linear programming reformulation of \eqref{minmaxcvar} and a standard row generation algorithm for min-max problems. 
        \subsubsection{Mixed-integer linear programming reformulation}\label{sec:dualreformulation}
By formulating the dual problem of the inner maximization problem in  \eqref{minmaxdual}, we obtain the following mixed-integer linear program:
         \begin{align}
             \underset{\mathbf{y},\mathbf{z}, {\bf w}, \eta, \alpha, \mathbf{\gamma}}{\operatorname{minimize}} \;\;& \eta + \alpha + \sum_{j \in \mathcal{J}}\lambda_j w_j  \notag \\
             \mbox{subject to} \;\;
             & \alpha + \gamma_k \ge \sum_{i \in \mathcal{I}} \sum_{k \in \mathcal{K}} d_{ik}z_{ik}, \quad k \in \mathcal{K} , \notag \\
             & \eta - \frac{\gamma_k}{1-\beta} + \sum_{\{j | k \in \mathcal{K}_j\}}w_j \ge 0, \quad  k \in \mathcal{K},  \label{model:robustcvardual}\\
             & \mathbf{z} \in \mathbf{Z}(\mathbf{y}), \; \mathbf{y} \in \mathbf{Y}, \notag \\
             & \gamma_k \ge 0, \quad k \in \mathcal{K}, \notag \\
                          & \alpha \ge 0. \notag 
             \end{align}
If $\mathcal{K}$ is relatively small, then the MIP formulation \eqref{model:robustcvardual} can be solved using standard integer programming solvers. However, this formulation can become computationally intractable if a large number of scenarios are used to discretize the uncertainty regions.

\subsubsection{Row generation}\label{sec:solutionminmax}
A well known approach to solving min-max problems like \eqref{minmaxcvar} is to use a Benders-based row generation algorithm  \citep{brown2006,brown2009,alderson2011}. The row generation algorithm is conceptually similar to our row-and-column generation method outlined in Section \ref{section:rowandcolumn}, in that it also decomposes the problem into a master and subproblem which provide lower and upper bounds, respectively. However, a key difference with our approach is that the row generation algorithm carries the full set of columns (i.e., all scenarios in $\mathcal{K}$) in the master problem in each iteration. The number of variables and constraints in the master problem of the row generation algorithm in the $s^{th}$ iteration is thus $O(|\mathcal{I}||\mathcal{K}|)$ instead of $O(|\mathcal{I}||\mathcal{R}|s)$. As a result, the row generation algorithm may have poor performance if $|\mathcal{K}|$ is large.

\subsection{Comparison of solution approaches}
To demonstrate the performance of the row-and-column generation algorithm, we generated and solved several random problem instances using rectangular uncertainty regions and a square lattice for the discretization. We use Euclidean distances for all parameters. We solve each of these instances using the three solution methods discussed above, and focus on the impact that the number of scenarios has on the total solution time. All problems were implemented in MATLAB R2011a using YALMIP as the modeling language and CPLEX 12.1 with default parameter settings as the solver, on a 2.4 GHz quad-core CPU. 

In Table \ref{table:bounds}, we report numerical results on how the solution times of the MIP formulation (MIP), row generation algorithm (Row) and row-and-column generation algorithm (R+C) are impacted by the number of scenarios. We report results only for the cases where $\beta = 0$ and $\beta = 0.9$, and note that similar trends were observed for other values of $\beta$. We also compute the error bound $\varepsilon$ and the associated upper bound on the continuous problem, $\bar{Z}_C$, using \eqref{Zbound}.  Note that the error bounds for the $\beta = 0.9$ cases are larger than the $\beta = 0$ case, which is expected based on the bound expression given in \eqref{Zbound}.

Unlike the MIP and row generation algorithm, the performance of the row-and-column generation algorithm scales extremely well with the number of scenarios, and in many cases outperforms the benchmark approaches by one to two orders of magnitude with respect to the solution time. Note that the performance improvement is most pronounced at larger values of $|\mathcal{K}|$, due to the poor scalability of the MIP and row generation approaches. These results suggest that, when using the row-and-column generation algorithm, if a problem can be solved to optimality using a small number of scenarios, then the problem tends to remain tractable when using a fine discretization and an extremely large number of scenarios (e.g., just under 100,000 scenarios are used in the final row of Table \ref{table:bounds}). 


\begin{landscape}
\begin{table}[htbp]
  \parbox{7.85in}{\caption{Computational results showing effect of number of scenarios $(|\mathcal{K}|)$ on error bound and solution times (in CPU seconds) of MIP, row generation and row-and-column generation algorithm. Asterisk (*) indicates instance did not solve in 10,000 CPU seconds. } \label{table:bounds}}
  \centering
    \begin{tabular}{rrrrrrrrrrrrrrrrrrr}
    \toprule \\[-1.0em]
       &&&& & & & \multicolumn{6}{c}{$\beta = 0$} & \multicolumn{6}{c}{$\beta = 0.9$}  \\[-1.0em] \\
     \cline{7-12}   \cline{14-19}
   \\[-0.8em] \multicolumn{1}{c}{$|\mathcal{I}|$} & \multicolumn{1}{c}{$|\mathcal{J}|$} & \multicolumn{1}{c}{$P$}     &\multicolumn{1}{c}{$\mathcal{|K|}$} & \multicolumn{1}{c}{$\sigma$} & & \multicolumn{1}{c}{$Z_D$} & \multicolumn{1}{c}{$\varepsilon$} &  \multicolumn{1}{c}{$\bar{Z}_C$}  & \multicolumn{1}{c}{MIP} & \multicolumn{1}{c}{Row} &\multicolumn{1}{c}{R+C} & & \multicolumn{1}{c}{$Z_D$} & \multicolumn{1}{c}{$\varepsilon$} &  \multicolumn{1}{c}{$\bar{Z}_C$} & \multicolumn{1}{c}{MIP} & \multicolumn{1}{c}{Row} &\multicolumn{1}{c}{R+C}   \\
    \midrule
25	&	5	&	2	&	456	&	5.0	&&	124.6	&	5.7\%	&	131.7	&	1	&	2	&	2	&&	180.0	&	39.3\%	&	250.7	&	12	&	35	&	44	\\
25	&	5	&	2	&	1,864	&	2.5	&&	126.0	&	2.8\%	&	129.5	&	2	&	3	&	2	&&	182.4	&	19.4\%	&	217.8	&	56	&	45	&	42	\\
25	&	5	&	2	&	7,446	&	1.3	&&	126.4	&	1.4\%	&	128.2	&	9	&	9	&	3	&&	182.8	&	9.7\%	&	200.4	&	107	&	135	&	49	\\
25	&	5	&	2	&	29,621	&	0.6	&&	126.9	&	0.7\%	&	127.8	&	93	&	311	&	8	&&	183.4	&	4.8\%	&	192.2	&	*10,000	&	5,359	&	64	\\
    \midrule
50	&	10	&	5	&	886	&	5.0	&&	105.2	&	6.7\%	&	112.3	&	6	&	12	&	4	&&	178.0	&	39.7\%	&	248.7	&	37	&	342	&	230	\\
50	&	10	&	5	&	3,563	&	2.5	&&	107.4	&	3.3\%	&	110.9	&	94	&	104	&	13	&&	181.8	&	19.4\%	&	217.2	&	390	&	1,257	&	290	\\
50	&	10	&	5	&	14,242	&	1.3	&&	108.3	&	1.6\%	&	110.0	&	*10,000	&	*10,000	&	33	&&	183.1	&	9.7\%	&	200.8	&	*10,000	&	*10,000	&	414	\\
50	&	10	&	5	&	57,016	&	0.6	&&	108.5	&	0.8\%	&	109.4	&	*10,000	&	*10,000	&	49	&&	181.3	&	4.9\%	&	190.2	&	*10,000	&	*10,000	&	603	\\
        \midrule
75	&	15	&	10	&	1,263	&	5.0	&&	83.6	&	8.5\%	&	90.7	&	23	&	1,875	&	174	&&	139.3	&	50.8\%	&	210.0	&	142	&	293	&	340	\\
75	&	15	&	10	&	5,112	&	2.5	&&	85.1	&	4.2\%	&	88.6	&	2,415	&	*10,000	&	523	&&	138.9	&	25.5\%	&	174.2	&	2,088	&	9,799	&	614	\\
75	&	15	&	10	&	20,470	&	1.3	&&	85.4	&	2.1\%	&	87.1	&	*10,000	&	*10,000	&	614	&&	139.1	&	12.7\%	&	156.8	&	*10,000	&	*10,000	&	1,534	\\
75	&	15	&	10	&	81,940	&	0.6	&&	85.8	&	1.0\%	&	86.7	&	*10,000	&	*10,000	&	778	&&	140.2	&	6.3\%	&	149.0	&	*10,000	&	*10,000	&	1,610	\\
            \midrule
100	&	20	&	20	&	1,511	&	5.0	&&	60.0	&	11.8\%	&	67.0	&	30	&	228	&	53	&&	84.6	&	83.6\%	&	155.3	&	9	&	174	&	87	\\
100	&	20	&	20	&	6,158	&	2.5	&&	61.4	&	5.8\%	&	64.9	&	2,617	&	2,579	&	485	&&	84.7	&	41.7\%	&	120.1	&	129	&	454	&	172	\\
100	&	20	&	20	&	24,743	&	1.3	&&	62.1	&	2.8\%	&	63.9	&	*10,000	&	*10,000	&	585	&&	85.7	&	20.6\%	&	103.4	&	*10,000	&	4,255	&	2,175	\\
100	&	20	&	20	&	99,032	&	0.6	&&	62.7	&	1.4\%	&	63.5	&	*10,000	&	*10,000	&	1,739	&&	73.0	&	12.1\%	&	81.8	&	*10,000	&	*10,000	&	3,252	\\
    \bottomrule
    \end{tabular}
\end{table}
\end{landscape}

\section{AED Deployment Study: Setup}\label{sec:casestudy}
In this section, we discuss the setup of our numerical study on the placement of AEDs in public locations, including an overview of the data used, calibration of the uncertainty set, benchmark models, and a simulation-based validation method to evaluate the performance of the AED deployments produced by each model. 

\subsection{Data}\label{section:experimentaldesign}

{\bf Historical cardiac arrests.} We consider the AED location problem in a densely populated region in the downtown core of Toronto, Canada, which roughly corresponds to the city's financial district. We obtained a cardiac arrest dataset through the Resuscitation Outcomes Consortium (ROC), which is a collaborative network of research institutions across the United States and Canada that collects cardiac arrest data for the purpose of improving clinical practice \citep{morrison2008}. The ROC dataset contains the geographic coordinates of all cardiac arrests in the City of Toronto from January 2006 to April 2013. During this period, a total of 43 cardiac arrests occurred in public locations within the study region. These cardiac arrest locations are used to estimate the $\lambda_j$ parameters in formulation \ref{model:robustcvar2}.

\begin{figure}[htbp]
  \centering
\includegraphics[scale = 0.60, clip = true, trim = 35mm 10mm 10mm 65mm]{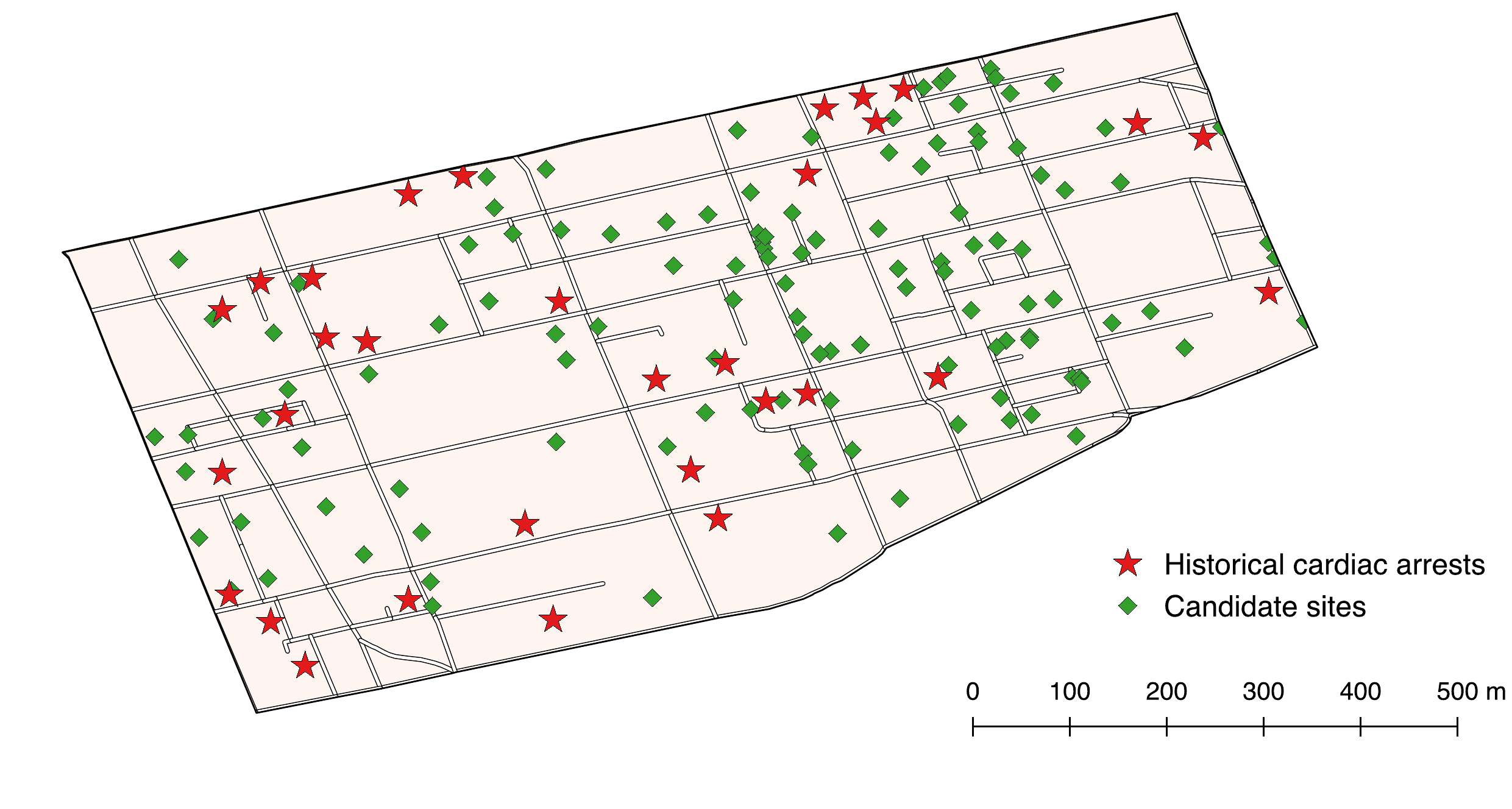}
   \caption{Locations of historical cardiac arrests and candidate sites for AED placement.}
\label{fig:CAMap}
\end{figure}

{\bf Candidate sites.} We constructed a set of candidate sites for AED placement by randomly selecting 120 public spaces (restaurants, cafes, shops, etc.) within the downtown region we consider. The addresses of these public spaces were obtained through the City of Toronto Employment Survey \citep{employmentsurvey2010}, which collects data on all Toronto businesses on an annual basis. We use publicly accessible buildings such as restaurants and cafes since they represent the type of location where one might reasonably expect an AED to be accessible to the public through a public access defibrillation program. Figure \ref{fig:CAMap} provides a visualization of the study region, historical cardiac arrests and the candidate sites for AED placement.

{\bf Distance parameters.} To calculate the distance parameters, we converted the addresses of the candidate sites into Universal Transverse Mercator (UTM) coordinates. Coordinates in the UTM system are given in units of meters (unlike latitude and longitude coordinates), which allows for straightforward calculation of distances within the study region. We constructed the scenario set $\mathcal{K}$ by intersecting the study region with a square grid with a spacing of 20 meters, which resulted in approximately 2,600 scenarios. We generated the $d_{ik}$ parameters by calculating the Euclidean distance between each candidate site and scenario pair. Although the Euclidean distance is only an approximation of the one-way distance a lay responder would have to travel to retrieve an AED, previous studies have shown Euclidean distance to be highly correlated with actual road distance, especially in urban areas where the road network is relatively dense \citep{bach1981problem,jones2010spatial,boscoe2012nationwide}. Moreover, the use of Euclidean distances is particularly justified in the context of AED placement, since the typical distance between a cardiac arrest and the nearest AED is low enough that lay responders can often take shortcuts that do not necessarily appear on the road network (e.g., through buildings, unmarked pedestrian paths, or parking lots).

As a robustness check, we repeat a subset of experiments using rectilinear distances and actual walking distances. The walking distances are obtained using the Google Maps Distance Matrix Application Programming Interface (API), implemented in Python. The Distance Matrix API is a query tool that accepts pairs of locations as input (as latitude-longitude coordinates or street addresses) and returns the distance between the two locations for a given mode of transportation (e.g., walking or driving) \citep{Maps}.

\subsection{Uncertainty regions and estimation of $\lambda_j$} \label{sec:MLE}
One of the key modeling decisions in our approach is to construct a set of uncertainty regions from the full service area. One simple approach is to partition the service area into a set of identically shaped cells and allow each cell to represent an uncertainty region. We propose that in the context of public access defibrillation, it may be most intuitive from a policy making perspective to use pre-existing geographic divisions within the city, such as city blocks, postal/zip codes, or census tracts. Further, in settings where historical cardiac arrest data is unavailable, one may have to rely on alternate demographic information that are only recorded in aggregate as a proxy for cardiac arrest incidence. For example, the City of Toronto collects socio-economic data at the neighbourhood level for planning purposes \citep{nhood}, which naturally implies a set of uncertainty regions.  

Once a set of uncertainty regions is selected, the corresponding probabilities $\lambda_1,\ldots,\lambda_{|\mathcal{J}|}$ must be estimated from historical cardiac arrest data. We propose the following approach, which produces consistent (i.e., asymptotically optimal) estimates of the uncertainty region probabilities.
\begin{proposition}\label{prop:MLE}
 Let $a_1,\ldots,a_{n} \in \mathcal{A}$ be the locations of $n$ cardiac arrest events, and assume that $a_1,\ldots,a_{n}$ are drawn i.i.d. from the distribution $\mu$. For each $j \in \mathcal{J}$, let
\begin{align} \label{MLE}
\hat{\lambda}_j^{n} =  \frac{1}{n} \sum_{c = 1}^{n} \mathbbm{1}\{a_c \in \mathcal{A}_j \}.
\end{align}
Then $|\hat{\lambda}_j^n - \lambda_j| \longrightarrow 0$ as $n \longrightarrow \infty$ almost surely. Further, if the sets $\Xi'_1,\ldots,\Xi'_{|\mathcal{R}|}$ are all non-empty, then the uncertainty set  ${\bf U}$ constructed from $\hat{\lambda}_j^n$, $j \in \mathcal{J}$ is guaranteed to be non-empty.
\end{proposition}
The estimator in Proposition \ref{prop:MLE} is simply the fraction of cardiac arrest events that have occurred in each uncertainty region. Since the uncertainty regions are fixed, the convergence result in Proposition \ref{prop:MLE} follows directly from the strong law of large numbers (to be precise, note that the arrival of cardiac arrests in each region can, with some modification, be modeled as a multinomial trials process, which has a well-known maximum likelihood estimator \citep{bickel2015mathematical}). The assumption that the cardiac arrest locations are drawn i.i.d. from a common distribution is a standard statistical assumption, and can be interpreted as requiring the probability of a cardiac arrest event in different parts of the city to be stable over time. Fortunately, cardiac arrests in Toronto have been observed to exhibit this temporal stability \citep{chan2014}.

 Note also that these properties hold for {\it any} chosen configuration of the uncertainty regions. In other words, once the uncertainty regions are selected, the estimate $\hat{\lambda}^n$ is guaranteed to converge to the true arrival probabilities. This convergence is generally known as statistical {\it consistency}, which is considered a fundamental requirement for an estimator \citep{bickel2015mathematical}. 

In addition to estimating the arrival probabilities, we can also use the Wilson score method as a heuristic to estimate an associated 95\% confidence interval for each arrival probability \citep{wilson1927probable}. These confidence intervals can also be used to guide the design of the uncertainty regions. For example, we might require that the uncertainty regions be selected such that estimation error of each of the arrival probabilities $\lambda_1,\ldots,\lambda_{|\mathcal{J}|}$ is no greater than some threshold $\delta$ with $95\%$ confidence. We emphasize, however, that the Wilson score method is only a heuristic for estimating the confidence intervals, meaning the estimation error is not always guaranteed to be less than $\delta$ with $95\%$ confidence. As an alternative, we also note that the problem of designing uncertainty regions based on data is closely related to the well-studied problem in statistics of optimally selecting histogram bin widths (e.g., \cite{wand1997data}), although we do not formally explore that connection in this paper.

The service area in our numerical study spans only two census tracts and three postal codes zones. Since using too few uncertainty regions may fail to capture important aspects of the cardiac arrest distribution, we instead identified uncertainty regions by using the road network to divide the study area into 15 approximately equal-area uncertainty regions. The estimates of the arrival probabilities for the 15 uncertainty regions obtained using this approach are reported in Section \ref{estimateprobabilities} of the electronic companion, along with their approximate 95\% confidence intervals, calculated using the Wilson score method \citep{wilson1927probable}.

\subsection{Benchmarks: Sample average approximation and ex-post models}
We compare the performance of the AED deployment produced by the robust formulation \ref{model:robustcvar2} with two benchmark models. The first benchmark is a {\it nominal} model, which optimizes only over historical cardiac arrest locations without taking uncertainty into account. The nominal model represents the the sample average approximation approach to the AED problem, since it optimizes directly with respect to a historical sample. The nominal model is given by formulation \ref{model:muknownprogram}, where the scenarios $\{\xi_1,\ldots,\xi_{|\mathcal{K}|}\}$ are given by the historical cardiac arrest locations, and $u_k = 1/|\mathcal{K}|$ for all $k \in \mathcal{K}$. The second benchmark is an {\it ex-post} model, which optimizes directly over a set of simulated cardiac arrests, and thus has perfect foresight. The ex-post model is also given by formulation \ref{model:muknownprogram}, except instead of optimizing over historical cardiac arrest locations, the ex-post model optimizes the AED deployment with respect to a simulated set of cardiac arrest locations, which we also evaluate the performace of the nominal and robust models against. The ex-post model can be thus interpreted as the ``best possible" model we could solve if the locations of all future cardiac arrests were known a priori. In all instances, we set $\beta = 0.9$ and the number of AEDs ($P$) to 30.  All instances for all models were implemented using MATLAB R2011a via YALMIP, and solved using CPLEX 12.1 with default parameter settings on a single node of a computing cluster with a 2.9 GHz quad-core CPU. 

\subsection{Model validation}
To assess the performance of the AED deployments produced by \ref{model:robustcvar2} and \ref{model:muknownprogram}, we used simulation to generate hypothetical cardiac arrest events, and computed various performance metrics associated with the distance distribution induced by each AED deployment and the set of simulated cardiac arrests. We then solve the ex-post model based on the simulated cardiac arrest locations, to gain insight into the best-possible deployment, for the given model parameterization and objective function.

To construct the validation sets, we first used kernel density estimation to estimate the underlying demand density, which is a well known semi-parametric technique for estimating a probability density function from a finite sample \citep{sheather1991,terrell1992}. Intuitively, kernel density involves centering a continuous density function (the kernel) at each data point, and then aggregating and normalizing all kernel functions to obtain a single probability density function. Kernel density estimation requires the specification of a kernel function as well as a parameter $h$, known as the {\it bandwidth}, which is typically proportional to the standard deviation of the kernel distribution. In short, the bandwidth represents the degree to which the data set is smoothed to obtain the estimated density, with a higher bandwidth resulting in more aggressive smoothing. We note that using kernel density estimation to estimate the distribution of future cardiac arrests is supported by the fact that cardiac arrest locations in Toronto have been shown to exhibit temporal stability \citep{chan2014}.

\begin{figure}
\centering
\subfigure[h = 50m.]{%
\includegraphics[scale=0.35, clip = true, trim = 40mm 15mm 40mm 65mm]{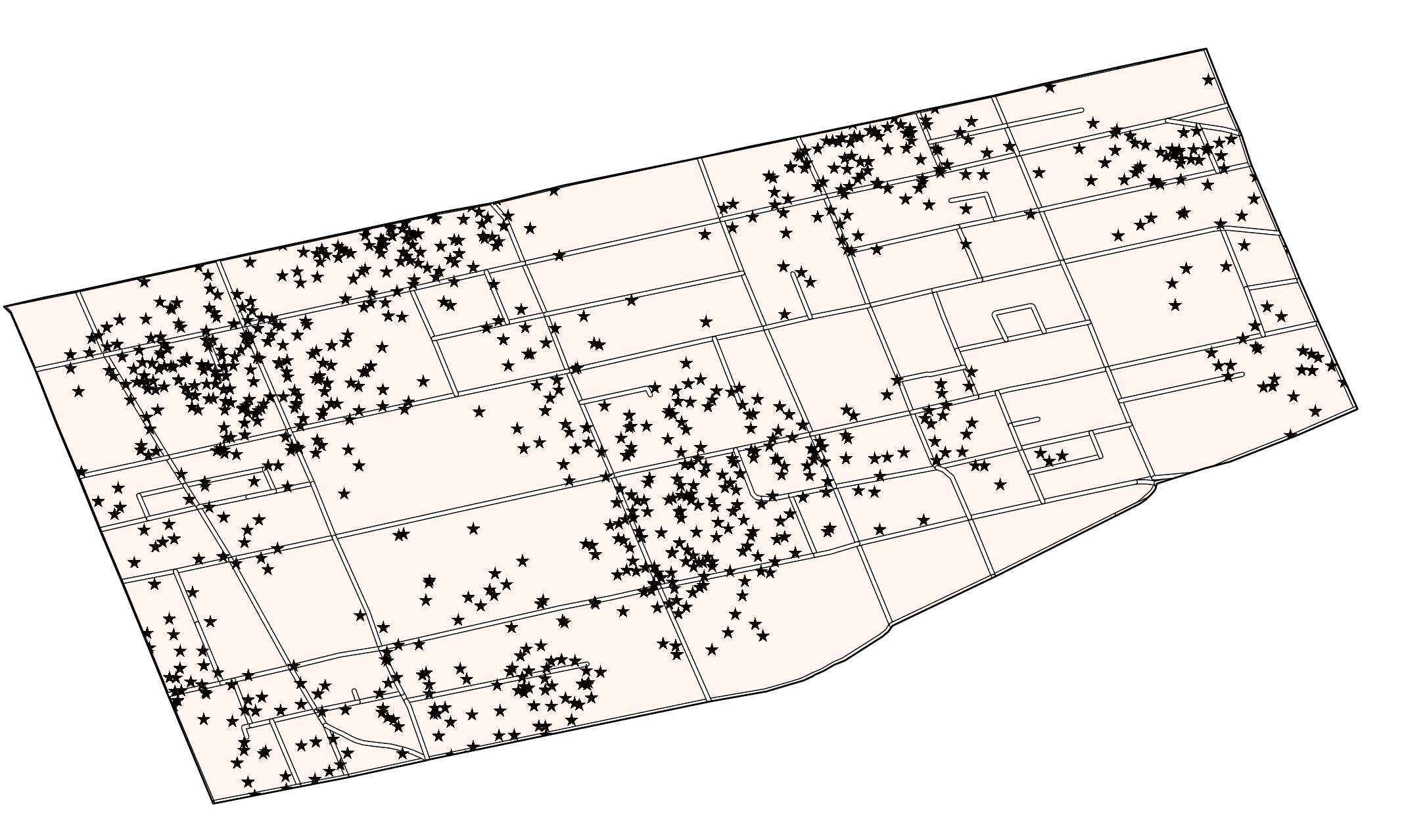}}
\label{fig:subfigure1}
\hspace{5mm}
\subfigure[h = 100m.]{%
\includegraphics[scale=0.35, clip = true, trim = 40mm 15mm 40mm 65mm]{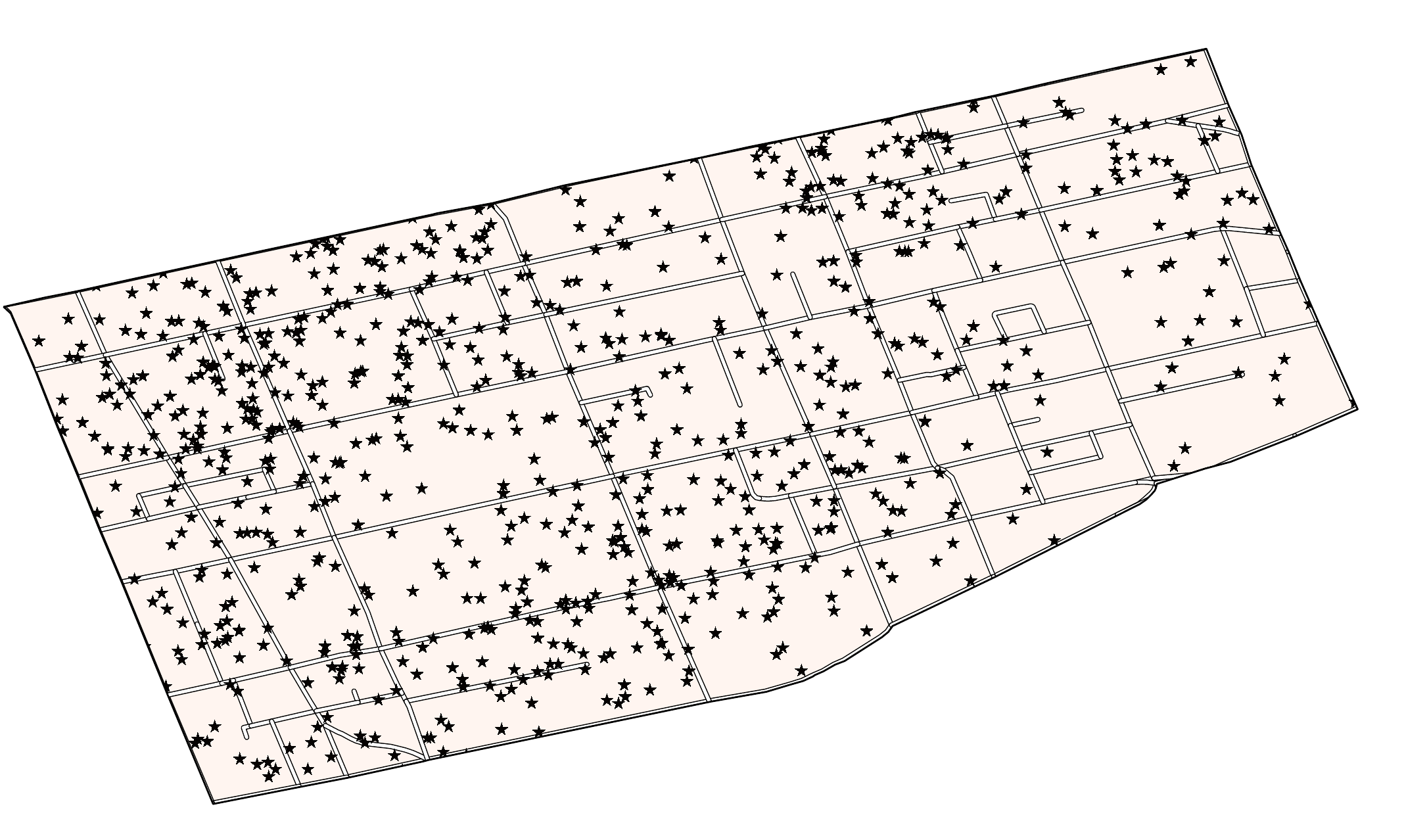}
\label{fig:subfigure2}}
\vspace{2mm}
\caption{Locations of 1,000 simulated cardiac arrests based on kernel density estimation using bandwidths of 50m and 100m.}
\label{fig:KDE}
\end{figure}

To assess model performance under different values of the bandwidth parameter $h$, we constructed four different validation sets. Specifically, we used Gaussian kernels with $h = 10$, 50, 100 and 150 meters to estimate four demand densities, each of which is used to simulate a separate validation set. The case where $h = 10$ can be interpretted as a ``low uncertainty" environment, since the low bandwidth causes a majority of simulated events to fall near a historical cardiac arrest location. Conversely, the case where $h = 150$ can be interpretted as a ``high uncertainty" environment, since the larger bandwidth leads to some cardiac arrests being simulated far from any historical location. By varying the bandwidth parameter from $10$ to $150$, we are able to evaluate the performance of the AED deployment under a variety of possible demand distributions. For each of the four estimated densities, we simulated a total of 50 validation sets containing 100 cardiac arrests each. We do not impose any constraints on where cardiac arrests may be simulated (i.e., a small fraction may occur beyond the service area). Figure \ref{fig:KDE} illustrates the effect of the bandwidth on the estimated demand density. Note that the distribution of the simulated cardiac arrests is notably more uniform for a bandwidth of 100 meters compared to the distribution for a bandwidth of 50 meters, due to increased smoothing.

 \section{AED Deployment Study: Results}\label{sec:results}

The numerical results are organized as follows. In Section \ref{sec:performancemetrics}, we compute the average value (over 50 validation sets) of four performance metrics: the mean, VaR, CVaR and maximum of the distance distribution induced by each AED deployment. We also report worst-case performance of each deployment, which is given by the maximum value of each performance metric over 50 validation sets. In Section \ref{section:coverage}, we evaluate the models using a different but related measure of AED deployment performance, known as cardiac arrest coverage. Lastly, in Section \ref{sec:alt} we do a robustness check by presenting results under alternate distance measures, namely, actual walking distances obtained via Google Maps and the rectilinear distance metric.

\subsection{Performance metrics}\label{sec:performancemetrics} 
\begin{figure}
\centering
\subfigure[Mean.]{%
\includegraphics[scale=0.4, clip = true, trim = 15mm 65mm 20mm 70mm]{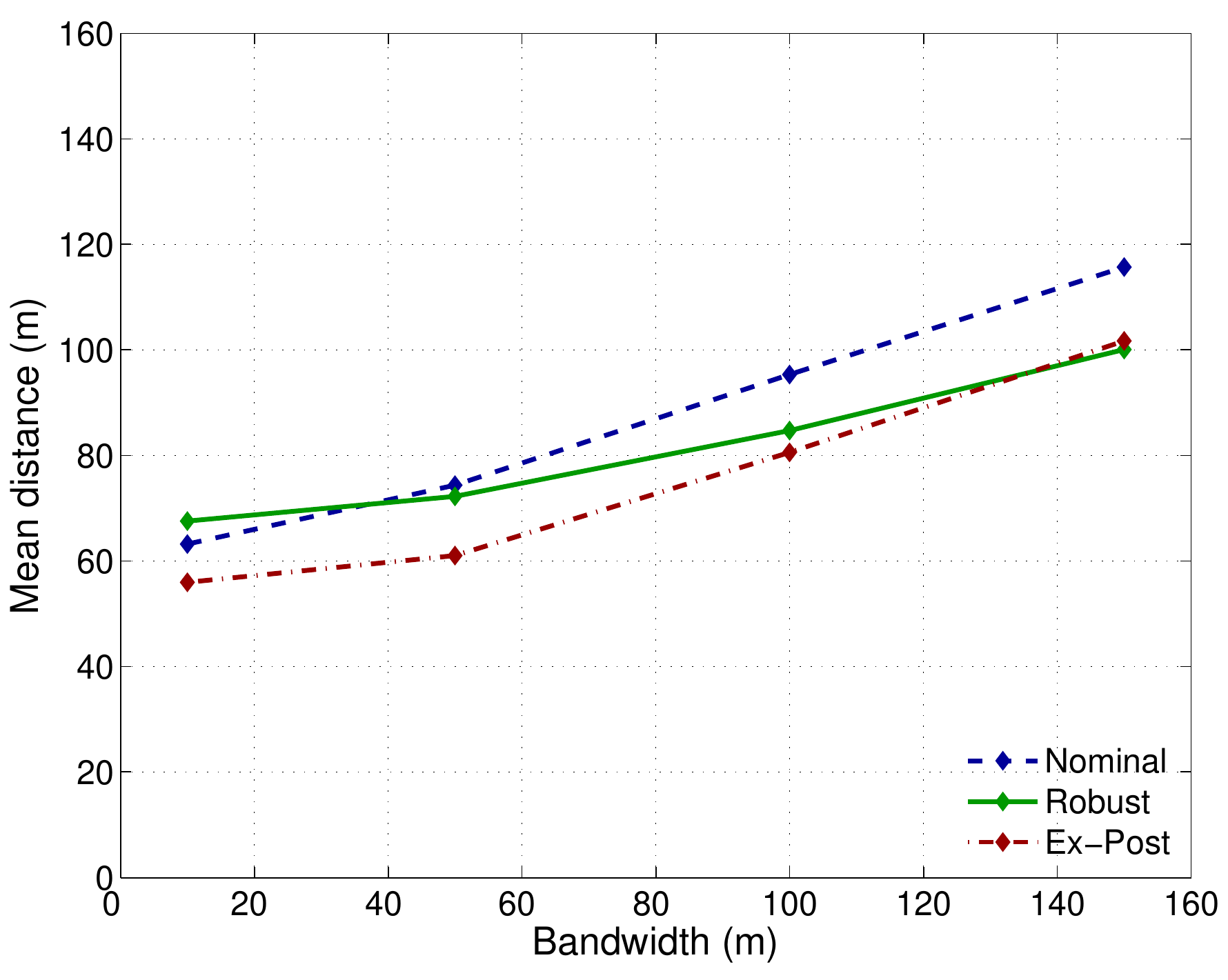}
\label{fig:subfigure1}}
\subfigure[VaR.]{%
\includegraphics[scale=0.4, clip = true, trim = 15mm 65mm 20mm 70mm]{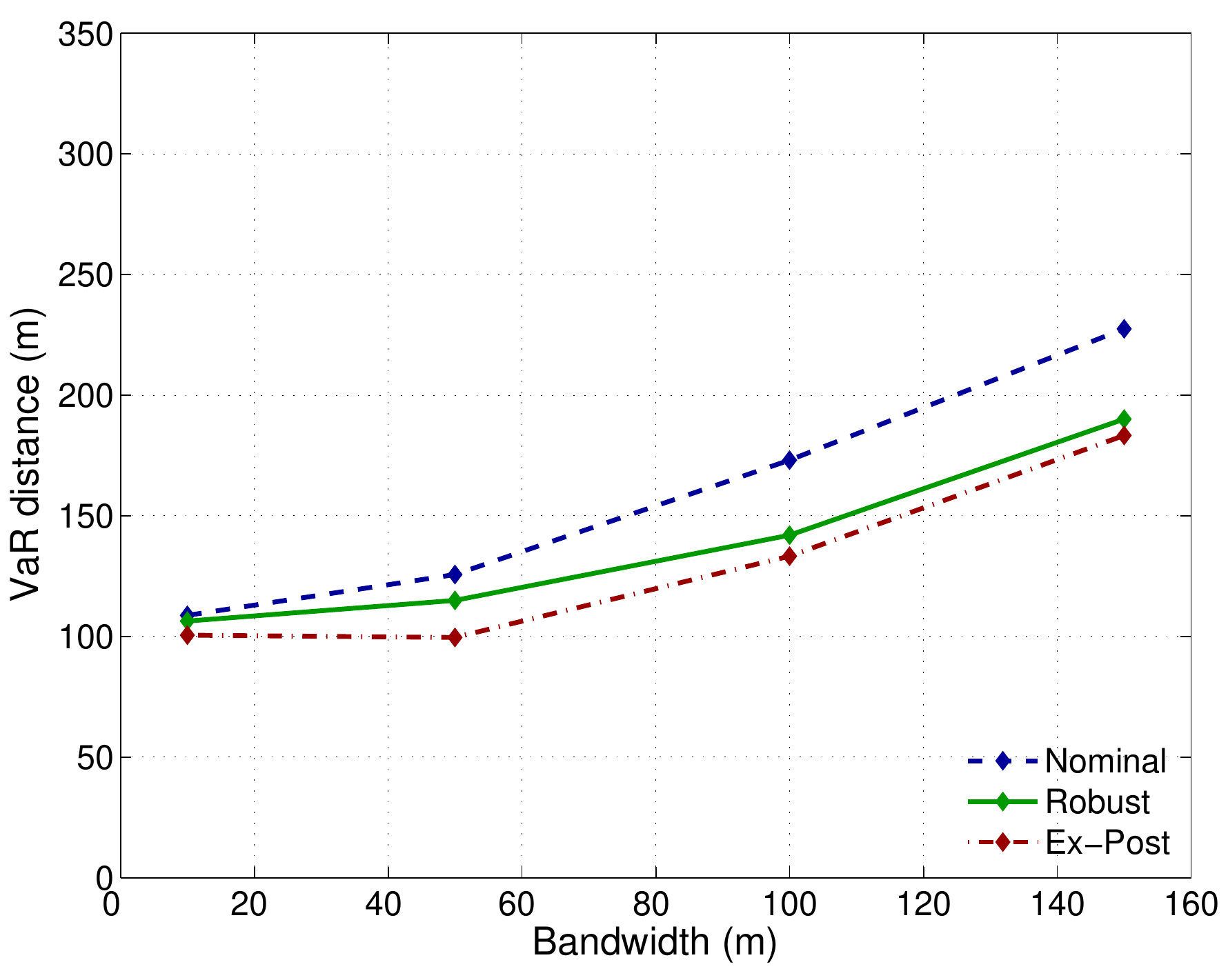}
\label{fig:subfigure2}}
\vspace{10mm}
\subfigure[CVaR.]{%
\includegraphics[scale=0.4, clip = true, trim = 15mm 65mm 20mm 70mm]{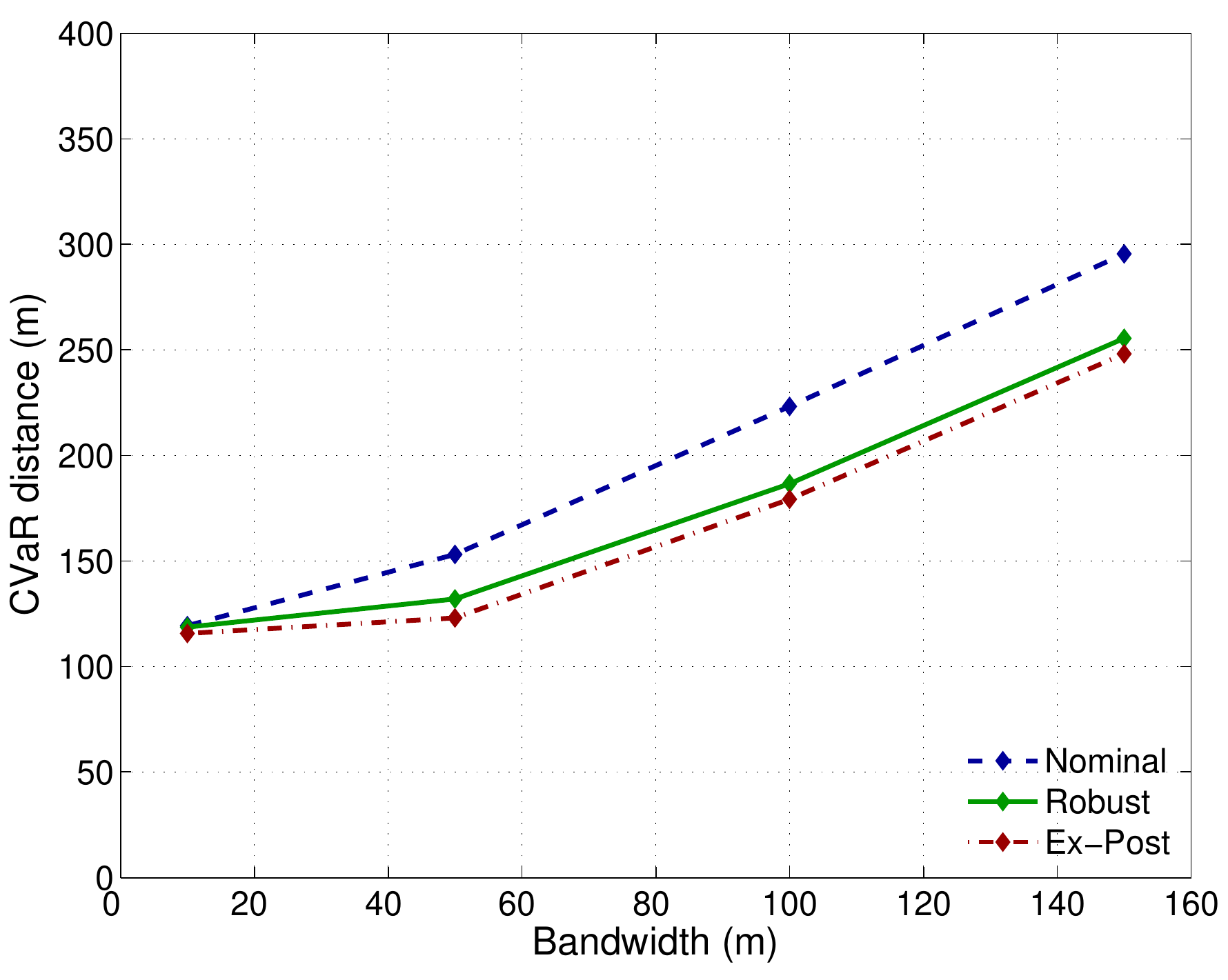}
\label{fig:averagecvar}}
\subfigure[Maximum.]{%
\includegraphics[scale=0.4, clip = true, trim = 15mm 65mm 20mm 70mm]{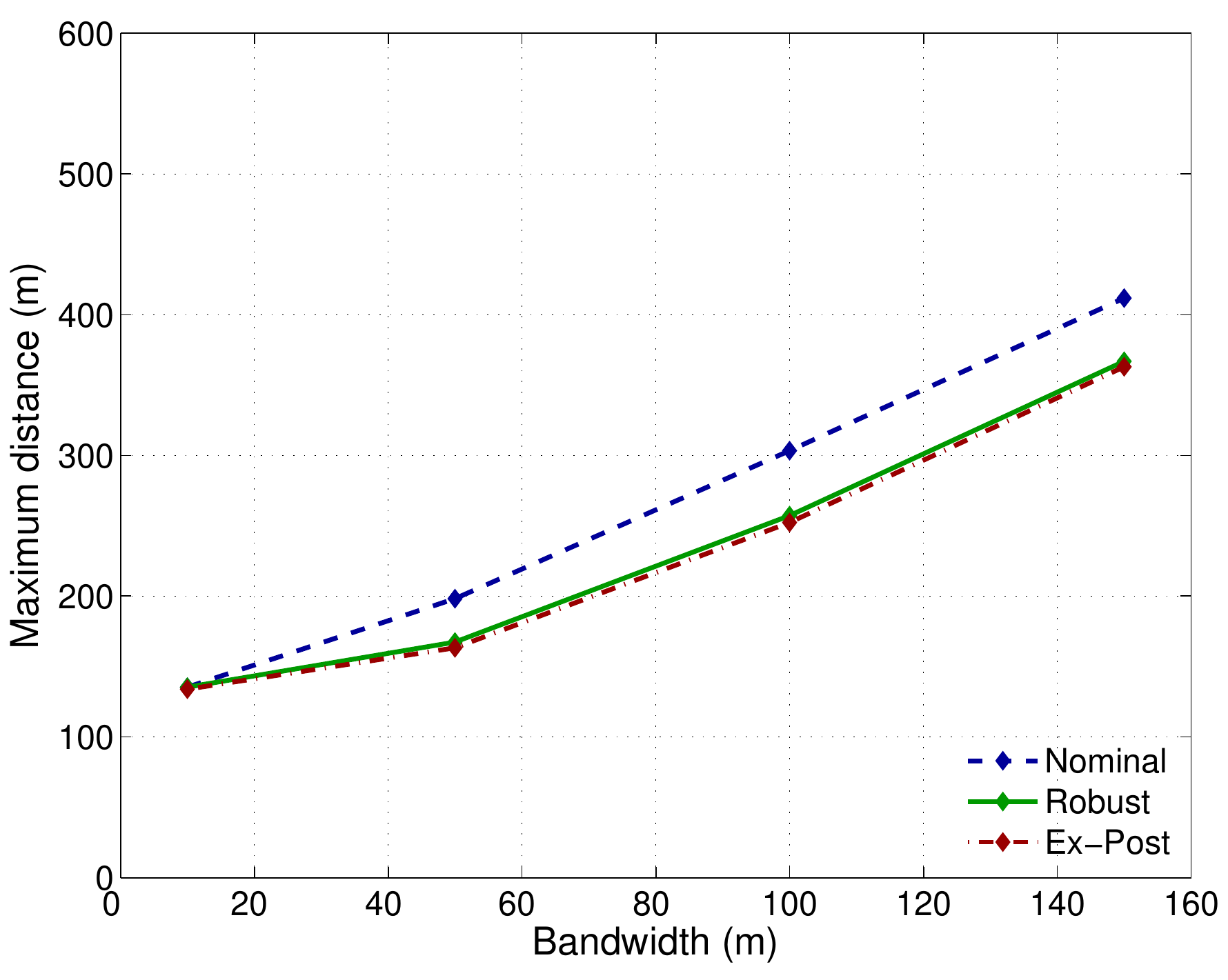}
\label{fig:subfigure2}}
\vspace{-6mm}
\caption{Average performance over 50 validation sets of nominal, robust and ex-post solutions under increasing kernel bandwidth ($\beta = 0.9$).}
\label{fig:average}
\end{figure}

Figure \ref{fig:average} shows the mean, VaR, CVaR and maximum value of the distance distribution for the nominal, robust and ex-post models (note that the y-axis differs between the plots).  The performance metrics are computed for each of the four demand distributions (i.e., four bandwidth parameters), and are averaged over 50 validation sets. Since Figure \ref{fig:average} shows the average value of each performance metric, we can interpret it as depicting the performance of each model under ``typical" cardiac arrest locations.

In the cases where $h = 10$, the AED deployments produced by all three models perform similarly. This is not a surprising observation, since in the $h = 10$ case, which represents a low uncertainty environment, the vast majority of simulated cardiac arrests occur within the vicinity of one of the 43 historical cardiac arrest locations. However, it is worth noting that the robust solution performs at least as well as the nominal solution in the $h = 10$ case. In other words, the performance of the robust model does not appear to suffer when the level of uncertainty is low, despite optimizing AED locations to protect against uncertainty. Moreover, in the cases where $h \ge 50$, the robust solution substantially outperforms the nominal solution on all four performance metrics. For example, in the case where $h = 100$, the robust solution outperforms the nominal solution on the mean distance metric by 9\% (86m vs 95m) and on the CVaR metric by 15\%  (223m vs 189m). The performance difference between the robust and nominal models was found to be statistically significant $(p < 0.05)$ in all comparisons where $h \ge 50$. Examples of AED deployments under each model for $\beta = 0$ and $\beta = 0.9$ are given in Section \ref{additional} of the electronic companion.

Note that the ex-post model achieves the best performance in all instances, which is expected since it optimizes the AED deployment directly over the simulated cardiac arrest locations (note that the ex-post model is solved independently for each of the 50 trials, since it takes the simulated cardiac arrest locations as input). The gap in performance between the nominal and ex-post models may therefore be interpreted as the degradation in performance due to the uncertainty in cardiac arrest locations. Figure \ref{fig:average} shows that the robust deployment substantially closes this performance gap in many cases. For example, for the CVaR metric with $h = 100$, using the robust deployment decreases the performance gap with the ex-post model from 44 meters (223m vs 179m) down to 10 meters (189m vs 179m), representing a gap improvement of 77\%. In this setting, it may be convenient to think of the performance gap as being analogous to the usual notion of optimality gap. However, we emphasize that the ex-post solution represents a theoretical lower bound on performance only for the CVaR metric, since all models optimize for CVaR only. For the remaining three performance metrics, it is possible for the nominal or robust model to outperform the ex-post model by chance.

Table \ref{table:gapav} shows the performance gap achieved by the nominal and robust deployment with respect to the ex-post deployment. For example, the first row and column of Table \ref{table:gapav} indicates that the performance of the nominal solution on the mean distance metric in the $h = 10$ case is within $12\%$ of the performance achieved by the ex-post solution, whereas the gap from the robust solution is within $7\%$ (averaged over 50 trials). The difference in performance between the nominal and robust solutions is most pronounced at the tail of the distribution, which is unsurprising since the models optimize for CVaR. Note that for the CVaR and maximum distance metrics, the performance gap of the robust solutions is at most 8\%, whereas the gap of the nominal solution can be as high as 20\%.

\begin{table}  
  \caption{Average performance gap of nominal and robust solutions with respect to ex-post solution.}
\centering
    \begin{tabular}{rrrrrrrrr}
\toprule 
            &\multicolumn{2}{c}{Mean} & \multicolumn{2}{c}{VaR} & \multicolumn{2}{c}{CVaR} & \multicolumn{2}{c}{Max} \\
            \cmidrule(l){2-3} \cmidrule(l){4-5} \cmidrule(l){6-7} \cmidrule(l){8-9}
           $h$ & \multicolumn{1}{c}{{Nom}.} & \multicolumn{1}{c}{Rob.} & \multicolumn{1}{c}{Nom.} & \multicolumn{1}{c}{Rob.}   & \multicolumn{1}{c}{Nom.} & \multicolumn{1}{c}{Rob.}  &\multicolumn{1}{c}{Nom.}  & \multicolumn{1}{c}{Rob.} \\
           \midrule
 10  & 12\%    & 7\%     & 8\% & 5\%     & 3\% & 3\%     & 1\% & 1\% \\
  50  & 18\%    & 11\%   & 21\%  & 13\%   & 20\% & 8\%    & 18\%  & 2\% \\
  100 & 15\%    & 6\%  & 23\%   & 9\%  & 20\%  & 5\%   & 17\%   & 3\% \\
  150 & 12\%    & 1\% & 19\%    & 4\% & 16\%    & 4\%  & 12\%    & 3\% \\  
    \bottomrule
        \end{tabular}
         \label{table:gapav}
\end{table}

\begin{figure}
\centering
\subfigure[Mean.]{%
\includegraphics[scale=0.40, clip = true, trim = 15mm 65mm 20mm 70mm]{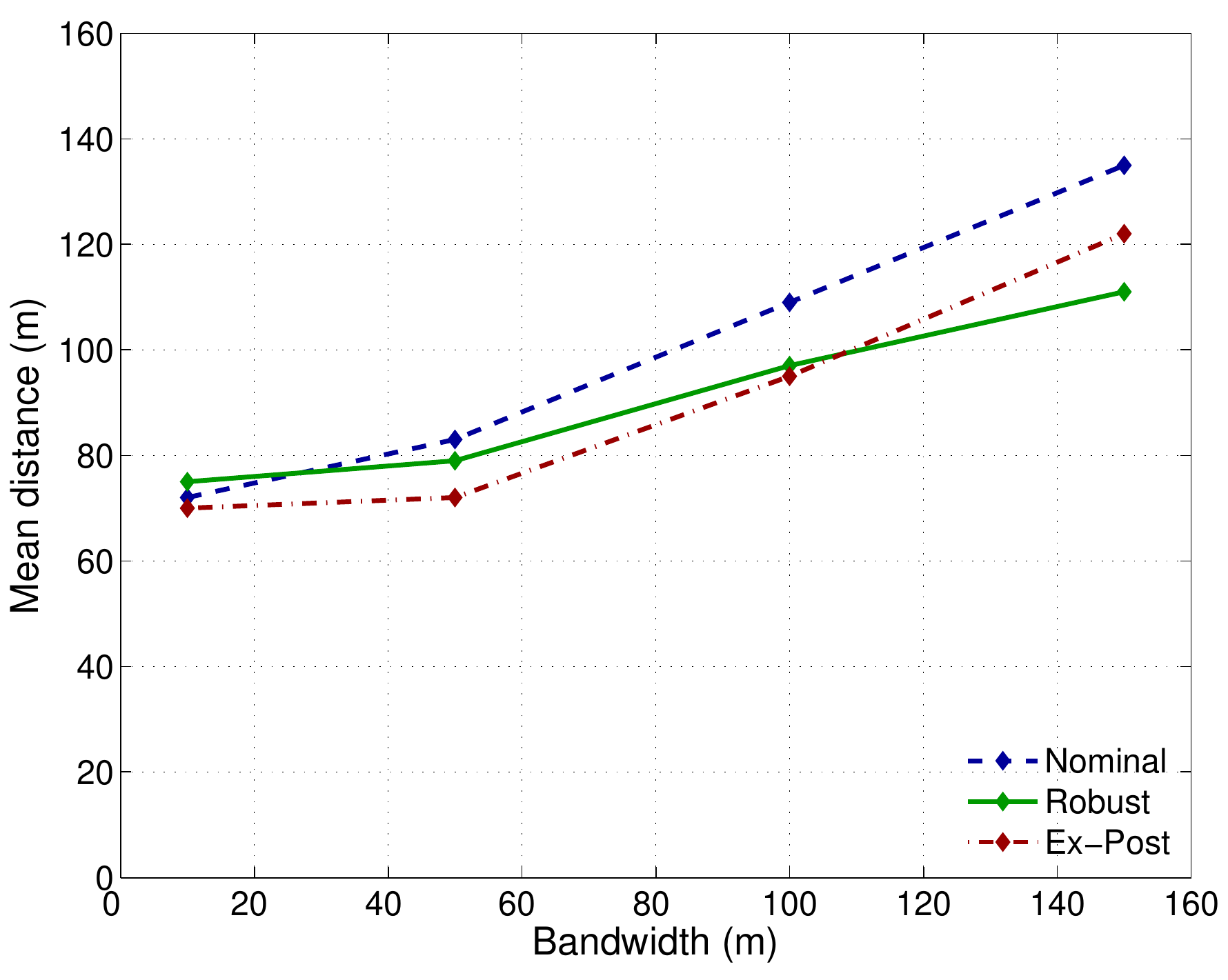}
\label{fig:subfigure1}}
\subfigure[VaR.]{%
\includegraphics[scale=0.40, clip = true, trim = 15mm 65mm 20mm 70mm]{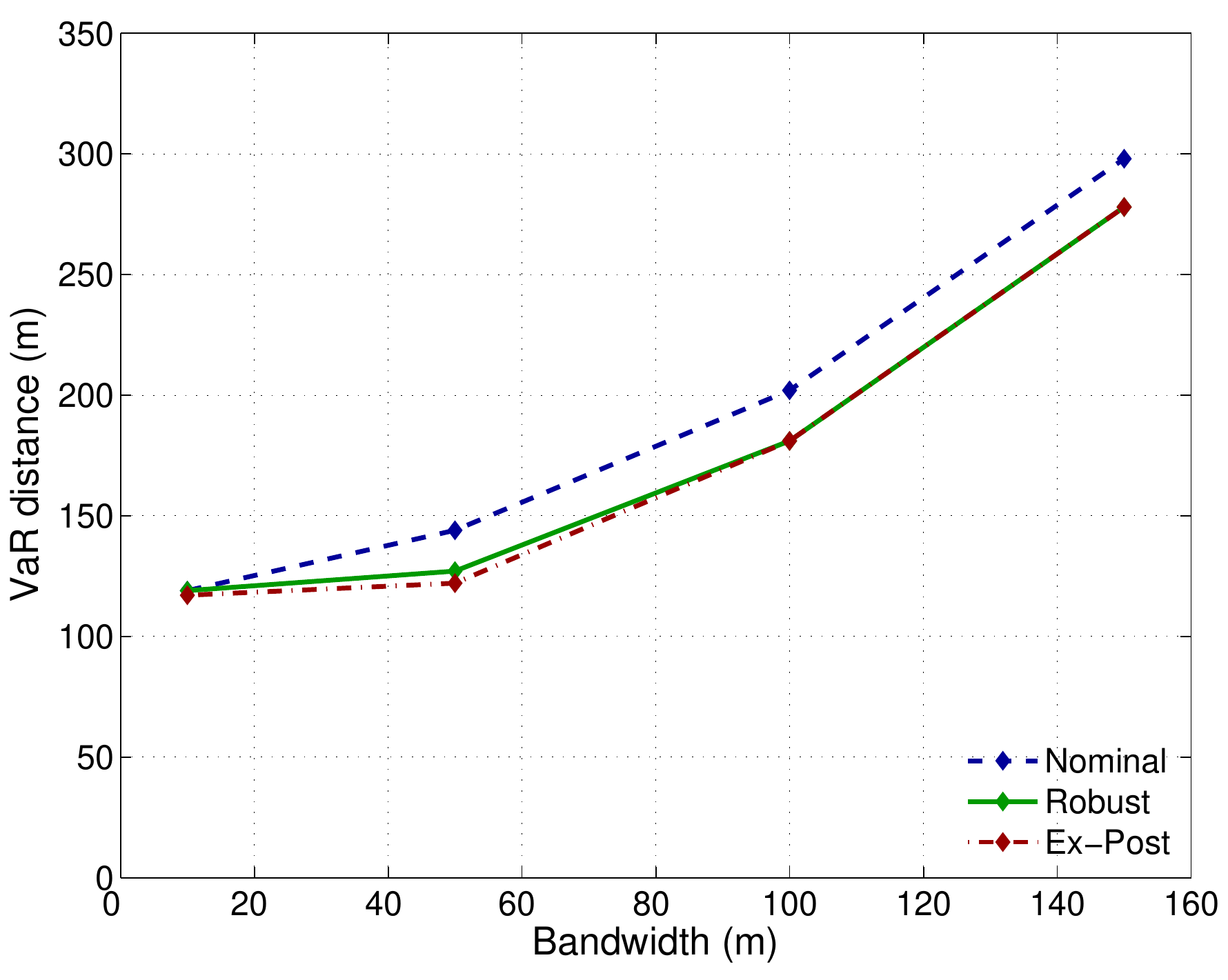}
\label{fig:subfigure2}}
\vspace{10mm}
\subfigure[CVaR.]{%
\includegraphics[scale=0.40, clip = true, trim = 15mm 65mm 20mm 70mm]{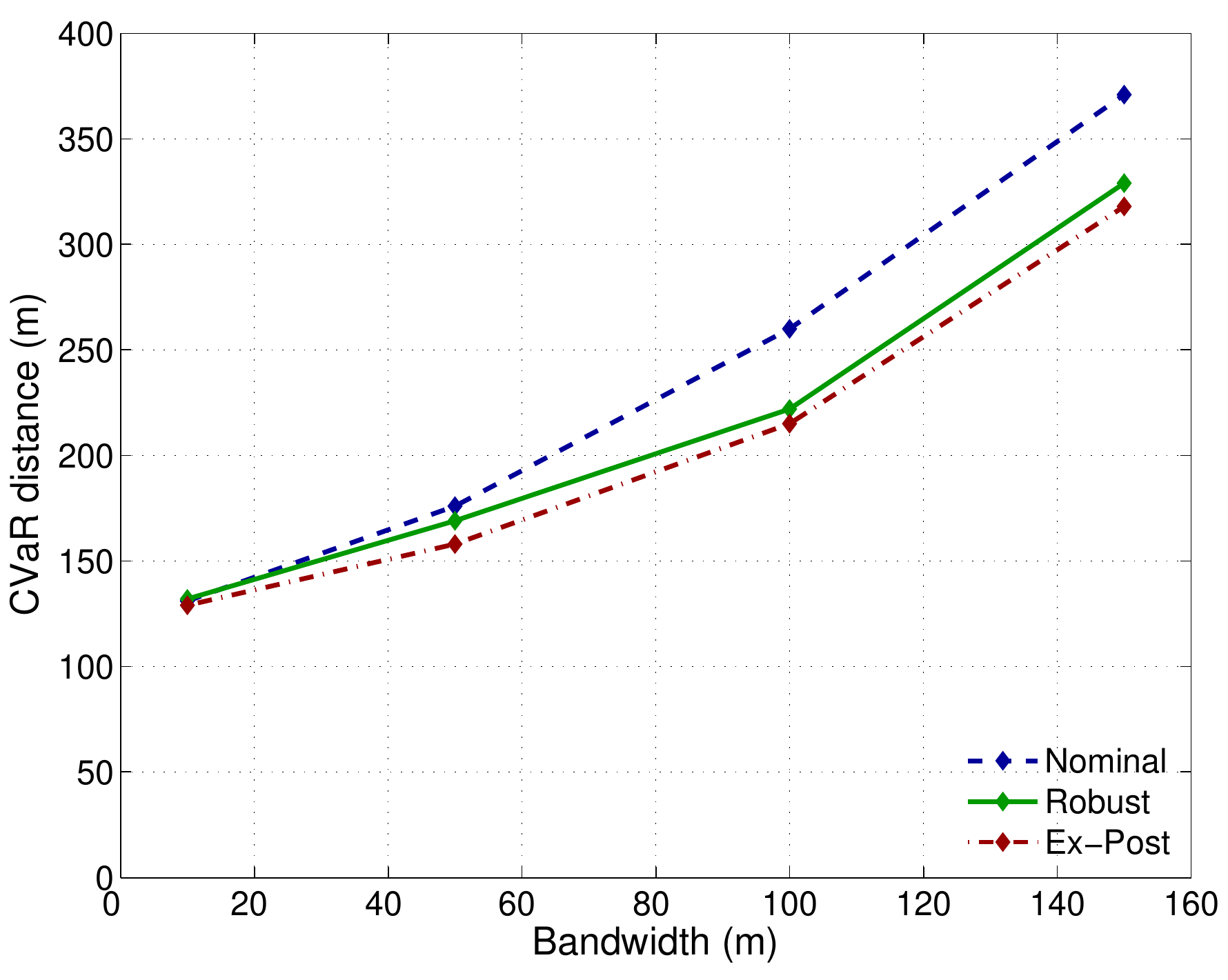}
\label{fig:subfigure1}}
\subfigure[Maximum.]{%
\includegraphics[scale=0.40, clip = true, trim = 15mm 65mm 20mm 70mm]{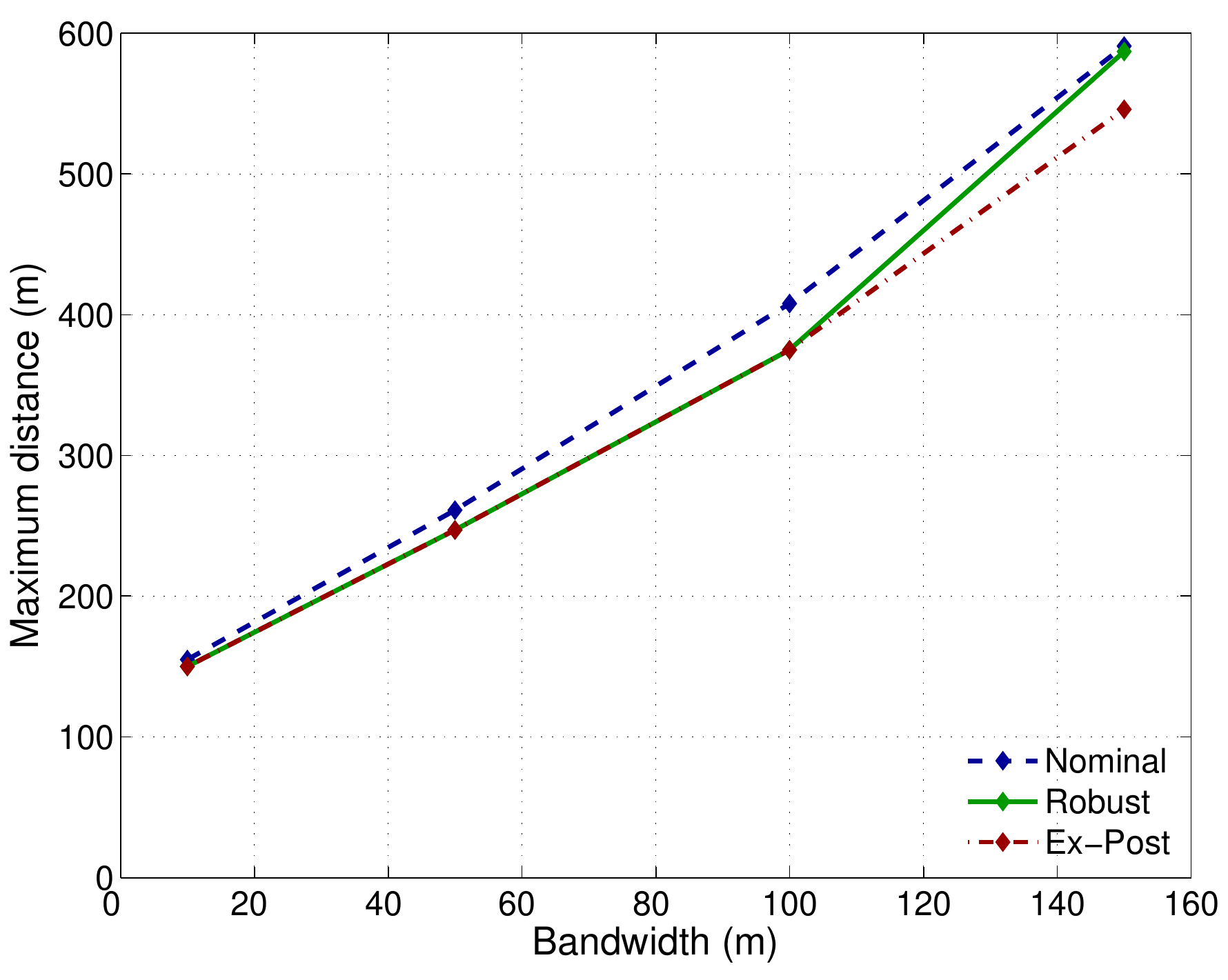}
\label{fig:subfigure2}}
\vspace{-6mm}
\caption{Worst-case performance of nominal, robust and ex-post solutions on four distance metrics under increasing kernel bandwidth.}
\label{fig:WC}
\end{figure}

We also consider the worst-case performance of each model by computing the maximum value of each performance metric over the 50 trials. As shown in Figure \ref{fig:WC}, the robust model generally outperforms the nominal model in the majority of high uncertainty instances ($h \ge 50$), and performs no worse than the nominal model in the low uncertainty instances ($h = 10$). Table \ref{table:gapwc} presents the worst-case performance gaps, which are obtained by computing the maximum relative gap in performance over 50 trials. As in Table \ref{table:gapav}, the robust solution substantially closes the worst-case performance gap in the vast majority of instances. For example, on the CVaR metric with $h = 100$, the worst-case performance gap is improved by 41\% (27\% down to 16\%).

\begin{table}  
\caption{Worst case performance gap of nominal and robust solutions with respect to ex-post solution.}
\centering
    \begin{tabular}{rrrrrrrrr}
\toprule 
            &\multicolumn{2}{c}{Mean} & \multicolumn{2}{c}{VaR} & \multicolumn{2}{c}{CVaR} & \multicolumn{2}{c}{Max} \\
            \cmidrule(l){2-3} \cmidrule(l){4-5} \cmidrule(l){6-7} \cmidrule(l){8-9}
           $h$ & \multicolumn{1}{c}{Nom.} & \multicolumn{1}{c}{Rob.} & \multicolumn{1}{c}{Nom.} & \multicolumn{1}{c}{Rob.}   & \multicolumn{1}{c}{Nom.} & \multicolumn{1}{c}{Rob.}  &\multicolumn{1}{c}{Nom.}  & \multicolumn{1}{c}{Rob.} \\
           \midrule
10 &30\%	&24\%	&28\%	    &19\%	 &10\%	&7\%	    &6\%	    &6\% \\
50 &29\%	&20\%	&34\%	    &29\%	 &31\%	&19\%	&41\%	&20\% \\
100 &29\%	&17\%	&40\%	    &26\%	 &27\%	&16\%	&31\%	&24\% \\
150 &20\%	&10\%	&31\%	    &20\%	 &24\%	&12\%	&26\%	&26\% \\
    \bottomrule
        \end{tabular}
         \label{table:gapwc}
\end{table}

\subsection{Cardiac arrest coverage}\label{section:coverage}
\begin{figure}
\centering
\subfigure[$h = 10$.]{%
\includegraphics[scale=0.40, clip = true, trim = 15mm 65mm 20mm 70mm]{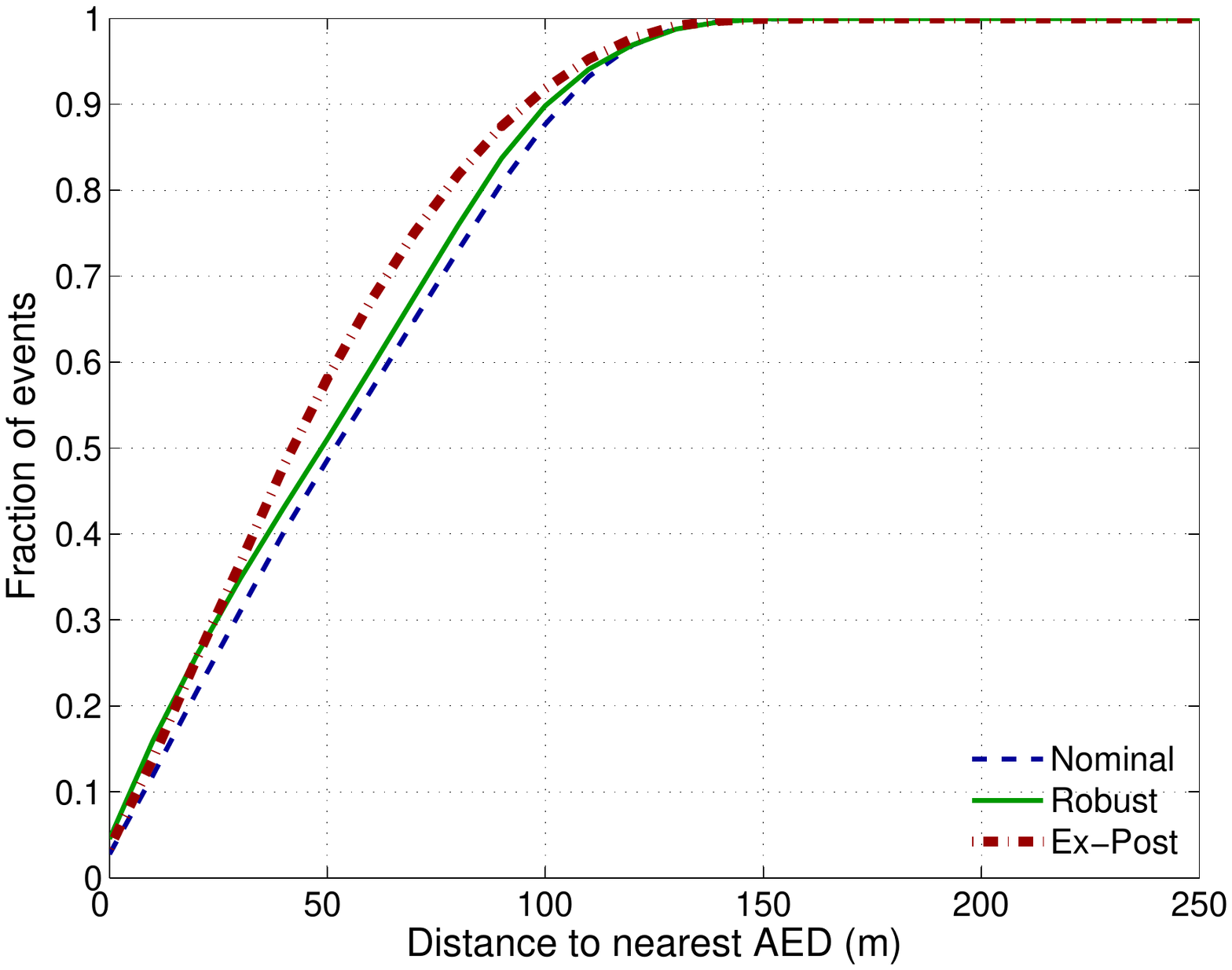}
\label{fig:subfigure1}}
\subfigure[$h = 50$.]{%
\includegraphics[scale=0.40, clip = true, trim = 15mm 65mm 20mm 70mm]{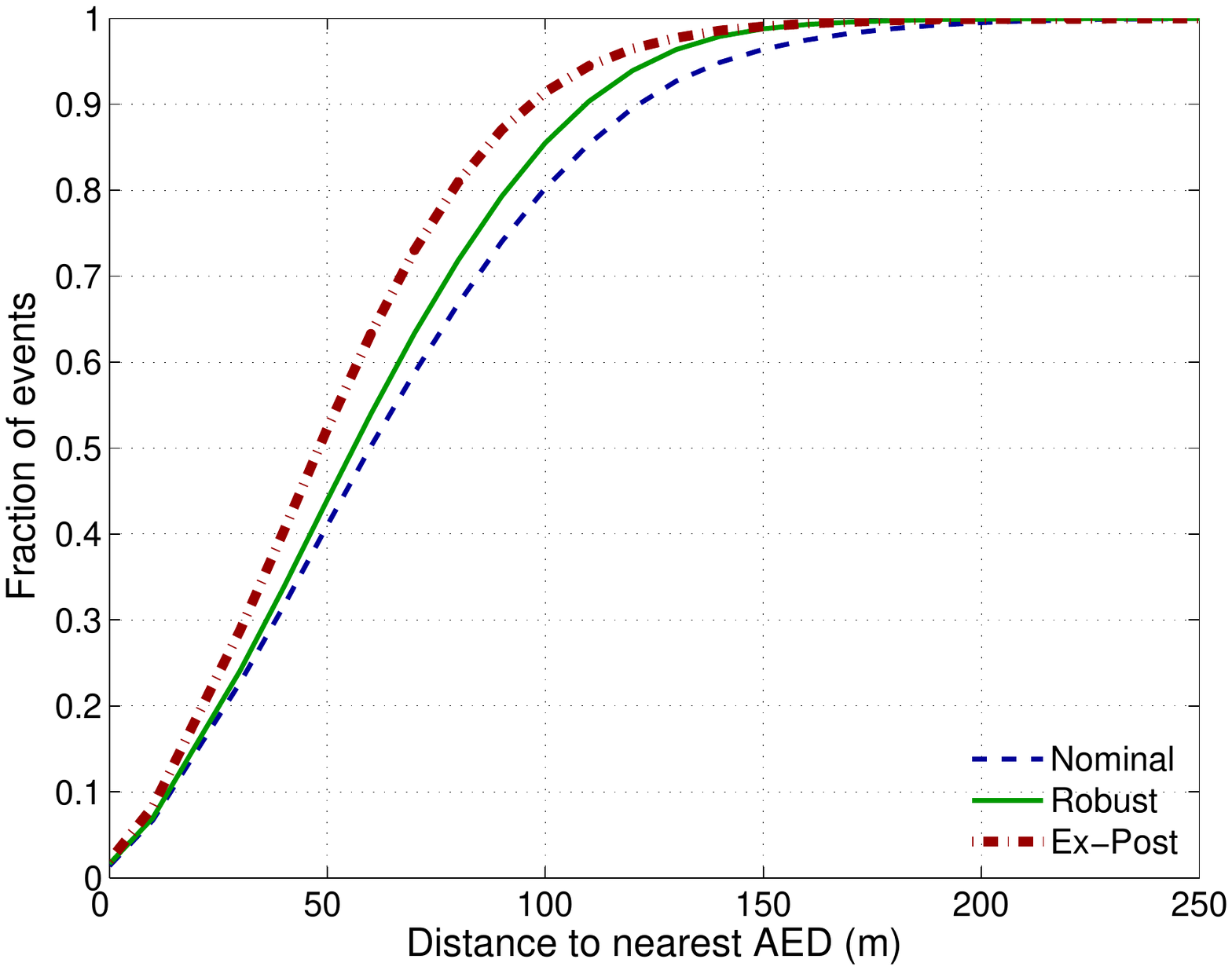}
\label{fig:subfigure2}}
\vspace{10mm}
\subfigure[$h = 100$.]{%
\includegraphics[scale=0.40, clip = true, trim = 15mm 65mm 20mm 70mm]{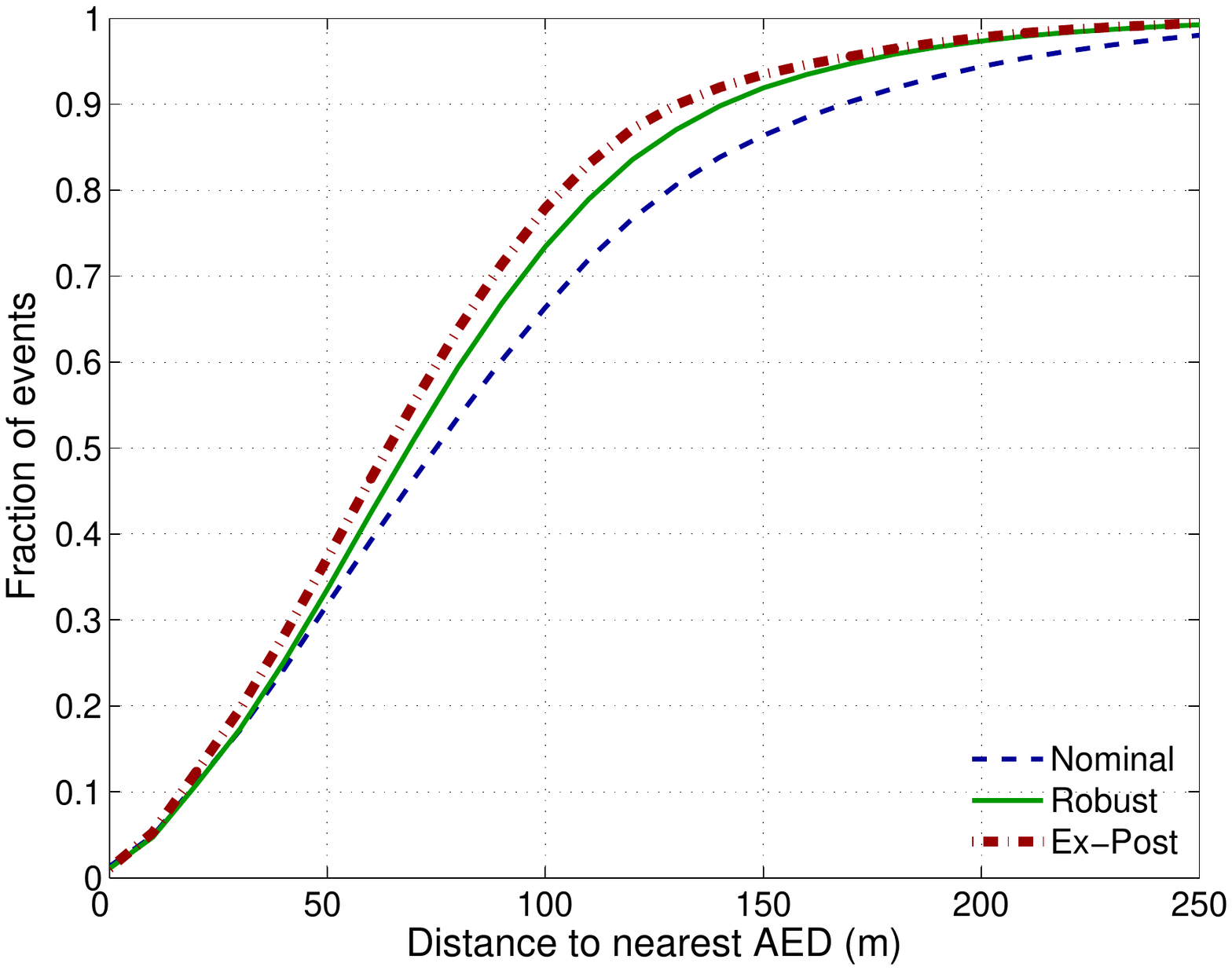}
\label{fig:subfigure1}}
\subfigure[$h = 150$.]{%
\includegraphics[scale=0.40, clip = true, trim = 15mm 65mm 20mm 70mm]{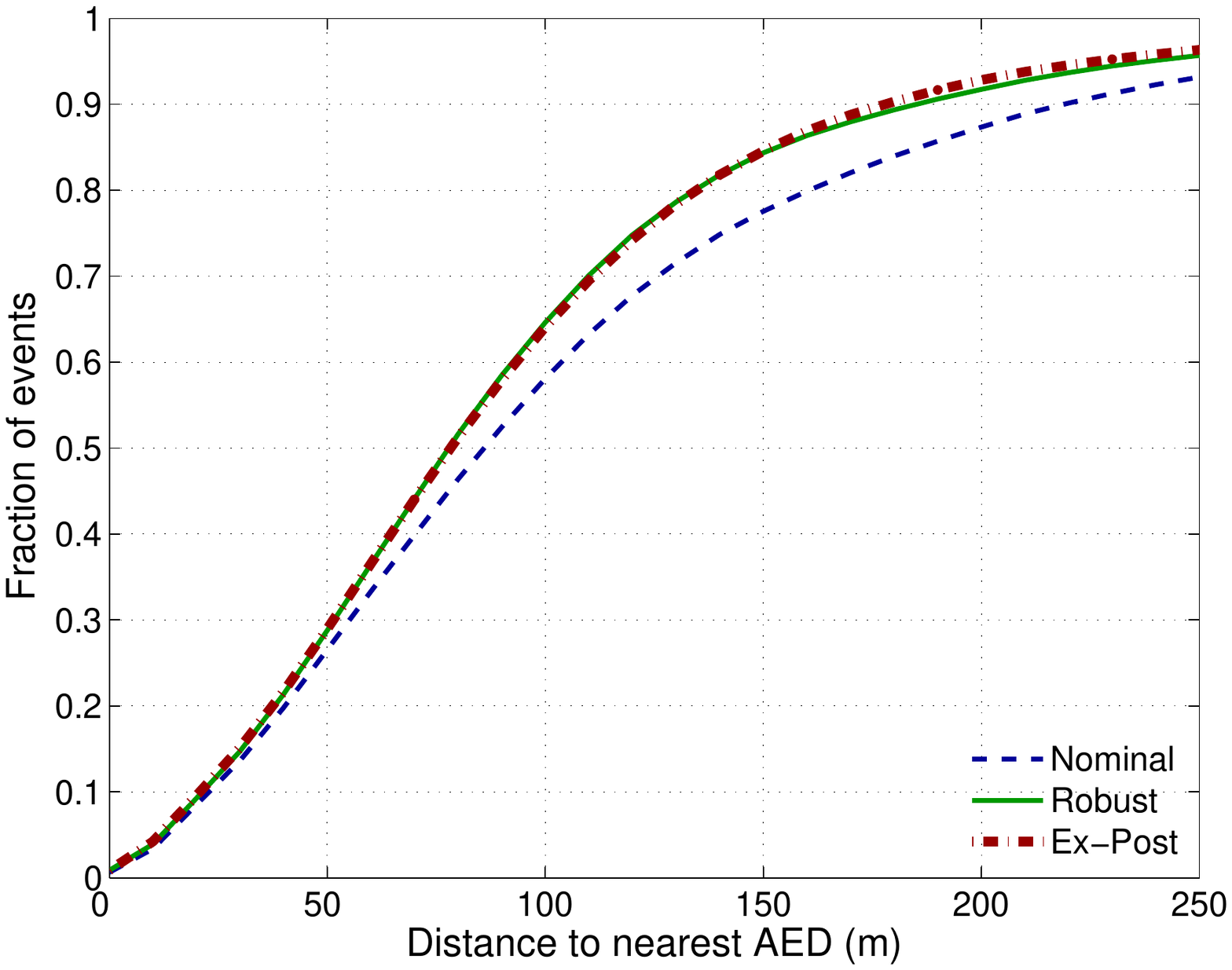}
\label{fig:subfigure2}}
\vspace{-6mm}
\caption{Cardiac arrest coverage for nominal, robust and ex-post solutions under four kernel bandwidths.}
\label{fig:coverage}
\end{figure}

An alternate but related measure of performance in public access defibrillation is cardiac arrest {\it coverage}, which is the fraction of cardiac arrests that occur within a given distance of at least one AED  \citep{chan2013}. Figure \ref{fig:coverage} shows the cardiac arrest coverage under the nominal, robust and ex-post solutions at all distances between 0 and 250 meters, for each of the four bandwidths. The curves shown in Figure \ref{fig:coverage} are obtained by computing the average coverage at each distance over all 50 trials. Note that in the case where $h = 10$, all three models attain similar coverage at all distances, which is consistent with the results shown in Figure \ref{fig:average}. Note that the robust solution dominates the nominal solution at all bandwidths, in the sense that it achieves a higher coverage level at all distances. Moreover, the coverage of the robust model is very close to the coverage obtained by the ex-post model, especially for high bandwidths. The improvement of the robust deployment over the nominal deployment appears most pronounced at higher distances ($\ge 100$m), which is expected behavior since the AED locations are optimized for the tail of the distance distribution.

\subsection{Alternate distance measures}\label{sec:alt}
We repeated a set of experiments using two alternate distance measures: the rectilinear distance, and actual walking distances based on city roads and walking paths. Table \ref{table:compare} shows the performance of each model when Euclidean, rectilinear and walking distances are used to compute the distances  (for both the $d_{ik}$ parameters and the performance metrics). When measured using actual walking distances, the robust and nominal deployments perform similarly on the mean distance metric. However, the robust model attains a lower VaR and CVaR, indicating improved performance by the robust deployment at the tail of the distance distribution.  Table \ref{table:comparedistancegap} provides the accompanying average and worst-case performance gaps of the nominal and robust solutions with respect to the ex-post solution. We found the results to be consistent with our main finding in Section \ref{sec:performancemetrics}: the AED deployments produced by the robust model outperform the nominal model on all four performance metrics, both on average and in the worst-case realization of cardiac arrest locations. On the CVaR performance metric (which we emphasize is the most relevant point of comparison since all three models optimize for CVaR), the robust model cuts the performance gap in half, compared to the nominal model. It should be noted that the walking distances from the Google Maps API are likely conservative estimates of the actual distance a lay responder would have to travel, since it constrains the lay responder to walk on recognized pedestrian paths. In practice, a lay responder may be able to take shortcuts that are unmarked in Google Maps (e.g., parking lots, buildings, unmarked crosswalks). We therefore expect the actual lay responder distance to be somewhere between the Euclidean and Google Maps distances.

\begin{table}
\caption{Average and worst-case performance of nominal, robust and ex-post solutions over 50 simulated validation sets using different distance measures ($\beta = 0.9$, $h = 100$).}
\centering
    \begin{tabular}{lrrrrrrrrrrrrr}
\toprule 
           & &\multicolumn{3}{c}{Mean} & \multicolumn{3}{c}{VaR} & \multicolumn{3}{c}{CVaR} & \multicolumn{3}{c}{Max} \\
            \cmidrule(l){3-5} \cmidrule(l){6-8} \cmidrule(l){9-11} \cmidrule(l){12-14}
          &  & \multicolumn{1}{c}{Nom.} & \multicolumn{1}{c}{Rob.} & \multicolumn{1}{c}{Ex.} & \multicolumn{1}{c}{Nom.} & \multicolumn{1}{c}{Rob.} & \multicolumn{1}{c}{Ex.}  & \multicolumn{1}{c}{Nom.} & \multicolumn{1}{c}{Rob.} &\multicolumn{1}{c}{Ex.}  &\multicolumn{1}{c}{Nom.}  & \multicolumn{1}{c}{Rob.}& \multicolumn{1}{c}{Ex.} \\ 
           \midrule
& Euclidean  &  95	& 86	& 81	& 173	&147 &133	&223	 &189	&179 	&303 	&260	&252 \\
Average & Rectilinear  &  121	& 107	& 101	& 216	&190	 &164	&280	 &248	&223	    &386	 &342	&316 \\
& Walking  &  134	& 133	& 119	& 233	&221	 &202	&281 &269	&256	    &367	 &365	&355\\
           \midrule
& Euclidean &   109	& 96	& 95	& 202	& 181	& 181	& 260	& 222	& 215	&408	&375	&375 \\
Worst-case \;\; & Rectilinear &    136	&120&	119&	263&	241&	215&	345&	311&	283&	553&	553&	498 \\
& Walking &  150	&148&	137&	256&	243&	236&	326&	315&	307&	558&	558&	572 \\
    \bottomrule
\label{table:compare}
    \end{tabular}%
\end{table}

\begin{table}[htbp]   
\caption{Average and worst-case performance gap of nominal and robust solutions with respect to ex-post solution under alternate distance measures.}
\centering
    \begin{tabular}{lrrrrrrrrr}
\toprule 
           & &\multicolumn{2}{c}{Mean} & \multicolumn{2}{c}{VaR} & \multicolumn{2}{c}{CVaR} & \multicolumn{2}{c}{Max} \\
            \cmidrule(l){3-4} \cmidrule(l){5-6} \cmidrule(l){7-8} \cmidrule(l){9-10}
         &   & \multicolumn{1}{c}{Nom.} & \multicolumn{1}{c}{Rob.} & \multicolumn{1}{c}{Nom.} & \multicolumn{1}{c}{Rob.}   & \multicolumn{1}{c}{Nom.} & \multicolumn{1}{c}{Rob.}   &\multicolumn{1}{c}{Nom.}  & \multicolumn{1}{c}{Rob.} \\
           \midrule
&Euclidean           &15\%	&6\%	&23\%	    &9\%	     &20\%	&5\%	    &17\%	&3\% \\
Average & Rectilinear  &17\%	&6\%	&24\%	    &13\%	 &21\%	&10\%	&18\%	&8\% \\
&Walking         &12\%	&11\%	&13\%	    &8\%	     &9\%	&5\%    	&7\%	&3\% \\
           \midrule
&Euclidean            &29\%	&17\%	&40\%	    &26\%	 &27\%	&16\%	&31\%	&24\% \\
Worst-case &Rectilinear &28\%	&18\%	&41\%	    &27\%	 &32\%	&20\%	&41\%	&23\% \\
&Walking          &21\%	&20\%	&31\%	    &24\%	 &20\%	&12\%	&24\%	&22\% \\
    \bottomrule
     \label{table:comparedistancegap}
        \end{tabular}
\end{table}
As an additional check to compare Euclidean and actual walking distances, we performed the following experiment. We randomly generated 1,000 pairs of locations within the downtown region of the City of Toronto, which includes the service area in our study. For each pair of locations, we computed the Euclidean distance (based on the UTM coordinates of the locations) and the walking distance (via the Google Maps Distance Matrix API). We then performed an ordinary least squares regression on the distance data, which yielded the following model: 
$$Dist_{\text{Walk}} = 1.341\times Dist_{\text{Euclid}} + 21.1.$$ 
The $R^2$ coefficient of this model was 0.904, suggesting a strong linear relationship between Euclidean and actual walking distances. Both coefficients were statistically significant ($p < 0.05$). This finding is consistent with previous studies that have also found very high correlation between Euclidean and road network distances \citep{bach1981problem,jones2010spatial,boscoe2012nationwide}. These results suggest that Euclidean distances serve as a reasonable (constant factor) approximation for actual walking distances.

\section{Discussion \& Policy Implications}
Interestingly, our results indicate that the robust deployment performed better than the nominal deployment under {\it typical} demand realizations (Figure \ref{fig:average}), in addition to the worst-case realizations (Figure \ref{fig:WC}). This might seem like a counterintuitive result, given that robust optimization models can perform poorly in typical realizations of the uncertainty due to their emphasis on optimizing a worst-case objective function. However, in the context of AED deployment (and facility location more generally), the observed performance gap between the robust and nominal models can likely be attributed to the nominal model overfitting to the set of historical cardiac arrest locations. In other words, the nominal model represents a naive sample average approximation in which AED placement is optimized only with respect to a small set of historical cardiac arrest locations. Our results suggest that if there is a nontrivial amount of uncertainty in the locations of future cardiac arrest, a sample average approximation based on limited historical data can lead to AED deployments that perform poorly when measured against out-of-sample cardiac arrests. On the other hand, our robust formulation produces an improved deployment by accounting for the possibility of cardiac arrests occuring in new locations not reflected in the historical data.

{\bf Survival outcomes.} Our results suggest that a robust approach to AED deployment can improve cardiac arrest response by improving bystander accessibility to AEDs, especially for those patients who collapse far away from the nearest defibrillator. While a rigorous analysis of survival outcomes is beyond the scope of this paper, the time sensitive nature of cardiac arrest treatment suggests that even modest improvements in AED accessibility (through shorter distances) can have a significant impact on survival outcomes. We also emphasize that all numerical results presented in Section \ref{sec:results} are with respect to the one-way distance to the nearest AED. However, since retrieval of an AED by a lay responder likely requires a round-trip, decreasing the distance between the cardiac arrest victim and the nearest AED by 30 meters implies that the trip made by the lay responder is shortened by 60 meters, which can have a material impact on the total time between collapse and treatment.

{\bf Comparison to best practices.} Historically, best practices for AED deployment (e.g, as prescribed by the American Heart Association) have been purely retrospective, and do not account for the uncertainty in future cardiac arrest locations, which our results suggest can help improve the accessibility of AEDs. Existing guidelines also emphasize the importance of identifying building {\it types} that are at risk of cardiac arrest, and placing AEDs in those locations. Our proposed approach is markedly different from the current AHA recommendations. We instead focus on characterizing the distributional uncertainty in cardiac arrest locations, and optimally deploying AEDs in a manner that does not depend on building type, as long as the AED can be accessed by a lay responder. We also note that placing AEDs in a small set of high risk buildings ignores the shape of the distance distribution, and may lead to AED deployments where a large number of cardiac arrests occur far away from any AED. In contrast, our approach considers the entirety of the (worst-case) distance distribution, and allows the modeler to directly optimize various aspects of the distribution via the CVaR objective function.

{\bf Crowdsourcing lay responders.} The existing focus on identifying high-risk buildings for AED placement is partially due to the current low rate of bystander intervention during cardiac arrests, and a prevailing belief that an AED must be located at the precise location of a cardiac arrest in order to be used. However, recent advances in mobile phone applications show promise in recruiting lay responders (who may be up to 500m away) to perform CPR on victims of cardiac arrest \citep{ringh2015mobile, brooks2016pulsepoint,brooks2014implementation}. This ``crowdsourcing" approach to cardiac arrest interventions can be extended to incorporate AEDs as well, by notifying lay responders of the locations of  both the cardiac arrest victim and a nearby AED. We posit that in a setting where lay responders are recruited to intervene during a cardiac arrest, it is even more critical that AEDs are tactically located throughout an urban area, so that they can be quickly located and transported to a victim of cardiac arrest by a lay responder. Thus, our approach to strategically deploying AEDs can work synergistically with other innovations in the chain of survival to improve the overall response to cardiac arrest.  

{\bf Risk measures in public access defibrillation.} Our model also highlights the potential role that risk measures like CVaR can play in public access defibrillation. While VaR-based response time targets are common in ambulance dispatching \citep{pons2002}, no similar guidelines exist in public access defibrillation. Given the potential impact on survival of protecting against large distances between cardiac arrest victims and nearby AEDs, it may be worthwhile to explicitly incorporate risk-based targets in the decision making process when determining AED locations. This paper formally introduces the notion of risk to the AED deployment problem, by allowing decision makers to emphasize the tail of the distance (or response time) distribution when planning AED locations.

\section{Conclusion}\label{sec:conclusion}
We develop a data-driven optimization approach for deploying AEDs in public spaces while accounting for uncertainty in the locations of future cardiac arrests. Our approach involves constructing a distributional uncertainty set that consists of several uncertainty regions, where the parameters of the uncertainty set are calibrated based on historical cardiac arrest data. To approximate the arrival of cardiac arrests in continuous space, we discretize the service area into an extremely large set of scenarios, and develop a row-and-column generation algorithm that exploits the structure of the uncertainty set and scales gracefully in the number of scenarios. As an auxilliary result, we show that our formulation subsumes a large class of facility location problems, and unifies the classical $p$-median and $p$-center problems, as well as their robust analogues.

Our numerical results suggest that hedging against cardiac arrest location uncertainty can lead to improved AED accessibility under both typical and worst-case realizations of the uncertainty. In particular, we found that our robust AED deployment outperformed a nominal (sample average approximation) deployment by between 9-20\% with respect to AED retrieval distances, depending on the performance metric and underlying demand distribution. Further, we found that in many cases, the robust deployment performed nearly as well as an ex-post model with perfect foresight, and in many instances improved the performance gap with respect to the ex-post model by 40-70\%. This finding suggests that our robust approach manages to avoid the performance loss due to uncertainty that is suffered by the sample average approach.

We highlight a few potential directions for future work. First, we have assumed throughout this paper that the set of uncertainty regions is given. It may be fruitful to investigate data-driven approaches to identifying an effective configuration of the uncertainty regions. Similarly, we note that the uncertainty set proposed in this paper can be generalized to permit the arrival probabilities to be uncertain as well. For example, we might allow each $\lambda_j$ to reside within a 95\% confidence interval that is also estimated from the data, which would retain the polyhedral structure of the uncertainty set. Our model can also be extended to incorporate additional considerations in public access defibrillation, such as lay responder behavior. However, we expect that our main finding in this paper -- that accounting for cardiac arrest location uncertainty improves AED accessibility -- would persist under such modifications to the model. Lastly, the main features of our modeling approach -- a distributional uncertainty set that induces sparse worst-case distributions coupled with an efficient row-and-column generation algorithm -- may find relevance in other applications with a similar problem structure as well.

\ACKNOWLEDGMENT{The authors thank the area editor Chung Piaw Teo, the associate editor, and three anonymous referees for providing valuable comments that significantly improved the paper.  The authors also thank the Natural Sciences and Engineering Research Council of Canada (NSERC) and the Ontario Ministry of Training, Colleges and Universities for financial support, and Dr. Laurie Morrison and the Resuscitation Outcomes Consortium for the cardiac arrest data.\\}

\noindent{\bf \large Authors}\\

\noindent{\bf Timothy Chan} is the Canada Research Chair in Novel Optimization and Analytics in Health, an Associate Professor in the department of Mechanical and Industrial Engineering and the Director of the Centre for Healthcare Engineering at the University of Toronto. His primary research interests are in robust and inverse optimization, and the application of optimization methods to problems in healthcare and sports.

{\bf Zuo-Jun (Max) Shen} is a Chancellor's Professor in the Department of Industrial Engineering and Operations Research and the Department of Civil and Environmental Engineering at UC Berkeley. He is also an honorary professor at Tsinghua University. He has been active in the following research areas: integrated supply chain design and management, design and analysis of optimization algorithms, energy system and transportation system planning and optimization.

{\bf Auyon Siddiq} is a Ph.D. candidate in the Department of Industrial Engineering and Operations Research at the University of California, Berkeley. His current research interests are in data-driven optimization, healthcare operations and incentive design.\\

\bibliographystyle{ormsv080} 
\bibliography{articles} 

\ECSwitch


\ECHead{Electronic companion for ``Robust Defibrillator Deployment Under Cardiac Arrest Location Uncertainty via Row-and-Column Generation".}

\section{Proofs.}
We first present two lemmas that are helpful in the proof of Theorem 1.

\begin{lemma}\label{pj}
Let ${\bf y}$ be fixed. Define $ \mathcal{U}^* = \underset{\mu \in\mathcal{U}}{\textnormal{argmax}} \; \textnormal{CVaR}_\mu[d(\xi,{\bf y}(\xi))].$ Then for any set of facility locations ${\bf y}$, there exists a worst-case distribution $\tilde{\mu} \in \mathcal{U}^*$ supported on $a_1,\ldots,a_{|\mathcal{R}|}$, where $a_r  = \underset{a \in \mathcal{A}'_r}{\textnormal{argmax}} \; d(a,{\bf y}(a))$, for all $r = 1,\ldots,|\mathcal{R}|$.
\end{lemma}
\noindent {\bf Proof.} 
For any $\bar{\mathcal{A}} \subset \mathcal{A}$, define
\begin{equation}\label{deflambda}
\lambda(\mu,\bar{\mathcal{A}}) = P_{\mu}(\xi \in \bar{\mathcal{A}}) = \int_{ \bar{\mathcal{A}}} f_{\mu}(\xi)d\xi, 
\end{equation}
where $f_{\mu}$ is the density function associated with distribution $\mu$. Pick a distribution $\mu^* \in \mathcal{U}^*$. Now construct $\tilde{\mu}$ to be the distribution with mass $\lambda(\mu^*,\mathcal{A}'_1),\ldots,\lambda(\mu^*,\mathcal{A}'_{|\mathcal{R}|})$ placed at the locations $a_1,\ldots,a_{|\mathcal{R}|}$. We now show that $\tilde{\mu}$ is also in $\mathcal{U}^*$. We do so by first showing that $\tilde{\mu}$ is in $\mathcal{U}$, and then by showing that the CVaR of the distance distribution under the constructed $\tilde{\mu}$ is at least as large as the CVaR under ${\mu^*}$. For the first step, let $\mathcal{R}_j$ index the subregions that comprise $\mathcal{A}_j$. Now note that for each $j \in \mathcal{J}$, we have $$P_{\tilde{\mu}}(\xi \in \mathcal{A}_j) = \sum_{r \in \mathcal{R}_j} P_{\tilde{\mu}}(\xi \in \mathcal{A}'_r) = \sum_{r \in \mathcal{R}_j} \lambda(\tilde{\mu},\mathcal{A}'_r) =  \sum_{r \in \mathcal{R}_j} P_{\mu^*}(\xi \in \mathcal{A}'_r) = P_{\mu^*}(\xi \in \mathcal{A}_j)= \lambda_j.$$ 
The first equality is due to the fact that $\bigcup_{r \in \mathcal{R}_j}\mathcal{A}'_r = \mathcal{A}_j$. The second equality follows from the definition of $\lambda(\mu,\bar{\mathcal{A}})$ in \eqref{deflambda}. The third equality follows by the construction of $\tilde{\mu}$.  The fourth equality is again due to $\bigcup_{r \in \mathcal{R}_j}\mathcal{A}'_r = \mathcal{A}_j$. The final equality holds since $\mu^* \in \mathcal{U}$.  Since from the above equations, $P_{\tilde{\mu}}(\xi \in \mathcal{A}_j) = \lambda_j$, we have $\tilde{\mu} \in \mathcal{U}.$ It remains to show that $\text{CVaR}_{\mu^*} \le \text{CVaR}_{\tilde{\mu}}$. Writing CVaR explicitly, we have
\begin{align*}
\text{CVaR}_{\mu^*}[d(\xi,{\bf y}(\xi))] &= \min_\alpha \; \alpha + \frac{1}{1-\beta} \int_{\mathcal{A}} [d(\xi,{\bf y}(\xi)) - \alpha]^+f_{\mu^*}(\xi)d\xi  \\
&= \min_\alpha \; \alpha + \frac{1}{1-\beta} \sum_{r \in \mathcal{R}} \int_{\mathcal{A}'_r} [d(\xi,{\bf y}(\xi)) - \alpha]^+f_{\mu^*}(\xi)d\xi \\
&\le \min_\alpha \; \alpha + \frac{1}{1-\beta} \sum_{r \in \mathcal{R}} \int_{\mathcal{A}'_r} [d(a_r,{\bf y}(a_r)) - \alpha]^+f_{\mu^*}(\xi)d\xi  \\
& = \min_\alpha \; \alpha + \frac{1}{1-\beta} \sum_{r \in \mathcal{R}}  [d(a_r,{\bf y}(a_r)) - \alpha]^+\int_{\mathcal{A}'_r}f_{\mu^*}(\xi)d\xi  \\
& = \min_\alpha \; \alpha + \frac{1}{1-\beta} \sum_{r \in \mathcal{R}}  [d(a_r,{\bf y}(a_r)) - \alpha]^+  \lambda(\tilde{\mu},\mathcal{A}'_r) \\
& = \text{CVaR}_{\tilde{\mu}}[d(\xi,{\bf y}(\xi))]. 
\end{align*}
The first line follows from the definition of CVaR. The second line is due to the fact that the subregions $\mathcal{A}'_r$, $r  \in \mathcal{R}$ form a partition of $\mathcal{A}$. The third line follows from the definition of $a_r$ as a point in $\mathcal{A}'_r$ that maximizes distance to the nearest AED. The fourth line follows since $[d(a_r,{\bf y}(a_r)) - \alpha]^+$ is constant with respect to the variable $\xi$. The fifth line follows from the definition of $\lambda(\tilde{\mu},\mathcal{A}'_r)$ as the total probability mass in $\mathcal{A}'_r$. The final line follows by construction of $\tilde{\mu}$ and the definition of CVaR. \qed

For the second lemma, define $z_C({\bf y}) = \max_{\mu \in \mathcal{U}} \text{CVaR}_\mu[d(\xi,{\bf y}(\xi))]$ to be the worst-case CVaR under the solution ${\bf y}$ in the continuous case. Similarly, let $z_D({\bf y})$ be the worst-case CVaR in the discrete case (i.e., the optimal value of \eqref{model:maxmincvar} given ${\bf y}$). The next lemma bounds the difference between the the worst-case continuous and discrete CVaR, for a given deployment ${\bf y}$. Note that Lemma 1 is used in the proof.
\begin{lemma}\label{ellsigma}
For any ${\bf y}$, $z_C({\bf y}) - z_D({\bf y}) \le \ell(\sigma)/(1-\beta)$.
\end{lemma}
 {\bf Proof.} Let $\mu^*$ be a worst-case continuous distribution with mass $\lambda'_1,\ldots,\lambda'_{|\mathcal{R}|}$ on the locations $p_1,\ldots,p_{|\mathcal{R}|}$ (cf Lemma \ref{pj}). Construct a discrete distribution by placing mass $\lambda'_1,\ldots,\lambda'_{|\mathcal{R}|}$ on the locations $\xi_1,\ldots,\xi_{|\mathcal{R}|}$, where $\xi_r  = \underset{\xi \in \Xi'_r}{\textnormal{argmin}} \; d(\xi,p_r)$ for all $r$. We can now write
\begin{align*}
z_C({\bf y}) - z_D({\bf y}) &= \left( \min_\alpha \; \alpha + \frac{1}{1-\beta} \sum_{r \in \mathcal{R}}  \lambda'_r [d(a_r,{\bf y}(a_r)) - \alpha]^+\right)  - z_D({\bf y}) \\
&\le \left( \min_\alpha \; \alpha + \frac{1}{1-\beta} \sum_{r \in \mathcal{R}}  \lambda'_r [d(\xi_r,{\bf y}(\xi_r)) + \ell(\sigma) - \alpha]^+ \right)  - z_D({\bf y}) \\
&\le \left( \min_\alpha \;  \alpha + \frac{1}{1-\beta} \sum_{r \in \mathcal{R}}  \lambda'_r [d(\xi_r,{\bf y}(\xi_r)) - \alpha]^+\right) + \frac{1}{1-\beta}\sum_{r \in \mathcal{R}}\lambda'_r\ell(\sigma)  - z_D({\bf y}) \\
& = z_D({\bf y}) + \frac{1}{1-\beta}\sum_{r \in \mathcal{R}}\lambda'_r\ell(\sigma)  - z_D({\bf y}) \\
& = \frac{\ell(\sigma)}{1 - \beta}
\end{align*}
The first line follows from the definition of $z_C({\bf y})$. The second line follows from the triangle inequality and since $d(\xi_r,a_r) \le \ell(\sigma)$ by definition of $\ell(\sigma)$. The third line follows from pulling $\ell(\sigma)$ out of the CVaR expression. The fourth line follows by definition of $z_D({\bf y})$. The fifth line follows from cancelling out $z_D({\bf y})$ terms and noting that $\sum_{r \in \mathcal{R}}\lambda'_r = 1$. \qed \\

\noindent {\bf Proof of Theorem 1.}\\
Note that $z_D({\bf y}) \le z_C({\bf y})$ for any ${\bf y}$, since $\Xi'_r \subset \mathcal{A}'_r$ for all $r \in \mathcal{R}$. Now define ${\bf y}_C  \in \underset{{\bf y}}{\text{argmin}} \; z_C({\bf y})$ and ${\bf y}_D  \in \underset{{\bf y}}{\text{argmin}} \; z_D({\bf y})$. Using Lemma \ref{ellsigma}, we can now write
$$z_D({\bf y}_D) \le z_D({\bf y}_C) \le z_C({\bf y}_C) = Z_C \le z_C({\bf y}_D) \le z_D({\bf y}_D) + \ell(\sigma)/(1-\beta) = Z_D + \ell(\sigma)/(1-\beta),$$
which implies $Z_C - Z_D \le \ell(\sigma)/(1-\beta)$. Since $Z_D \le Z_C$, we have
$$ \Delta = \frac{Z_C - Z_D}{Z_C} \le \frac{Z_C - Z_D}{Z_D} \le \frac{\ell(\sigma)}{Z_D(1-\beta)}.$$
It follows that if $\ell(\sigma) \le \varepsilon(1-\beta)Z_D$, then $\Delta \le \varepsilon$. \; \qed \\

\noindent {\bf Proof of Proposition \ref{mastersize}}.\\
Note that the number of variables in the master problem formulation \ref{masterprob} is $|\mathcal{I}|(|\mathcal{K}_+|+1)+ 1$, and the number of constraints is $(2|\mathcal{I}| + 1)|\mathcal{K}_+| + |\mathcal{I}| + |\mathcal{S}| + 2$. Hence the number of variables and constraints in iteration $s$ is $O(|\mathcal{I}||\mathcal{K}_+|)$. It remains to show that $|\mathcal{K}_+| \le (|\mathcal{J}|+1)s$ in the $s^{th}$ iteration of the row-and-column generation algorithm. Since the new variables and constraints that are added to \textsf{R-AED-MP} in iteration $s$ correspond directly to non-zero elements of the vector ${\bf p}_s$ obtained at a solution to \textsf{R-AED-SP} (cf Algorithm 1), it suffices to show that any ${\bf p}$ obtained at a solution to \textsf{R-AED-SP} contains at most $|\mathcal{J}|+1$ non-zero elements.  We first prove that any ${\bf u}$ obtained at an optimal solution to \textsf{R-AED-SP} is also an optimal solution to the following problem:
\begin{align}\label{subproblemu}
\underset{{\bf u}}{\text{maximize}} \;\; &\sum_{k \in \mathcal{K}} d_{ik} \bar{z}_{ik} u_k \notag \\
\text{subject to} \;\; &\sum_{k \in \mathcal{K}} u_k = 1, \\
& \sum_{k \in \mathcal{K}_j} u_k = \lambda_j, \quad j \in \mathcal{J},  \notag \\
& u_k \ge 0. \notag
\end{align}
By way of contradiction, suppose a solution ${\bf p},{\bf u}$ is optimal to \textsf{R-AED-SP} but ${\bf u}$ is suboptimal to \eqref{subproblemu}. Then there exists $k_1$ and $k_2$ such that $\sum_{i \in\mathcal{I}} d_{i{k_2}}\bar{z}_{i{k_2}} > \sum_{i \in\mathcal{I}} d_{i{k_1}}\bar{z}_{i{k_1}}$. Further, there must exist $\epsilon > 0$ and a vector $\tilde{{\bf u}} \in \mathcal{U}$ such that $\tilde{u}_{k_1} = u_{k_1} - \epsilon$, $\tilde{u}_{k_2} = u_{k_2} + \epsilon$, and $\tilde{u}_k = u_k$ for all $k \in \mathcal{K} \setminus \{ k_1,k_2 \}$. Let $\tilde{\bf p}$ be the solution obtained by solving \textsf{R-AED-SP} with $\tilde{\bf u}$ as a parameter. Now note in \textsf{R-AED-SP}, since $\sum_{i \in\mathcal{I}} d_{i{k_2}}\bar{z}_{i{k_2}} > \sum_{i \in\mathcal{I}} d_{i{k_1}}\bar{z}_{i{k_1}}$, the constraints $p_{k_2} \le u_{k_2}/(1-\beta)$ and $\tilde{p}_{k_2} \le \tilde{u}_{k_2}/(1-\beta)$ must be binding (otherwise we could slightly increase $p_{k_2}$ while slightly decreasing $p_{k_1}$ by the same amount to improve the objective function, with an identical argument applying to $\tilde{p}_{k_2}$). Since $p_{k_2} = u_{k_2}/(1-\beta)$ and $\tilde{p}_{k_2} = \tilde{u}_{k_2}/(1-\beta)$, by the construction of $\tilde{{\bf u}}$ it follows that $\tilde{p}_2 - p_2 \ge \epsilon / (1-\beta)$ and $p_1 - \tilde{p}_1 \le \epsilon / (1-\beta)$. Since $\sum_{i \in\mathcal{I}} d_{i{k_2}}\bar{z}_{i{k_2}} > \sum_{i \in\mathcal{I}} d_{i{k_1}}\bar{z}_{i{k_1}}$, the solution $\tilde{\bf p},\tilde{\bf u}$ improves the objective by $\left( \sum_{i \in\mathcal{I}} d_{i{k_2}}\bar{z}_{i{k_2}} - \sum_{i \in\mathcal{I}} d_{i{k_1}}\bar{z}_{i{k_1}} \right) \frac{\epsilon}{1-\beta}> 0$ over the solution ${\bf p},{\bf u}$, which contradicts the optimality of ${\bf p},{\bf u}$ with respect to \textsf{R-AED-SP}. Hence ${\bf u}$ must be optimal to \eqref{subproblemu} as well.

Since \eqref{subproblemu} is a standard form linear program with $|\mathcal{J}|+1$ equality constraints and $|\mathcal{K}|$ variables, at most $|\mathcal{J}|+1$ elements of a basic feasible solution to \eqref{subproblemu} can be non-zero \citep{bertsimas1997introduction}. Further, since by assumption there are no ties in the $d_{ik}$ parameters, every optimal solution to \eqref{subproblemu} must be a basic feasible solution. It follows that every ${\bf u}$ obtained at a solution of \textsf{R-AED-SP} has no more than $|\mathcal{J}|+1$ non-zero elements. Since $p_k > 0$ implies $u_k > 0$ for any feasible ${\bf p},{\bf u}$, it follows that any ${\bf p}$ obtained at a solution to \textsf{R-AED-SP} cannot have more than $|\mathcal{J}|+1$ non-zero elements. \qed  \\

\noindent {\bf Proof of Proposition \ref{prop:MLE}. }\\
Let the subregions $\mathcal{A}'_1,\ldots,\mathcal{A}'_{|\mathcal{R}|}$ be the partitioning of $\mathcal{A}$ implied by the uncertainty regions $\mathcal{A}_1,\ldots,\mathcal{A}_{|\mathcal{J}|}$. Let $\lambda'_r = P(\xi \in \mathcal{A}'_r)$ for all $r \in \mathcal{R}$ under the true distribution. Since the subregions are disjoint, the number of historical cardiac arrests that fall in each of the regions $\mathcal{A}'_1,\ldots,\mathcal{A}'_{|\mathcal{R}|}$ is simply the outcome of $n$ independent trials from a multinomial distribution with parameters $\lambda_1',\ldots,\lambda'_{|\mathcal{R}|}$. Now define 
$$ \hat{\lambda}'_r = \frac{1}{n} \sum_{c = 1}^{n} \mathbbm{1} \{a_c \in \mathcal{A}'_r \}, \quad r \in \mathcal{R},$$
is a maximum likelihood estimator of $\lambda'_1,\ldots,\lambda'_r$, and that $|\hat{\lambda}'_r - \lambda_r'| \longrightarrow 0$ with probability 1 as $n \longrightarrow \infty$ (e.g., Example 1.6.7. in \cite{bickel2015mathematical}). Now let $\mathcal{R}_j$ be the set of subregions that comprise region $\mathcal{A}_j$. We now have $\hat{\lambda}_j = \sum_{r \in \mathcal{R}_j}\hat{\lambda}'_r$, and thus $|\hat{\lambda}_j - \lambda_j|= | \sum_{r \in \mathcal{R}_j} (\hat{\lambda}_r' - \lambda_r') |\le  \sum_{r \in \mathcal{R}_j} | \hat{\lambda}_r' - \lambda_r'|$ for each $j \in \mathcal{J}$, which implies $|\hat{\lambda}_j - \lambda_j| \longrightarrow 0$ with probability 1 for each $j \in \mathcal{J}$ as $n \longrightarrow \infty$. To see that $\hat{\lambda}$ ensures that ${\bf U}$ is non-empty, construct a distribution ${\bf u} \in {\bf U}$ as follows. First, assign each of the cardiac arrest locations $a_1,\ldots,a_{n}$ to the nearest scenario location within the same subregion. Then for each $k \in \mathcal{K}$, set $u_k$ equal to the proportion of cardiac arrests that are assigned to location $k$. Thus ${\bf u}$ is non-negative and $\sum_{k \in \mathcal{K}} u_k = 1$ by construction. Lastly, for each $r \in \mathcal{R}$, $\sum_{k \in \mathcal{K}_r'} u_k$ is equal to the proportion of cardiac arrests within subregion $\mathcal{A}'_r$, which is also equal to $\hat{\lambda}'_r$ by definition. Thus $\sum_{k \in \mathcal{K}_j} u_k = \sum_{r \in \mathcal{R}_j} \sum_{k \in \mathcal{K}'_r} u_k = \sum_{r \in \mathcal{R}_j} \hat{\lambda}'_r = \hat{\lambda}_j$.  \qed \\

\section{Unifying the $p$-median and $p$-center problems}\label{sec:generality}
The $p$-median and $p$-center problems are two of the most well-studied location models, having served as the foundation for a significant portion of the existing facility location literature \citep{owen1998,snyder2006,melo2009facility}. In this section, we discuss how a non-stochastic interpretation of \ref{model:robustcvar2} specializes to robust variants of the $p$-median and $p$-center problems, depending on how $\beta$ is selected. A straightforward corollary of this result is that if the $\mathcal{A}_j$ are all singletons, then formulation \ref{model:robustcvar2} unifies the classical {\it p}-median and {\it p}-center problems (i.e. without any uncertainty). 

We begin by defining robust variants of the $p$-median and $p$-center problems, and then prove an equivalence between these problems and formulation \ref{model:robustcvar2}. Let $\mathcal{J}$ index a finite and fixed set of demand points. In the classical $p$-median problem, one wishes to site $P$ facilities so to minimize the (weighted)  total distance between all demand points and their nearest facilities.  Suppose now that the $j^{th}$ demand point is only known to reside within an uncertainty region $\mathcal{A}_j$. We formulate a generalization of the $p$-median problem that aims to minimize the {\it worst-case} total distance. Let each uncertainty region be discretized into a set of locations $\Xi_j \subset \mathcal{A}_j$, so that the demand point $j$ is only known to be at one of the locations in $\Xi_j$. Let $\mathcal{K}_j$ index the locations in $\Xi_j$. Let ${\bf Z}({\bf y})$ and ${\bf Y}$ retain their definitions from Section \ref{sec:model}. Let $\lambda_j$ be the usual notion of ``weight" placed on demand point $j$. Now define ${\bf X} = {\bf X}_1 \times \ldots \times {\bf X}_{|\mathcal{J}|}$, where
\begin{align}\label{Xset}
{\bf X}_j = \left\{ x_{jk} \ge 0, \; k \in \mathcal{K} ~\Bigg|~ \sum_{k \in \mathcal{K}_j} x_{jk} = 1\right\}, \quad j \in \mathcal{J}.
\end{align}
We can now write the {\it robust $p$-median problem} to minimize the worst-case total distance as
         \begin{align}
           \min_{\mathbf{y},\mathbf{z}} \max_{\mathbf{x}} \;\; & \sum_{i \in \mathcal{I}} \sum_{k \in \mathcal{K}_j} \sum_{j \in \mathcal{J} } \lambda_j d_{ik} z_{ik} x_{jk}   \label{robustmedian2} \\
           \mbox{\textnormal{subject to}}  \;\; &\mathbf{x} \in \mathbf{X}, \; \mathbf{z} \in \mathbf{Z}(\mathbf{y}), \; \mathbf{y} \in \mathbf{Y}. \notag
           \end{align}
In formulation \eqref{robustmedian2}, the set ${\bf X}$ can be interpreted as the decision space of an adversary which, after observing the facility locations, realizes each demand point within its uncertainty region at a location that maximizes the distance to the nearest facility. Using the same definitions as above, we can also define the {\it robust $p$-center problem}, which minimizes the worst-case {\it maximum} distance between a demand point and its nearest facility:
      \begin{align}
           \min_{\mathbf{y},\mathbf{z}} \max_{\mathbf{x}} \min_t \;\; & t   \label{robustcenter2} \notag\\
           \mbox{\textnormal{subject to}}  \;\; & t \ge \sum_{i \in \mathcal{I}} \sum_{k \in \mathcal{K}_j}  d_{ik}  z_{ik} x_{jk}, \quad   j \in \mathcal{J} ,  \\
           &\mathbf{x} \in \mathbf{X}, \; \mathbf{z} \in \mathbf{Z}(\mathbf{y}), \; \mathbf{y} \in \mathbf{Y}. \notag
           \end{align}
We now show that the robust $p$-median and $p$-center problems are unified by formulation \ref{model:robustcvar2}.
      \begin{proposition} \label{prop:pmedian}
Assume that each scenario $k$ belongs to exactly one of the sets $\mathcal{K}_1,\ldots,\mathcal{K}_{|\mathcal{J}|}$, and $\lambda_1 = \ldots = \lambda_{|\mathcal{J}|} = 1/n$. If $\beta= 0$, then formulation \ref{model:robustcvar2} is equivalent to the robust p-median problem \eqref{robustmedian2} with equally weighted demand points. If $\beta \ge (n-1)/n$, then formulation \ref{model:robustcvar2} is equivalent to the non-weighted robust $p$-center problem \eqref{robustcenter2}. \\
      \end{proposition}

\noindent {\bf Proof of Proposition \ref{prop:pmedian}.}\\
Since each scenario $k$ belongs to exactly one set $\mathcal{K}_j$, in the context of formulation \ref{model:robustcvar2} this implies that the uncertainty regions are disjoint. The set $\mathcal{U}$ is then separable in $j$, meaning the uncertainty in the discrete distribution ${\bf u}$ can be interpretted equivalently as uncertainty in $|\mathcal{J}|$ {\it conditional} distributions of ${\bf u}$ -- one for each uncertainty region. We may now interpret $x_j^k$ as the unknown probability that the next cardiac arrest arrives at location $k$, conditioned on the event $\{\xi \in \Xi_j\}$. Note that due to our choice of the uncertainty set $\mathcal{U}$, each vector ${\bf x}_j$ is only known to reside in the $|\mathcal{K}_j - 1|$-dimensional probability simplex, which is given precisely by the sets ${\bf X}_1,\ldots,{\bf X}_{|\mathcal{J}|}$. Also, since each scenario $k$ only belongs to one of the sets $\mathcal{K}_1,\ldots,\mathcal{K}_{|\mathcal{J}|}$, it follows that $u_k = \sum_{j \in \mathcal{J}} \lambda_j x_{jk}$, for each $k \in \mathcal{K}$. By substituting $\sum_{j \in \mathcal{J}} \lambda_j x_{jk}$ for $u_k$, we can write formulation \ref{model:robustcvar2} equivalently as
  \begin{align} 
  \min_{{\bf y},{\bf z}} \max_{{\bf x}} \min_{\alpha} \;\; & \alpha + \frac{1}{(1-\beta)} \sum_{j \in \mathcal{J}} \sum_{k \in \mathcal{K}_j} \lambda_jx_{jk} \max \left\{ \sum_{i \in \mathcal{I}} d_{ik} z_{ik} - \alpha, 0  \right\} \notag  \\
 \text{subject to} \;\; & {\bf x} \in {\bf X}, {\bf z} \in {\bf Z}({\bf y}), {\bf y} \in {\bf Y}, \label{twostage} \\
 & \alpha \ge 0. \notag
  \end{align}
Due to the separability of the scenario sets in $j$, we may take the $x_{jk}$ terms inside the $\max$ operator in the objective to obtain another reformulation of \ref{model:robustcvar2}:
  \begin{align} 
  \min_{{\bf y},{\bf z}} \max_{{\bf x}} \min_{\alpha} \;\; & \alpha + \frac{1}{(1-\beta)} \sum_{j \in \mathcal{J}} \lambda_j \max \left\{ \sum_{i \in \mathcal{I}} \sum_{k \in \mathcal{K}_j} d_{ik} z_{ik} x_{jk} - \alpha, 0  \right\} \notag  \\
 \text{subject to} \;\; & {\bf x} \in {\bf X}, {\bf z} \in {\bf Z}({\bf y}), {\bf y} \in {\bf Y}, \label{twostage2} \\
 & \alpha \ge 0. \notag
  \end{align}
We now show equivalence of \eqref{twostage2} and the robust $p$-median problem. Let $(\textbf{y}^*,\textbf{z}^*,\textbf{x}^*)$ be an optimal solution for formulation \eqref{twostage2}. Let $Z(\alpha)$ be the objective value of \eqref{twostage2} for the solution $(\textbf{y}^*,\textbf{z}^*,\textbf{x}^*, \alpha)$. Since the feasible sets of \eqref{twostage2} and \eqref{robustmedian2} are identical, all that remains is to show that their optimal values are equal. For conciseness define $g({\bf x}_j,{\bf z}_j) := \sum_{i \in \mathcal{I}} \sum_{k \in \mathcal{K}_j}  d_{ik} z_{ik} x_{jk}$. We first observe that for $\beta = 0$ and any $\alpha \le \underset{j \in \mathcal{J} }{\min} \; g({\bf x}^*_j,{\bf z}^*_j)$, the objective value of \eqref{twostage2},
\begin{align*}
 Z(\alpha) &= \alpha + \sum_{j \in \mathcal{J} } \lambda_j \left( g(\textbf{x}^*_j,\textbf{z}^*_j) - \alpha\right) \\
 &= \alpha + \sum_{j \in \mathcal{J} } \lambda_j g(\textbf{x}^*_j,\textbf{z}^*_j) - \frac{1}{n}(n\alpha) \\
 &=  \sum_{j \in \mathcal{J} } \lambda_j g(\textbf{x}^*_j,\textbf{z}^*_j) \\
 &=  \sum_{i \in \mathcal{I} } \lambda_j \sum_{j \in \mathcal{J}} \sum_{k \in \mathcal{K}_j} z_{ik}^* d_{ik} x_{jk}^*,
\end{align*}
is equivalent to the objective value of \eqref{robustmedian2}. Note that in the second equality above we use the fact that $\lambda_j = 1/n$ for all $j \in \mathcal{J}$. Thus to show equivalence in the optimal values of formulations \ref{model:robustcvar2} and \eqref{robustmedian2} when $\beta = 0$, it suffices to show that an optimal $\alpha^*$ satisfies $\alpha^* \le \underset{j \in \mathcal{J} }{\min} \; g(\textbf{x}^*_j,\textbf{z}^*_j)$. Suppose that $\alpha^* > \underset{j \in \mathcal{J} }{\min} \; g(\textbf{x}^*_j,\textbf{z}^*_j)$. Then, \begin{align*}
  Z(\alpha^*) &= \alpha^* + \sum_{j \in \mathcal{J} }\lambda_j \max \left\{ g(\textbf{x}^*_j,\textbf{z}^*_j) - \alpha^*,0\right\}  \\
  &> \alpha^* + \sum_{j \in \mathcal{J} } \lambda_j \left(  g(\textbf{x}^*_j,\textbf{z}^*_j) - \alpha^* \right)\\
  &=  \sum_{j \in \mathcal{J} } \lambda_j g(\textbf{x}^*_j,\textbf{z}^*_j) \\
   &=  \sum_{i \in \mathcal{I} } \sum_{j \in \mathcal{J}} \sum_{k \in \mathcal{K}_j} \lambda_j z_{ik}^{*} d_{ik} x_{jk}^{*},
 \end{align*}
which contradicts the optimality of $\alpha^*$. For the robust $p$-center problem, let $(\textbf{y}^*,\textbf{z}^*,\textbf{x}^*)$ again be an optimal solution to formulation \eqref{twostage2}.  Let $Z(\alpha)$ be the objective value of \eqref{twostage2} for the solution $(\textbf{y}^*,\textbf{z}^*,\textbf{x}^*,\alpha)$.  If $\alpha = \underset{j \in \mathcal{J} }{\max} \;\{ g(\textbf{x}^*_j,\textbf{z}^*_j) \}$, then the objective of formulation \eqref{twostage2},
\begin{align*}
 Z(\alpha) &= \alpha + \frac{1}{(1-\beta)} \sum_{j \in \mathcal{J} } \lambda_j \max \; \left\{ g(\textbf{x}^*_j,\textbf{z}^*_j) - \alpha, 0\right\}  \\
 &= \underset{j \in \mathcal{J} }{\max} \;\{ g(\textbf{x}^*_j,\textbf{z}^*_j) \} + \frac{1}{(1-\beta)} \sum_{j \in \mathcal{J} }  \lambda_j \max \; \left\{ g(\textbf{x}^*_j,\textbf{z}^*_j) - \underset{j \in \mathcal{J} }{\max} \;\{ g(\textbf{x}^*_j,\textbf{z}^*_j) \},  0 \right\}  \\
 &= \underset{j \in \mathcal{J} }{\max} \;\{ g(\textbf{x}^*_j,\textbf{z}^*_j) \} \\
 & = \underset{j \in \mathcal{J}}{\text{max}} \left\{ \sum_{i \in \mathcal{I} } \sum_{k \in \mathcal{K}_j} z_{ik}^{*} d_{ik} x_{jk}^* \right\},
\end{align*}
 is equivalent to the objective of \eqref{robustcenter2}. We now wish to show that an optimal $\alpha^*$ satisfies $\alpha^* = \underset{j \in \mathcal{J} }{\max} \;\{ g(\textbf{x}^*_j,\textbf{z}^*_j) \}$ for $\beta > (n-1)/n$.
We do so by showing that neither  $\alpha^* > \underset{j \in \mathcal{J} }{\max} \;\{ g(\textbf{x}^*_j,\textbf{z}^*_j) \}$ nor  $\alpha^* < \underset{j \in \mathcal{J} }{\max} \;\{ g(\textbf{x}^*_j,\textbf{z}^*_j) \}$ can hold. First, suppose $\alpha^* > \underset{j \in \mathcal{J} }{\max} \;\{ g(\textbf{x}^*_j,\textbf{z}^*_j) \}$. Then,
\begin{align*}
 Z(\alpha^*) &= \alpha^* + \frac{1}{(1-\beta)} \sum_{j \in \mathcal{J} } \lambda_j \max \;  \left\{ g(\textbf{x}^*_j,\textbf{z}^*_j) - \alpha^*, 0\right\}\\
 & = \alpha^* \\
&> \underset{j \in \mathcal{J} }{\max} \;\{ g(\textbf{x}^*_j,\textbf{z}^*_j) \},
\end{align*}
which contradicts the optimality of $\alpha^*$. Now suppose $\alpha^* < \underset{j \in \mathcal{J} }{\max} \;\{ g(\textbf{x}^*_j,\textbf{z}^*_j) \}$, and let $\delta =  \underset{j \in \mathcal{J} }{\max} \;\{ g(\textbf{x}^*_j,\textbf{z}^*_j) \} - \alpha^*.$ It follows that $\sum_{j \in \mathcal{J}} \max\{g({\bf x}^*_j,{\bf z}_j^*) - \alpha^* , 0\} \ge \delta$. Since $\lambda_j = 1/n$ for all $j \in \mathcal{J}$, we can write
\begin{align*}
 Z(\alpha^*) &= \alpha^* + \frac{1}{(1-\beta)n} \sum_{j \in \mathcal{J} }  \max \;  \left\{ g(\textbf{x}^*_j,\textbf{z}^*_j) - \alpha^*, 0\right\}\\
 & \ge \alpha^* + \frac{1}{(1-\beta)n} \delta \\
 &= \alpha^* + \delta + \frac{1}{(1-\beta)n} \delta  - \delta\\
 &=  \underset{j \in \mathcal{J} }{\max} \;\{ g(\textbf{x}^*_j,\textbf{z}^*_j) \} + \left(\frac{1}{(1-\beta)n} - 1 \right) \delta \\
 & >  \underset{j \in \mathcal{J} }{\max} \;\{ g(\textbf{x}^*_j,\textbf{z}^*_j) \},
\end{align*}
where the final strict inequality follows from $\beta > (n-1)/n$. Again, this contradicts the optimality of $\alpha^*$. \; \; \qed 

Based on Proposition \ref{prop:pmedian}, we can interpret formulation \ref{model:robustcvar2} as subsuming a spectrum of location models with the robust $p$-median problem at one end ($\beta = 0$) and the robust $p$-center problem at the other ($\beta$ close to 1). While the specialization of CVaR to the mean and maximum of a general discrete distribution is well known, this explicit unification of the $p$-median and $p$-center problems via a CVaR objective does not appear to have been discussed previously in the facility location literature. We note also that the assumption in Proposition \ref{prop:pmedian} that each $k$ belong to exactly one $\mathcal{K}_j$ is needed to ensure the equivalence between the main formulation \ref{model:robustcvar2} and the intermediate model \eqref{twostage} used in the proof. If this assumption is relaxed, then formulation \eqref{twostage} still serves as a unifying model for the $p$-median and $p$-center problems, although the equivalence with \ref{model:robustcvar2} no longer holds. The {\it ordered median problem} \citep{kalcsics2002,nickel1999} also unifies median and center problems, but requires an ordering of the demands according to their distance to the sited facilities so that a set of appropriate parameters can be identified. By contrast, our approach involves the adjustment of a single parameter and requires no such ordering.

To the best of our knowledge, \eqref{robustmedian2} and \eqref{robustcenter2} are the first robust formulations of the general $p$-median and $p$-center problem with non-stochastic uncertainty in the demand location (i.e., in the spirit of robust optimization rather than scenario-based stochastic optimization). Previous work in robust facility location has considered only one facility \citep{averbakh2000,averbakh2003} or uncertainty in demand level rather than location \citep{baron2011}. As discussed in Section 3.2, the formulation in \citep{baron2011} is well equipped to handle uncertainty in the level of demand at each of a set of known demand points. The robust $p$-median and $p$-center formulations presented above may be more appropriate in applications where there is continuous uncertainty in future demand locations (e.g., emergency response or peer-to-peer ridesharing applications). \\

\FloatBarrier

\section{Visualization of AED deployments}\label{additional}
\begin{figure}
\vspace{10mm}
\centering
\subfigure[Nominal.]{%
\includegraphics[scale=0.45, clip = true, trim = 25mm 15mm 10mm 65mm]{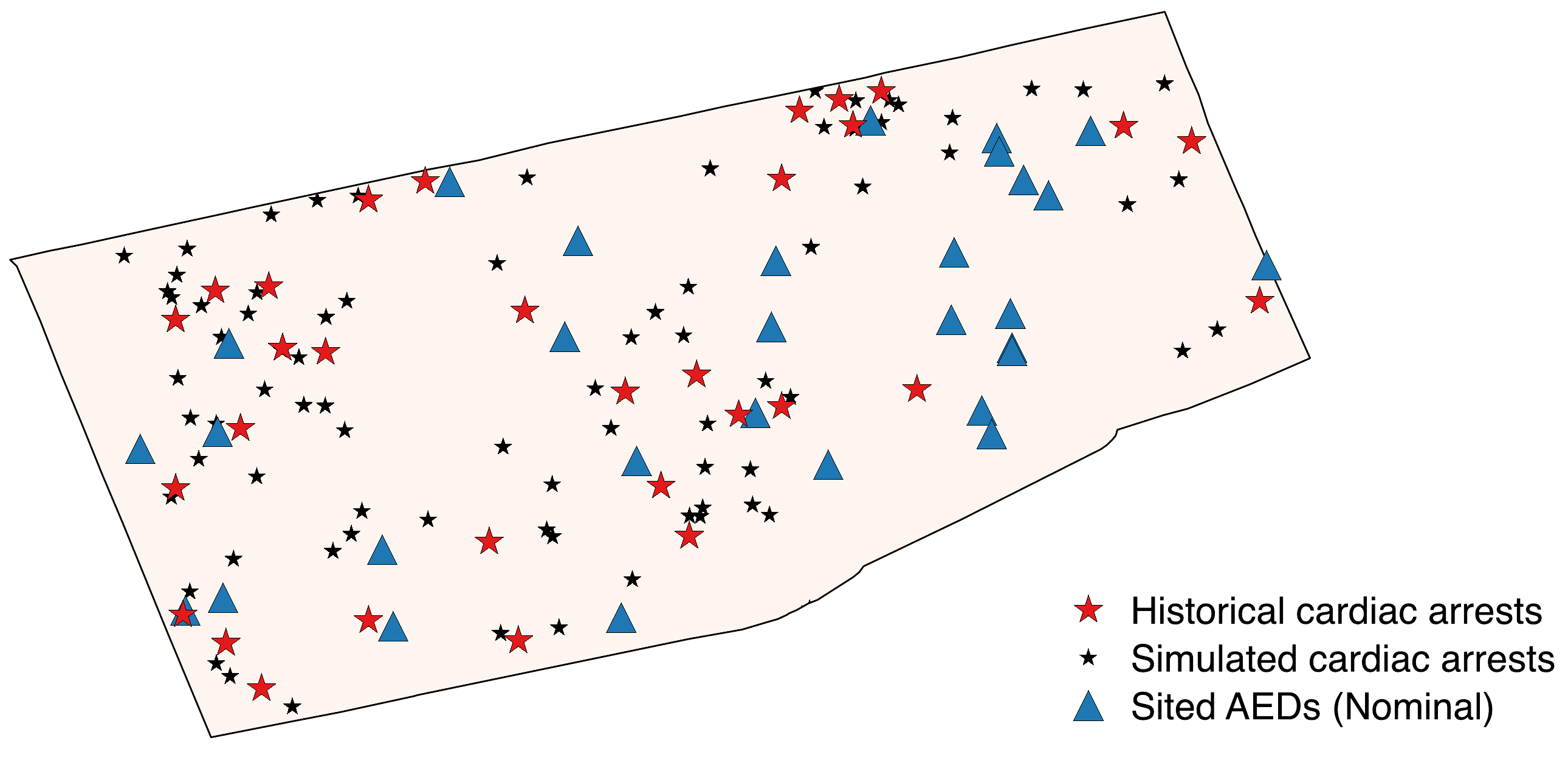}
\label{fig:subfigure1}}
\subfigure[Robust.]{%
\includegraphics[scale=0.45, clip = true, trim = 25mm 15mm 10mm 65mm]{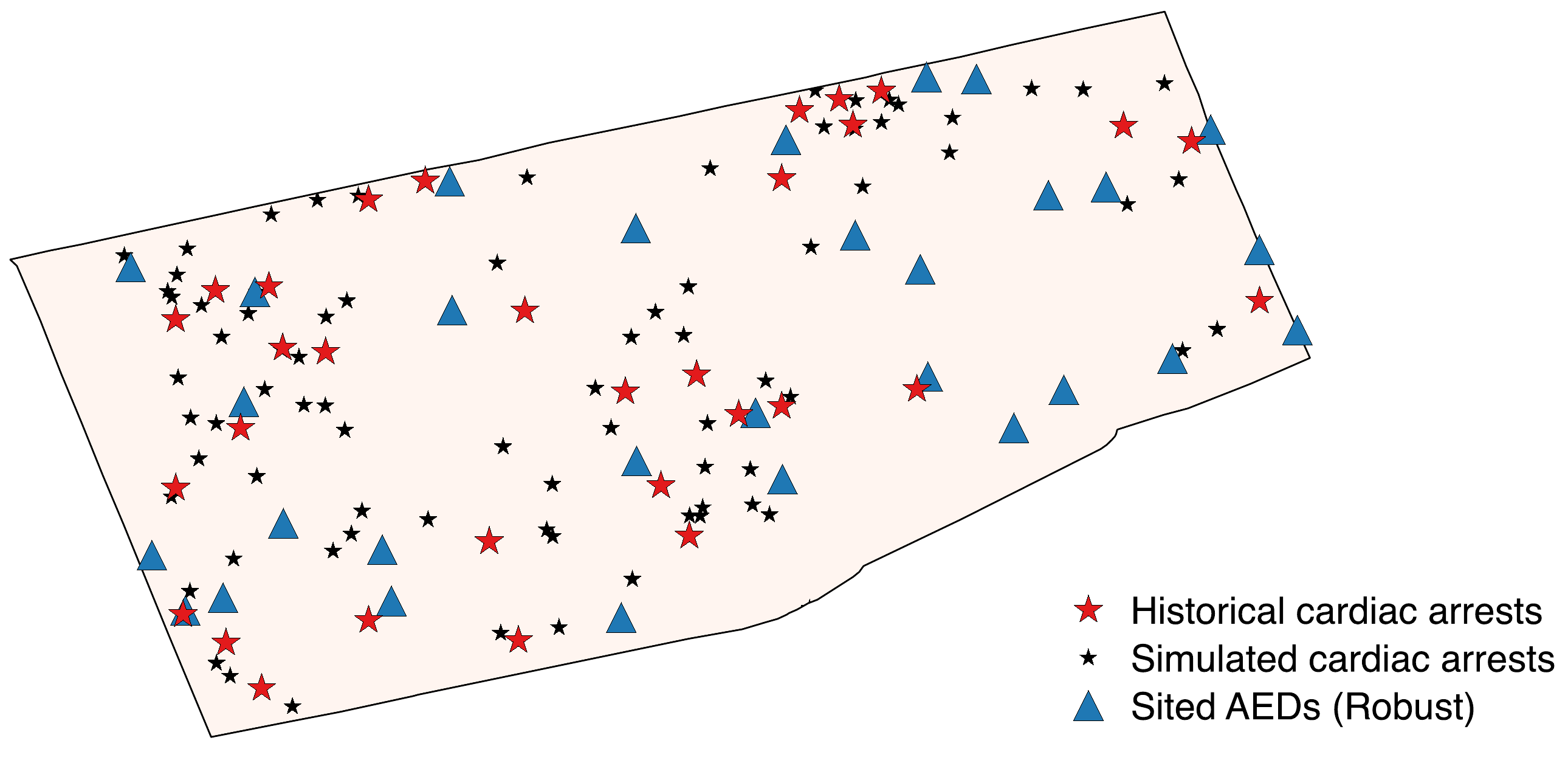}
\label{fig:subfigure2}}
\subfigure[Ex-post.]{%
\includegraphics[scale=0.45, clip = true, trim = 25mm 15mm 10mm 65mm]{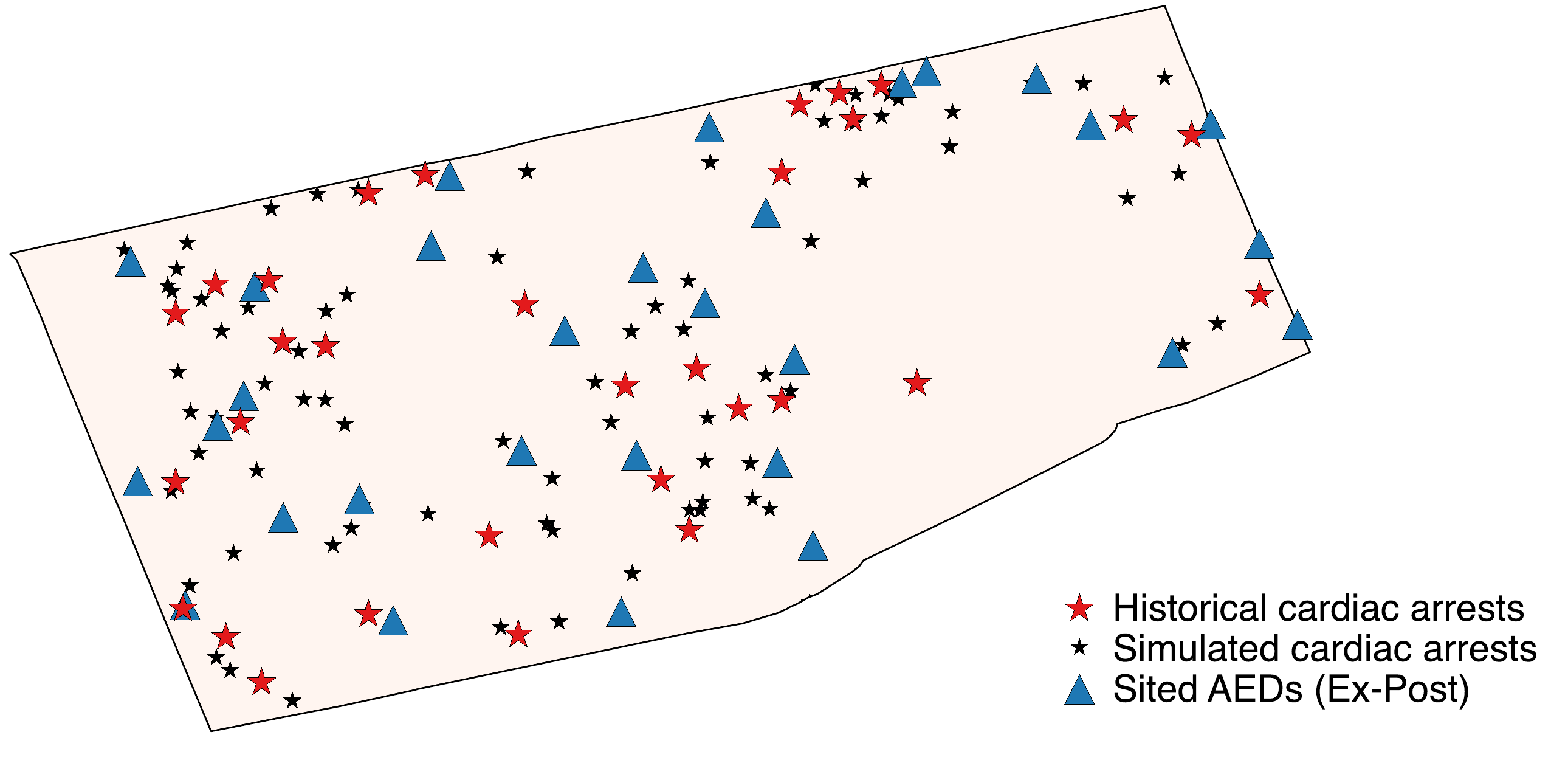}
\label{fig:subfigure1}}
\vspace{2mm}
\caption{Example of nominal, robust and ex-post AED deployments $(P = 30$, $h = 100$, $\beta = 0.9)$.}
\label{fig:AEDvisualize}
\end{figure}

Figure \ref{fig:AEDvisualize} depicts an instance of the AED deployments produced by the nominal (SAA), robust and ex-post models, along with a set of simulated cardiac arrests (corresponding to the ex-post deployment) and the 43 historical cardiac arrest locations. Note that in this instance, the nominal model produces small ``clusters" of AEDs, whereas the robust and ex-post models appear to position AEDs more efficiently with respect to the simulated cardiac arrests. This behavior is expected, since the nominal model only optimizes the tail of the empirical distance distribution, whereas the robust model accounts for potential deviation of the simulated cardiac arrests from the historical locations. Note also that the ex-post model achieves the most effective deployment of AEDs with respect to the simulated cardiac arrests, since it optimizes directly with respect to the simulated data.

\FloatBarrier 

\section{Row-and-column generation example: Total scenarios and optimality gap}\label{totalscenarios}
The plot below depicts an example of the total number of scenarios generated in the master problem \textsf{R-AED-MP} and the associated optimality gap at each iteration of the row-and-column generation algorithm. The plot corresponds to the instance given in the second-last row of Table 1 ($\mathcal{I} = 100$, $\mathcal{J} = 20$, $P = 20$, $\mathcal{K} = 24,743$).

\FloatBarrier   
\begin{figure}
\begin{center}
\includegraphics[scale=0.7, clip = true, trim = 10mm 65mm 10mm 65mm]{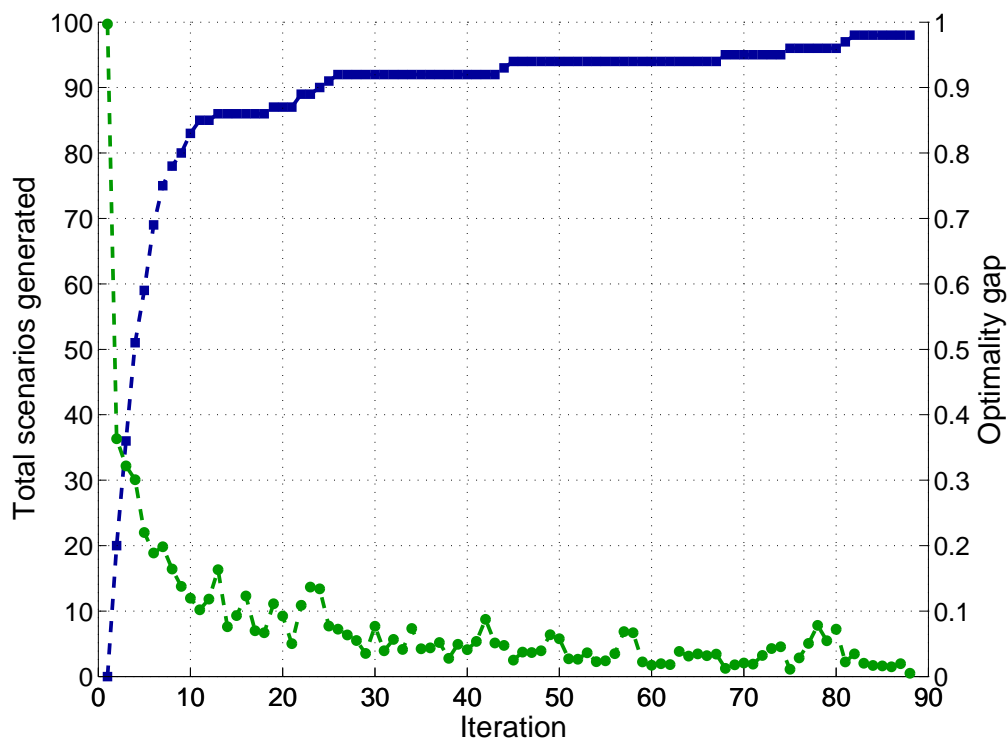}
\end{center}
\caption{Example of total number of scenarios generated (blue squares) and optimality gap (green circles) for each iteration of the row-and-column generation algorithm ($|\mathcal{I}| = 100$, $|\mathcal{J}| = 20$, $P = 20$, $|\mathcal{K}| = 24,743$).}
\vspace{50mm}
\end{figure}

\newpage 

\section{Estimates of uncertainty region arrival probabilities}\label{estimateprobabilities}
\begin{table}[htbp]
  \centering
  \caption{Estimates of probabilities $\lambda_1,\ldots,\lambda_{|\mathcal{J}|}$ and approximate 95\% confidence intervals for 15 uncertainty regions.}
  \vspace{2mm}
    \begin{tabular}{c|cc}
    \toprule
   $j$ & $\hat{\lambda}_j$ & 95\% C.I.   \\
    \midrule
 1 &  0.12 & (0.05, 0.24) \\
 2 &  0.16 & (0.07, 0.30) \\
 3 &  0.09 & (0.03, 0.22) \\
 4 &  0.12 & (0.04, 0.25) \\
 5 &  0.07 & (0.01, 0.19) \\
 6 &  0.05 & (0.01, 0.15) \\
 7 &  0.05 & (0.01, 0.15) \\
 8 &  0.02 & (0.00, 0.12) \\
 9 &  0.02 & (0.00, 0.12) \\
 10 &  0.09 & (0.03, 0.22) \\
 11 &  0.02 & (0.00, 0.12) \\
 12 &  0.02 & (0.00, 0.12) \\
 13 &  0.07 & (0.01, 0.18) \\
 14 &  0.05 & (0.01, 0.15) \\
 15  &  0.05 & (0.01, 0.15) \\
    \bottomrule
    \end{tabular}
  \label{tab:addlabel}
\end{table}

\end{document}